\pdfoutput=1

\documentclass{amsart}

\usepackage[dvipsnames]{xcolor}
\usepackage[utf8]{inputenc}
\usepackage[T1]{fontenc}

\usepackage{amsmath}
\usepackage{amsthm}
\usepackage{amsfonts}
\usepackage{amssymb}
\usepackage{mathtools}
\usepackage{braket}
\usepackage{mathrsfs}
\usepackage{mleftright}
\usepackage{thm-restate}

\makeatletter
\@ifundefined{newcounteralias}{}{%
  \renewcommand\thmt@autorefsetup{%
    \@xa\def
    \csname\thmt@envname autorefname\@xa\endcsname
    \@xa{\thmt@thmname}%
  }%
}
\makeatother

\usepackage{ifthen}
\usepackage{here}
\usepackage{array}
\usepackage{comment}
\usepackage{verbatim}
\usepackage{todonotes}
\usepackage{enumitem}

\newlist{steps}{enumerate}{1}
\setlist[steps,1]{label=Step \arabic*:}

\usepackage{graphicx}
\usepackage{tikz}
\usetikzlibrary{patterns,decorations.pathreplacing,calligraphy}
\usepackage[all]{xy}
\usepackage{amscd}
\usepackage{tikz-cd}
\usepackage[abs]{overpic}
\usepackage{pinlabel}
\usepackage{rotating}

\usepackage{geometry}
\usepackage[
  pagebackref,
  hypertexnames=false,
  colorlinks=true,
  linkcolor=NavyBlue,
  citecolor=NavyBlue,
  urlcolor=NavyBlue
]{hyperref}

\usepackage[alphabetic,backrefs,msc-links]{amsrefs}
\usepackage{cleveref}

\makeatletter
\DeclareRobustCommand\widecheck[1]{{\mathpalette\@widecheck{#1}}}
\def\@widecheck#1#2{%
   \setbox\z@\hbox{\m@th$#1#2$}%
   \setbox\tw@\hbox{\m@th$#1%
      {%
          \vrule\@width\z@\@height\ht\z@
          \vrule\@height\z@\@width\wd\z@}$}%
   \dp\tw@-\ht\z@
   \@tempdima\ht\z@
   \advance\@tempdima2\ht\tw@
   \divide\@tempdima\thr@@
   \setbox\tw@\hbox{%
      \raise\@tempdima\hbox{%
         \scalebox{1}[-1]{\lower\@tempdima\box\tw@}}}%
   {\ooalign{\box\tw@ \cr \box\z@}}}
\makeatother

\theoremstyle{plain}

\newtheorem{thm}{Theorem}[section]
\crefname{thm}{Theorem}{Theorems}
\Crefname{thm}{Theorem}{Theorems}

\newtheorem{prop}[thm]{Proposition}
\crefname{prop}{Proposition}{Propositions}
\Crefname{prop}{Proposition}{Propositions}

\newtheorem{lem}[thm]{Lemma}
\crefname{lem}{Lemma}{Lemmas}
\Crefname{lem}{Lemma}{Lemmas}

\newtheorem{cor}[thm]{Corollary}
\crefname{cor}{Corollary}{Corollaries}
\Crefname{cor}{Corollary}{Corollaries}

\newtheorem{rem}[thm]{Remark}
\crefname{rem}{Remark}{Remarks}
\Crefname{rem}{Remark}{Remarks}

\crefname{claim}{Claim}{Claims}
\Crefname{claim}{Claim}{Claims}

\crefname{property}{Property}{Properties}
\Crefname{property}{Property}{Properties}

\crefname{problem}{Problem}{Problems}
\Crefname{problem}{Problem}{Problems}

\newtheorem{question}[thm]{Question}
\crefname{question}{Question}{Questions}
\Crefname{question}{Question}{Questions}

\crefname{conjecture}{Conjecture}{Conjectures}
\Crefname{conjecture}{Conjecture}{Conjectures}

\theoremstyle{definition}

\newtheorem{defn}[thm]{Definition}
\crefname{defn}{Definition}{Definitions}
\Crefname{defn}{Definition}{Definitions}

\crefname{notation}{Notation}{Notations}
\Crefname{notation}{Notation}{Notations}

\crefname{convention}{Convention}{Conventions}
\Crefname{convention}{Convention}{Conventions}

\crefname{cond}{Condition}{Conditions}
\Crefname{cond}{Condition}{Conditions}

\crefname{assum}{Assumption}{Assumptions}
\Crefname{assum}{Assumption}{Assumptions}

\theoremstyle{remark}

\crefname{rem}{Remark}{Remarks}
\Crefname{rem}{Remark}{Remarks}

\crefname{section}{Section}{Sections}
\Crefname{section}{Section}{Sections}

\crefname{subsection}{Subsection}{Subsections}
\Crefname{subsection}{Subsection}{Subsections}

\crefname{figure}{Figure}{Figures}
\Crefname{figure}{Figure}{Figures}

\newcommand{\Z}{\mathbb{Z}}

\newcommand{\C}{\mathbb{C}}
\newcommand{\s}{\mathfrak{s}}

\newcommand{\R}{\mathbb R}

\def\det{\mathrm{det}}

\def\dim{\mathrm{dim}}

\def\id{\mathrm{Id}}

\newcommand{\mbar}[1]{%
{\ooalign{\hfil#1\hfil\crcr\raise.167ex\hbox{--}}}}

\def\cR{\mathcal{R}}

\def\wt{\widetilde}
\def\H{\mathbb{H}}
\def\Pin{\mathrm{Pin}}

\def\HMt{%
   \setbox0=\hbox{$\widehat{\mathit{HM}}$}
   \setbox1=\hbox{$\mathit{HM}$}
   \dimen0=1.1\ht0
   \advance\dimen0 by 1.17\ht1
   \smash{\mskip2mu\raise\dimen0\rlap{%
      \begin{turn}{180}
          {$\widehat{\phantom{\mathit{HM}}}$}
      \end{turn}}\mskip-2mu
      \mathit{HM}
   }{\vphantom{\widehat{\mathit{HM}}}}{}}

\title{Exotic diffeomorphisms and psc metrics on the 4-sphere}

\author{Sungkyung Kang}
\address{Department of Pure Mathematics and Mathematical Statistics, University of Cambridge, United Kingdom}
\email{sungkyung.kang@dpmms.cam.ac.uk}

\author{JungHwan Park}
\address{Department of Mathematical Sciences and RIM, Seoul National University, Republic of Korea}
\email{jungpark0817@snu.ac.kr}

\author{Masaki Taniguchi}
\address{Department of Mathematics, Kyoto University, Japan}
\email{taniguchi.masaki.7m@kyoto-u.ac.jp}

\begin{document}

\begin{abstract}
We show that all homotopy $4$-spheres and homotopy $K3$ surfaces admit exotic diffeomorphisms detected by family Seiberg--Witten theory. Using its proof, we also show that the space of positive scalar curvature metrics on the standard $4$-sphere is disconnected, each of its connected components has nontrivial fundamental group, and its component that contains the standard round metric has nontrivial homotopy groups in infinitely many degrees.
\end{abstract}

\maketitle
\section{Introduction}

An orientation-preserving self-diffeomorphism of a smooth $4$-manifold is called \emph{exotic} if it is topologically isotopic but not smoothly isotopic to the identity. Ruberman constructed the first examples on closed simply connected $4$-manifolds using one-parameter Donaldson theory \cite{Ruberman:1998}. For subsequent work on closed $4$-manifolds, see, for example, \cite{Ruberman:1999,ruberman2002positive,Auckly-Kim-Melvin-Ruberman:2015,Baraglia-Konno:2020,Kronheimer-Mrowka:2020-1,Lin:2023-1,Konno:2024-moduli,Qiu:2024-exotic,Auckly-Ruberman:2025,KLMM:2025-Lefschetz,Konno-Lin:2026,Baraglia:2026-mapping,Konno-Mallick-Taniguchi:2026,Baraglia-Tomlin:2026,baraglia2026exotic,Baraglia-Konno:2026,KLMM:2026-aspherical}. For compact $4$-manifolds with boundary, see \cite{KLMM:2024-Milnor,Miyazawa:2024-boundary,KLMM:2024-monodromy,Iida-Konno-Mukherjee-Taniguchi:2025,Konno-Mallick-Taniguchi:2025-absolute,kang2025exotic,KPT26,Krushkal-Mukherjee-Powell-Warren:2026,Konno-Mallick-Taniguchi:2026}.

For the standard $4$-sphere, Gabai, Gay, and Hartman announced the existence of an exotic diffeomorphism in \cite[Theorem~0.11]{gabai2025pseudo}, deferring the proof to a sequel. Their approach uses pseudo-isotopy theory and invariants of loops of embedded $2$-spheres. We use families Seiberg--Witten theory to detect an exotic diffeomorphism of every smooth $\mathbb F_2$-homology $4$-sphere, including the standard $S^4$.

\begin{restatable}{thm}{MainSphereTheorem}\label{thm: main}
Let $X$ be a smooth $\mathbb F_2$-homology $4$-sphere. Then the following hold:
\begin{itemize}
\item $X$ admits an exotic diffeomorphism.
\item The space $\mathcal R_+(X)$ of positive scalar curvature metrics on $X$ is either empty or disconnected.
\end{itemize}
\end{restatable}

Since an exotic diffeomorphism of $S^4$ can be smoothly isotoped
to be the identity near a point, restricting it to the complement
of a small standard $4$-ball yields an exotic diffeomorphism of
$D^4$ relative to the boundary.
Thus \Cref{thm: main} also has a consequence for spaces of
symplectic forms, giving a negative answer to
\cite[Question~3.1]{Wang:2026}.

\begin{cor}\label{cor: symplectic forms on D4}
Let $\omega_{\mathrm{std}}$ be the standard symplectic form on
$D^4$, and fix a collar neighborhood $U$ of $\partial D^4$.
The space of symplectic forms on $D^4$ that agree with
$\omega_{\mathrm{std}}$ on $U$, equipped with the $C^\infty$
topology, is disconnected.
\end{cor}
This follows from Gromov's theorem that the group
$\mathrm{Symp}_c(D^4,\omega_{\mathrm{std}})$ of symplectomorphisms
supported in the interior of the standard symplectic $4$-ball
is weakly contractible \cite{Gromov:1985}, together with \cite[Lemma~2.3]{Wang:2026}. 
See also
\cite{Kronheimer:families,Watanabe:2018,
Lin-Wu-Xie-Zhang:2025,Wang:2026}
for related results on the topology of spaces of symplectic
forms in dimension four.

On the other hand, since the round metric on $S^4$ has positive scalar curvature, we also obtain the following corollary.

\begin{cor}\label{cor: psc on S4}
The space $\mathcal R_+(S^4)$ is disconnected.\qed
\end{cor}

We next consider the homotopy groups of its individual components.

\begin{thm}\label{thm: pi 1 of psc metrics on S4}
Every connected component of $\mathcal R_+(S^4)$ has nontrivial fundamental group. The component containing the standard round metric also has nontrivial homotopy groups in infinitely many higher degrees.
\end{thm}

The study of spaces of psc metrics on spheres goes back at least to Hitchin, who proved that $\mathcal R_+(S^{8k})$ and $\mathcal R_+(S^{8k+1})$ are disconnected for $k\geq1$ \cite{hitchin1974harmonic}. Gromov and Lawson proved that $\mathcal R_+(S^7)$ has infinitely many connected components \cite{gromov1983positive}, and Carr established the same conclusion for $\mathcal R_+(S^{4k-1})$ for every $k\geq2$ \cite{carr1988construction}. The work of Botvinnik--Ebert--Randal-Williams \cite{BERW17}, together with Perlmutter's extension to dimension five \cite{Per17}, shows that for $n\geq5$, every connected component of $\mathcal R_+(S^n)$ has nontrivial homotopy groups in infinitely many degrees. The second assertion of \Cref{thm: pi 1 of psc metrics on S4} gives a four-dimensional analogue for the component containing the round metric.

In lower dimensions, Rosenberg and Stolz proved that $\mathcal R_+(S^2)$ is contractible \cite{stolz2014metrics}. For $\mathcal R_+(S^3)$, Marques proved path-connectedness \cite[Corollary~1.1]{marques2012deforming}, and Bamler and Kleiner subsequently proved contractibility \cite{bamler2019ricci}. Although Ruberman constructed closed simply connected $4$-manifolds whose spaces of psc metrics have infinitely many connected components \cite{ruberman2002positive}, it remained open whether $\mathcal R_+(S^4)$ is connected; see \cite{walsh2018aspects}. \Cref{cor: psc on S4} answers this question negatively.

It is natural to ask whether every smooth structure on a given smoothable topological $4$-manifold admits an exotic diffeomorphism. By \Cref{thm: main}, this holds for the topological $4$-sphere. Our next theorem shows that it also holds for the topological $K3$ surface. In particular, the case of the standard smooth structure answers the $K3$ part of \cite[Problem~4.68(b)]{BKR2026} affirmatively.

\begin{restatable}{thm}{MainKThreeTheorem}\label{thm: K3}
Every smooth homotopy $K3$ surface admits an exotic diffeomorphism.
\end{restatable}

Our exotic mapping class on $S^4$ is represented by the deck involution of a double cover branched along a $2$-knot in $S^4$. To detect this class, we use the reduced $\mathrm{Pin}(2)$-equivariant families Bauer--Furuta invariant of Lin and Mukherjee \cite{lin2025family}. For an oriented smooth $\mathbb F_2$-homology $4$-sphere $X$, evaluating this invariant on spin mapping tori of spin lifts of orientation-preserving diffeomorphisms of $X$ defines a group homomorphism
\[
\mathfrak{LM}_X\colon\pi_0\mathrm{Diff}^+(X)\longrightarrow\Z_2,
\]
which we call the \emph{Lin--Mukherjee invariant}. Although the $\mathrm{Pin}(2)$-equivariant Bauer--Furuta invariant of $X$ itself is the identity, the reduced invariant of a mapping torus need not vanish; see \Cref{lem: BF is identity,thm: defining Lin}.

We compute the Lin--Mukherjee invariant of deck involutions using real Seiberg--Witten theory. This theory studies solutions fixed by a compatible anti-linear lift of an involution to the spinor bundles. Its development includes work of Tian--Wang on real almost complex manifolds, Nakamura on free involutions, and Kato on involutions with fixed points \cite{Tian-Wang:2009,Nakamura:2013,Kato:2022}. Konno--Miyazawa--Taniguchi subsequently developed real Seiberg--Witten Floer homotopy types and relative Bauer--Furuta invariants \cite{Konno-Miyazawa-Taniguchi:2025,Konno-Miyazawa-Taniguchi:2023}.

Miyazawa applied real Bauer--Furuta theory to branched double covers to define invariants of embedded surfaces \cite{miyazawa2023gauge}. He used these invariants to distinguish infinitely many embeddings of $\mathbb{RP}^2$ in $S^4$ that are all topologically isotopic to the standard embedding. For a $2$-knot $F$ in $S^4$, the real Bauer--Furuta map of the branched double cover is a stable map between spheres of the same dimension, and the absolute value of its degree is an invariant of $F$. We extend this numerical invariant to smoothly embedded $2$-spheres $F$ in oriented smooth $\mathbb F_2$-homology $4$-spheres $X$. We denote it by $\mathfrak M(X,F)$ and call it the \emph{Miyazawa invariant}; see \Cref{defn: generalized Miyazawa}. Let $\Sigma(F)$ denote the double cover of $X$ branched along $F$. The following theorem expresses the Lin--Mukherjee invariant of its deck involution in terms of $\mathfrak M(X,F)$.

\begin{restatable}{thm}{LinToMiyazawa}\label{thm: Lin to Miyazawa}
Let $X$ be an oriented smooth $\mathbb F_2$-homology $4$-sphere, and let $F$ be a smoothly embedded $2$-sphere in $X$. If $\tau_F$ denotes the deck involution of $\Sigma(F)$, then
\[
\mathfrak{LM}_{\Sigma(F)}(\tau_F)
=\left(\frac{\mathfrak M(X,F)^2-1}{8}\bmod 2\right)
=\begin{cases}
0 & \text{if }\mathfrak M(X,F)\equiv\pm1\pmod{8},\\
1 & \text{if }\mathfrak M(X,F)\equiv\pm3\pmod{8}.
\end{cases}
\]
\end{restatable}

To apply the formula on the standard $S^4$, we need a $2$-knot whose branched double cover is standard and whose Miyazawa invariant is congruent to $\pm3$ modulo $8$. We use a Cappell--Shaneson $2$-knot $F$ in $S^4$ \cite{cs1976inequivalent}; see \Cref{sec: CS example construction} for the construction. Its branched double cover $\Sigma(F)$ is a Cappell--Shaneson homotopy $4$-sphere \cite{cappellshaneson76} in a family shown to be standard by Akbulut \cite{akbulut2010cappell} and Gompf \cite{gompf2010more}.\footnote{It remains open whether every Cappell--Shaneson homotopy $4$-sphere is standard. For standardness results and related handle constructions, see \cite{Akbulut-Kirby:1979,Aitchison-Rubinstein:1984,Akbulut-Kirby:1985,Gompf:1991,akbulut2010cappell,Meier-Zupan:2022,Kim-Yamada:2023,Iwaki:2025}.} Surgery on $F$ gives the mapping torus $X_B$ of the automorphism of $T^3$ induced by a matrix $B\in SL_3(\Z)$. Its nontrivial double cover is $X_{B^2}$. We equip $X_{B^2}$ with a spin structure for which the deck involution admits a spin lift whose square is $-1$. Composing this lift with the quaternionic structure on spinors gives a compatible anti-linear involution and hence a real spin structure on $X_{B^2}$. Let $BF^R_{X_{B^2}}$ denote the associated real Bauer--Furuta invariant. Its underlying nonequivariant stable map $S^0\to S^0$ has an integer degree, whose absolute value is identified with the Miyazawa invariant by the surgery formula of Miyazawa--Park--Taniguchi \cite[Proposition~7.3]{miyazawa2025satellite}:
\[
\mathfrak M(S^4,F)=\bigl|\deg BF^R_{X_{B^2}}\bigr|.
\]

We compute this degree using metrics on $X_{B^2}$ that are invariant under the deck involution and stretch the base circle of the mapping torus. For sufficiently large stretching parameter, we show that every real monopole is gauge equivalent to a configuration with constant coefficients in an explicit trivialization. For our choice of $B$, the resulting equations have a unique solution, a degenerate reducible. Its Kuranishi map is given near the origin by an explicit cubic map $\C\to\C$ that is homotopic to $z\mapsto\overline{z}^{\,3}$ through maps whose only zero is the origin. With our reference orientations, its local degree is therefore $-3$. Since the reducible is the only unperturbed solution up to gauge, this local degree equals the ordinary degree of the real Bauer--Furuta map. Thus
\[
\mathfrak M(S^4,F)=|-3|=3.
\]
By \Cref{thm: Lin to Miyazawa}, the deck involution $\tau_F$ satisfies $\mathfrak{LM}_{S^4}(\tau_F)=1$, so it is not smoothly isotopic to the identity. By smooth isotopy extension, we may choose a diffeomorphism $f_0$ smoothly isotopic to $\tau_F$ and supported in the interior of a standard $4$-ball. Applying the Alexander trick in this ball shows that $f_0$, and hence $\tau_F$, is topologically isotopic to the identity. Thus $\tau_F$ is an exotic diffeomorphism of $S^4$.

For any smooth $\mathbb F_2$-homology $4$-sphere $X$, form the connected sum $f=f_0\#\id_X$ away from the ball supporting $f_0$. By the connected-sum formula in \Cref{lem: connected sum and Lin invariant}, we have $\mathfrak{LM}_X(f)=1$. The topological isotopy on the ball extends by the identity over its complement in $X$, so $f$ is topologically isotopic to the identity. Thus $f$ is exotic. The same construction gives an exotic diffeomorphism on every smooth homotopy $K3$ surface $Z$. In this case, we distinguish $f_0\#\id_Z$ from the identity by composing the reduced family Bauer--Furuta invariant of $f_0$ with the Bauer--Furuta invariant of $Z$; see the proof of \Cref{thm: K3}.

If $g\in\mathcal R_+(X)$ could be joined to $f^\ast g$ by a path of psc metrics, this path would give a fiberwise psc metric on the mapping torus of $f$. By \Cref{cor: fiberwise psc metric implies vanishing of LM}, this would imply $\mathfrak{LM}_X(f)=0$, contrary to the construction. Thus $g$ and $f^\ast g$ lie in different components of $\mathcal R_+(X)$.

The Lin--Mukherjee invariant also distinguishes $\tau_F$ from twists along half-unknotted Montesinos twins in $S^4$. Gay--Hartman proved that the mapping classes of these twists generate a cyclic subgroup of order at most two \cite[Theorem~4 and Lemma~6]{gay2025relations}. A generator is represented by the diffeomorphism $f_{GH}$ appearing in the announcement of Gabai--Gay--Hartman; see \cite[Theorem~0.5]{Gay2026PCMI}. Using fiberwise psc surgery, we show in \Cref{prop: Montesinos twist vanishing} that $\mathfrak{LM}_{S^4}$ vanishes on every such twist and hence on the entire subgroup. In particular,
\[
\mathfrak{LM}_{S^4}(f_{GH})=0.
\]
This argument does not rely on the announced nontriviality of $f_{GH}$. Since $\mathfrak{LM}_{S^4}(\tau_F)=1$, the deck involution $\tau_F$ is not smoothly isotopic to any composition of twists along half-unknotted Montesinos twins. Combined with the proof of \Cref{thm: main}, this answers \cite[Question 7]{gay2025diffeomorphisms} in the negative.

Finally, we explain a comparison between the full equivariant 2-knot invariant from the equivariant Bauer--Furuta invariant and the Miyazawa invariant.  Let $F$ be a $2$-knot in $S^4$, and choose a spin lift $\tilde\tau_F$ of the deck involution of its branched double cover $\Sigma(F)$. With respect to its unique spin structure and the chosen lift, $\Sigma(F)$ has a $G_{\mathrm{odd}}$-equivariant Bauer--Furuta invariant represented by a stable map
\[
BF_{\Sigma(F)}\colon S^0\longrightarrow S^0,
\qquad G_{\mathrm{odd}}=\mathrm{Pin}(2)\times_{\pm1}\Z_4.
\]

\begin{thm}\label{thm: Miyazawa determines equivariant BF for S4}
The Miyazawa invariant $\mathfrak M(S^4,F)$ determines the $G_{\mathrm{odd}}$-equivariant stable homotopy class $[BF_{\Sigma(F)}]\in\pi_0^{G_{\mathrm{odd}}}(S^0)$.
\end{thm}

The proof uses techniques similar to those in \cite[Appendix~A.2]{KPT-cables}, where $O(2)$-equivariant maps between representation spheres are studied. The corresponding classification of $O(2)$-equivariant Bauer--Furuta invariants in terms of Miyazawa's degree invariant is given in \cite[Theorem~1.7 and Remark~1.8]{KPT-cables}.

\subsection*{Questions} We conclude the intoduction with some questions. First of all, since we now know the existence of an exotic diffeomorphism on $S^4$ and the K3 surface, we can ask whether there also exist exotic diffeomorphisms on other small simply-connected 4-manifolds. The smallest natural candidates are $\C P^2$ and $S^2 \times S^2$.

\begin{question}
Does $\C P^2$ admit an exotic diffeomorphism? Does $S^2\times S^2$?
\end{question}

A related question concerns stabilization. By work of Kreck and Quinn \cite{Kreck1979,Quinn1986}, every exotic diffeomorphism of a closed simply connected smooth $4$-manifold becomes smoothly isotopic to the identity after sufficiently many stabilizations by $S^2\times S^2$; see \cite[Theorem~1.1]{GGHKP2026} for the completion of Quinn's proof. The same conclusion holds for stabilizations by $Q=\C P^2\#\overline{\C P^2}$. Auckly--Kim--Melvin--Ruberman \cite{Auckly-Kim-Melvin-Ruberman:2015} constructed families of exotic diffeomorphisms for which a single stabilization by $S^2\times S^2$ suffices, whereas Lin \cite{Lin:2023-1} showed that the Dehn twist along the neck of $K3\#K3$ remains exotic after one such stabilization.

\begin{question}
What is the least number of stabilizations by $S^2\times S^2$ needed to make the exotic diffeomorphism $\tau_F$ of $S^4$ constructed in this paper smoothly isotopic to the identity? What is the corresponding number for stabilizations by $\C P^2\#\overline{\C P^2}$?
\end{question}

\subsection*{Organization}

In \Cref{sec: Lin-Mukherjee invariant}, we introduce the Lin--Mukherjee invariant and establish some of its important properties. In \Cref{sec: relation to miyazawa}, we express this invariant for deck involutions in terms of the Miyazawa invariant. In \Cref{sec: CS degree computation}, we compute the Miyazawa invariant of special Cappell--Shaneson $2$-knots. In \Cref{sec: top construction}, we use this computation to prove the main theorems and compare our diffeomorphism of $S^4$ with twists along half-unknotted Montesinos twins. The appendix gives an explicit proof of the equivariant stable homotopy calculation in \Cref{lem: computation of T}.

\subsection*{Notation}

For homotopic morphisms $f$ and $g$, we write $f\sim g$; if the homotopy is $G$-equivariant, we write $f\sim_G g$. We write $X\simeq Y$ if $X$ and $Y$ are homotopy equivalent or stably isomorphic, and $X\simeq_G Y$ if the homotopy equivalence or stable isomorphism is $G$-equivariant. For strict equivalences, such as homeomorphisms, diffeomorphisms, or isomorphisms, we write $X\cong Y$, or $X\cong_G Y$ if the equivalence is $G$-equivariant. For a closed subgroup $H$ of a Lie group $G$, we write $W_G H=N_G H/H$ for the Weyl group of $H$ in $G$. Given an object $X$ and a base space $B$, we denote the trivial $X$-family over $B$, considered as a parametrized object over $B$, by $\underline{X}_B$.

\subsection*{Universes}

Throughout the paper, for any compact Lie group $G$, we use the \emph{maximal semifree $G$-universe} $\mathcal{U}_G$ that is defined as the maximal universe satisfying the following property: \emph{the induced action of the identity component of $G$ on any finite-dimensional $G$-subrepresentation $V \subset \mathcal{U}_G$ is semifree}. For example:
\begin{itemize}
    \item If $G$ is finite, then $\mathcal{U}_G$ is the complete $G$-universe;
    \item If $G=S^1$, then $\mathcal{U}_G \cong_G \R^\infty \oplus \C^\infty$, where $S^1$ acts trivially on $\R$ and by complex multiplication on $\C$.
\end{itemize}
Then, for any closed subgroup $H\subset G$, the following statements hold:
\begin{itemize}
    \item There exists an $H$-equivariant isometric embedding $\mathrm{Res}^G_H \mathcal{U}_G \hookrightarrow \mathcal{U}_H$, which can be chosen to be an isomorphism if $[G:H]<\infty$;
    \item There exists an $W_G H$-equivariant embedding $\mathcal{U}_G^H \hookrightarrow \mathcal{U}_{W_G H}$.
\end{itemize}
For some concrete examples of maximal semifree $G$-universes that are incomplete, we have
\[
\mathcal{U}_{\mathrm{Pin}(2)} = \R^\infty \oplus \wt\R^\infty \oplus \H^\infty,\quad \mathcal{U}_{\mathrm{Pin}(2)\times \Z_4/\langle (-1,[2])\rangle} = \R^\infty \oplus \wt\R_j^\infty \oplus \wt\R_\tau^\infty \oplus \wt\R_I^\infty \oplus \H^\infty.
\]
For $\mathcal{U}_{\mathrm{Pin}(2)}$, $j$ acts by left multiplication on $\H$ and $-1$ on $\wt\R$. For $\mathcal{U}_{\mathrm{Pin}(2)\times \Z_4/\langle (-1,[2])\rangle}$, its summands are defined as in \Cref{table: universe summands}.

\begin{table}[!h]
    \begin{center}
    \begin{tabular}{ |>{\centering\arraybackslash}p{2cm}|| >{\centering\arraybackslash}p{3cm} |>{\centering\arraybackslash}p{3cm}|>{\centering\arraybackslash}p{3cm}|>{\centering\arraybackslash}p{3cm}|  }
    \hline
    Summand & Action of $z\in S^1$ & Action of $j$ & Action of $(0,[1])$ & Action of $(1,[1])$ \\
    \hline
    $\R$ & $\mathrm{id}$ & $\mathrm{id}$ & $\mathrm{id}$ & $\mathrm{id}$ \\
    $\wt\R_j$ & $\mathrm{id}$ & $-1$ & $\mathrm{id}$ & $-1$ \\
    $\wt\R_\tau$ & $\mathrm{id}$ & $\mathrm{id}$ & $-1$ & $-1$ \\
    $\wt\R_I$ & $\mathrm{id}$ & $-1$ & $-1$ & $\mathrm{id}$ \\
    $\H$ & $w\mapsto zw$ & $w\mapsto jw$ & $w\mapsto -wj$ & $w\mapsto -jwj$ \\
    \hline
    \end{tabular}
    \vspace{0.5em}
    \caption{Summands of the universe $\mathcal{U}_{\mathrm{Pin}(2)\times \Z_4/\langle (-1,[2])\rangle}$} \label{table: universe summands}
    \end{center}
    \end{table}

\subsection*{AI declaration}

The authors used ChatGPT 6 Astra Ultra and ChatGPT 6 Pro to generate initial drafts of the proofs of \Cref{prop:cs-constant-real-solutions} and \Cref{lem: nonvanishing K3 class} and to proofread and polish the paper. The authors subsequently rewrote both proofs in full.

\subsection*{Acknowledgements}
 The first author is partially supported by the Royal Society University Research Fellowship URF\textbackslash R1\textbackslash 251501. The second author is partially supported by the Samsung Science and Technology Foundation (SSTF-BA2102-02) and by National Research Foundation of Korea (NRF) grants funded by the Korean government (MSIT) (RS-2025-00542968 and RS-2026-25506097). The third author was partially supported by JSPS KAKENHI Grant Number 22K13921.

\section{The Lin--Mukherjee invariant}\label{sec: Lin-Mukherjee invariant}

One of the main tools in this paper is the families Bauer--Furuta invariant. For its construction, see \cite[Section~2]{szymik2010characteristic} and \cite[Section~2.3]{baraglia2022bauer}. In the spin setting, we retain the quaternionic symmetry and use the $\mathrm{Pin}(2)$-equivariant construction of \cite[Section~2.1]{lin2025family}. For a family over $S^1$, we take the reduced component of the invariant, following the convention of \cite[Remark~10]{lin2025family}.

Throughout this section, we fix an oriented smooth $\Z_2$-homology $4$-sphere $X$ and denote its unique spin structure by $\s_X$. The spin $4$-manifold $(X,\s_X)$ has a $\mathrm{Pin}(2)$-equivariant Bauer--Furuta invariant \cite{manolescu2007gluing}, represented by a stable map
\[
BF_{X,\s_X}\colon S^0\longrightarrow S^0.
\]
We first compute this map up to $\mathrm{Pin}(2)$-equivariant stable homotopy. Here, $S^0$ denotes the degree-zero $\mathrm{Pin}(2)$-equivariant sphere spectrum.

\begin{lem}\label{lem: BF is identity}
We have $BF_{X,\s_X}\sim_{\mathrm{Pin}(2)}\mathrm{id}$.
\end{lem}
\begin{proof}
We identify $BF_{X,\s_X}$ with its equivariant stable homotopy class
\[
[BF_{X,\s_X}]\in\pi_0^{\mathrm{Pin}(2)}(S^0).
\]
As discussed in \cite[proof of Proposition~3.9]{Ma16}, taking $S^1$-fixed points gives $BF_{X,\s_X}^{S^1}\sim_{\Z_2}\mathrm{id}$, where $\Z_2=\mathrm{Pin}(2)/S^1$. Hence
\[
[\mathrm{id}]-[BF_{X,\s_X}]
\in\ker\bigl(\pi_0^{\mathrm{Pin}(2)}(S^0)\longrightarrow\pi_0^{\Z_2}(S^0)\bigr)
\cong\pi_0^{\mathrm{Pin}(2)}(E\mathrm{Pin}(2)_+),
\]
where the identification follows from \cite[Lemma~21]{lin2025family}. Since the adjoint representation of $\mathrm{Pin}(2)$ is $\wt\R$, the Adams isomorphism \cite[Equation~(16)]{lin2025family} gives
\[
\pi_0^{\mathrm{Pin}(2)}(E\mathrm{Pin}(2)_+)
\cong\pi_0^s\left(\frac{\wt\R^+\wedge E\mathrm{Pin}(2)_+}{\mathrm{Pin}(2)}\right)
\cong\pi_0^s\bigl(\mathrm{Th}(\xi)\bigr),
\]
where $\xi$ denotes the line bundle over $B\mathrm{Pin}(2)$ associated with $\wt\R$. By the stable Hurewicz theorem and the connectedness of $\mathrm{Th}(\xi)$, we have
\[
\pi_0^s\bigl(\mathrm{Th}(\xi)\bigr)
\cong\wt{H}_0(\mathrm{Th}(\xi);\Z)=0.
\]
Thus the kernel above is trivial, so $[\mathrm{id}]-[BF_{X,\s_X}]=0$. Therefore $BF_{X,\s_X}\sim_{\mathrm{Pin}(2)}\mathrm{id}$, as desired.
\end{proof}

Let $\operatorname{Fr}^+(TX)$ denote the oriented frame bundle of $X$, and represent the spin structure $\s_X$ by a double covering
\[
\rho\colon P\longrightarrow\operatorname{Fr}^+(TX).
\]
Here $P$ is a principal $\widetilde{GL}^+(4,\mathbb R)$-bundle, where $\widetilde{GL}^+(4,\mathbb R)$ is the connected double cover of $GL^+(4,\mathbb R)$. The topological group $\mathrm{Diff}^+(X,\s_X)$ of \emph{spin diffeomorphisms} is the group of principal bundle automorphisms $\tilde f\colon P\to P$ satisfying
\[
\rho\circ\tilde f=\operatorname{Fr}^+(df)\circ\rho
\]
for some $f\in\mathrm{Diff}^+(X)$, equipped with the $C^\infty$ topology. The diffeomorphism $f$ is uniquely determined by $\tilde f$. We call $\tilde f$ a \emph{spin lift} of $f$ and say that $\tilde f$ \emph{covers} $f$.

Let $f\in\mathrm{Diff}^+(X)$ be an orientation-preserving self-diffeomorphism of $X$. Since $X$ is connected and $\s_X$ is unique, we have a fiber bundle
\[
\Z_2\cong\mathrm{Map}(X,\Z_2)\longrightarrow
\mathrm{Diff}^+(X,\s_X)\longrightarrow
\mathrm{Diff}^+(X),
\]
where the second map is the forgetful map.
Thus $f$ has exactly two spin lifts, related by fiberwise multiplication by $-1$. Choose one lift and denote it by $\tilde f$; the other is denoted by $-\tilde f$.

Consider the spin mapping torus
\[
\mathrm{Torus}(\tilde f)=(X\times I)/{\sim};
\qquad
(x,1)\sim(f(x),0),
\]
equipped with the fiberwise spin structure obtained by gluing $P\times I$ via $\tilde f$. Projection to the second coordinate modulo $1$ makes this a smooth $X$-bundle over $S^1$ with a fiberwise spin structure.

We first describe the family Bauer--Furuta invariant as a parametrized stable map over $S^1$. Choose a fiberwise Riemannian metric and the usual auxiliary data. After Sobolev completion, the family monopole maps form a $\mathrm{Pin}(2)$-equivariant map between Hilbert bundles over $S^1$. An equivariant finite-dimensional approximation gives a fiberwise based $\mathrm{Pin}(2)$-equivariant map
\[
\Phi_{\tilde f}\colon
S_{S^1}(U)\longrightarrow S_{S^1}(V),
\]
where $U,V\to S^1$ are finite-dimensional $\mathrm{Pin}(2)$-equivariant vector bundles and $S_{S^1}(U)$ denotes the fiberwise one-point compactification of $U$. Here $\mathrm{Pin}(2)$ acts fiberwise and acts trivially on the base $S^1$.

The stable difference between $U$ and $V$ is determined by the family index of the linear part of the monopole map. We now describe the stable index identifications used in our setting. Let $D$ denote the associated family of spin Dirac operators. Its quaternionic index has virtual rank $ 
-\frac{\sigma(X)}{16}=0$.
After stabilization, we may write $ 
\operatorname{Ind}_{\H}(D)
=
[K]-[S^1\times\H^N]$,
where $K\to S^1$ is a quaternionic vector bundle of rank $N$. Since $\mathrm{Sp}(N)$ is connected and simply connected, $K$ admits a quaternionic trivialization unique up to homotopy. Thus the Dirac index has a stable quaternionic trivialization unique up to homotopy; see
\cite[Section~2.1 and Remark~12]{lin2025family}.

For the real part, let $\Omega^0_0(X_t;i\R)$ denote the imaginary-valued functions with integral zero on the fiber $X_t$. With the usual Sobolev completions understood, the operators
\[
L_t=d_t^*\oplus d_t^+\colon
\Omega^1(X_t;i\R)
\longrightarrow
\Omega^0_0(X_t;i\R)\oplus\Omega^+(X_t;i\R)
\]
form a family of isomorphisms, since $b_1(X)=b^+(X)=0$. These isomorphisms are $\mathrm{Pin}(2)$-equivariant and compatible with the mapping-torus gluing. We use them to identify the real source and target summands in the finite-dimensional approximation.

Together with the stable quaternionic trivialization above, these identifications allow us, after stabilization, to identify the source and target bundles of the finite-dimensional approximation. Equivalently, after adding the same auxiliary $\mathrm{Pin}(2)$-equivariant vector bundle to the source and target if necessary, we may regard $\Phi_{\tilde f}$ as a fiberwise stable endomorphism of a common sphere bundle. Fiberwise desuspension by this common sphere bundle gives a parametrized $\mathrm{Pin}(2)$-equivariant stable map
\[
\widehat{\Phi}_{\tilde f}\colon
\underline{S^0}_{S^1}
\longrightarrow
\underline{S^0}_{S^1}
\]
over $S^1$.

We next explain how this parametrized stable map determines a $\mathrm{Pin}(2)$-equivariant stable cohomotopy class. 
Suppose that after stabilization the common vector bundle is trivialized as $S^1\times W$. Then a representative of the parametrized map has the form
\[
S^1\times S^W\longrightarrow S^1\times S^W,
\qquad
(t,x)\longmapsto (t,\Phi_t(x))
\]
where $ 
\Phi_t\colon S^W\longrightarrow S^W$
is a based $\mathrm{Pin}(2)$-equivariant map for each $t\in S^1$. 
Thus this map descends to the quotient $
\frac{S^1\times S^W}{S^1\times\{\infty\}}
\longrightarrow S^W$.
Via a natural identification $
\frac{S^1\times S^W}{S^1\times\{\infty\}}
\cong S^1_+\wedge S^W$, we get a $\mathrm{Pin}(2)$-equivariant map $S^1_+\wedge S^W \to S^W$.
Desuspending by $W$ gives an element of $
[\Sigma^\infty S^1_+,S^0]^{\mathrm{Pin}(2)}$.

Then, the stable splitting $ 
\Sigma^\infty S^1_+
\simeq
S^0\vee\Sigma S^0$
gives
\[
[\Sigma^\infty S^1_+,S^0]^{\mathrm{Pin}(2)}
\cong
\pi_0^{\mathrm{Pin}(2)}(S^0)
\oplus
\pi_1^{\mathrm{Pin}(2)}(S^0).
\]
The first component is obtained by restricting the parametrized map to a fiber over a point of $S^1$, and is therefore the $\mathrm{Pin}(2)$-equivariant Bauer--Furuta invariant of $(X,\s_X)$. The second component, which we call the \emph{reduced component}, records the additional family information. We denote this component by
\[
FBF(\tilde f)\in\pi_1^{\mathrm{Pin}(2)}(S^0).
\]
\footnote{This is denoted
$\mathrm{FBF}^{\mathrm{Pin}(2)}(X,f,\tilde f)$
in \cite{lin2025family}.}
Its $\mathrm{Pin}(2)$-equivariant stable homotopy class depends only on the isomorphism class of
$\mathrm{Torus}(\tilde f)\to S^1$
as a smooth bundle equipped with a fiberwise spin structure, and hence only on the isotopy class of the spin diffeomorphism $\tilde f$.

Taking its $S^1$-fixed points gives
\[
\left[FBF(\tilde f)^{S^1}\right]=0\in\pi_1^{\Z_2}(S^0),
\]
as shown in \cite[Lemma~31(i)]{lin2025family}. Thus, by \cite[Proposition~26]{lin2025family}, we have
\[
\left[FBF(\tilde f)\right]\in
\ker\left(
\pi_1^{\mathrm{Pin}(2)}(S^0)
\overset{\Phi^{S^1}}{\longrightarrow}
\pi_1^{\Z_2}(S^0)
\right)
=\Z_2\langle\alpha\rangle,
\]
where $\alpha$ denotes the unique nonzero element of $\pi_1^{\mathrm{Pin}(2)}(S^0)$ contain in the the kernel of the $S^1$-fixed point map to $\pi_1^{\Z_2}(S^0)$; see \cite[Theorem~8]{lin2025family}. We then define the invariant $\mathcal{LM}$ in \Cref{thm: defining Lin} as follows.

\begin{defn}
We define the \emph{Lin--Mukherjee invariant} of $f$ by
\[
\mathfrak{LM}_X(f):=[FBF(\tilde f)]\in\Z_2\langle\alpha\rangle\cong\Z_2.
\]
Equivalently,
\[
\mathfrak{LM}_X(f)=
\begin{cases}
0 & \text{if } FBF(\tilde f)\sim_{\mathrm{Pin}(2)}0,\\
1 & \text{if } FBF(\tilde f)\sim_{\mathrm{Pin}(2)}\alpha.
\end{cases}
\]
\end{defn}
A priori, the above definition depends on the choice of spin lift $\tilde f$, not just on $X$ and the isotopy class of $f$. We must show that $\tilde f$ and $-\tilde f$ give the same value and establish additivity. We prove these properties in the remainder of this section. Until independence of the lift is established, we write $\mathfrak{LM}_X(\tilde f)$ instead of $\mathfrak{LM}_X(f)$.

\begin{lem}\label{lem: isotopy invariance}
The value of $\mathfrak{LM}_X(\tilde f)$ depends only on the isotopy class of the spin diffeomorphism $\tilde f$.
\end{lem}
\begin{proof}
An isotopy of spin diffeomorphisms induces an isomorphism of their mapping tori preserving the fiberwise spin structures. The argument of \cite[Theorem~6.2]{szymik2010characteristic} therefore applies in the $\mathrm{Pin}(2)$-equivariant setting. The required quaternionic trivializations are unique up to homotopy: since $\mathrm{Sp}(n)$ is $2$-connected for every $n>0$, every quaternionic vector bundle over $S^k$, for $k=1,2$, admits a trivialization unique up to homotopy. See also \cite[Remark~12]{lin2025family}.
\end{proof}

\begin{lem}\label{lem: naturality of BF}
Let $B$ and $B'$ be finite 1-dimensional CW complexes, and let $f\colon B'\to B$ be a continuous map. Let $(W,\s_W)$ be a closed smooth spin $4$-manifold with $b_1(W)=0$, and let $E\to B$ be a smooth $W$-bundle equipped with a fiberwise spin structure $\s_E$ restricting to $\s_W$ on each fiber. Then the family Bauer--Furuta invariants satisfy
\[
f^\ast BF_{E,\s_E}\sim_{\mathrm{Pin}(2)}BF_{f^\ast(E,\s_E)}
\]
as parametrized stable maps over $B'$, where $f^\ast$ also pulls back the index bundles.
\end{lem}

\begin{proof}
This is the $\mathrm{Pin}(2)$-equivariant version of \cite[Theorems~2.1 and~3.1]{szymik2010characteristic}. Note that, since $B$ and $B'$ are 1-dimensional, the given families $(E,\s_E)$ and $f^\ast (E,\s_E)$ admit unique trivializations of their quaternionic Dirac index bundles; these trivializations are used to define $BF_{E,\s_E}$ and $BF_{f^\ast (E,\s_E)}$. Pulling back the family monopole map and its equivariant finite-dimensional approximation gives the corresponding construction for $f^\ast(E,\s_E)$. Since pullback preserves the quaternionic symmetry, the resulting identification is $\mathrm{Pin}(2)$-equivariant. 
\end{proof}

\begin{lem}\label{lem: pre-additivity of Lin--Mukherjee invariant}
For any $\tilde f,\tilde g\in\mathrm{Diff}^+(X,\s_X)$, we have
\[
\mathfrak{LM}_X(\tilde f\circ\tilde g)
=
\mathfrak{LM}_X(\tilde f)+\mathfrak{LM}_X(\tilde g).
\]
\end{lem}
\begin{proof}
 Since we have a stable splitting
\[
\Sigma^\infty (S^1 \vee S^1)_+ \simeq \Sigma^\infty S^0 \vee \Sigma^\infty S^1 \vee \Sigma^\infty S^1,
\]
we get
\[
[\Sigma^\infty (S^1 \vee S^1)_+,S^0]^{\mathrm{Pin}(2)}
\cong
\pi_0^{\mathrm{Pin}(2)}(S^0)
\oplus
\pi_1^{\mathrm{Pin}(2)}(S^0)
\oplus
\pi_1^{\mathrm{Pin}(2)}(S^0).
\]
Let $E\to S^1\vee S^1$ be an $(X,\s_X)$-bundle whose monodromies along the two circles are $\tilde f$ and $\tilde g$, respectively. Since $S^1 \vee S^1$ is a graph, its (rank 0) quaternionic index bundle admits a stable trivialization that is unique up to homotopy, so the family invariant $BF_E$ is well-defined up to stable homotopy. By restriction to the two circles, we see that (the stable homotopy class of) $BF_E$ corresponds to the triple
\[
\bigl(BF_{X,\s_X},FBF(\tilde f),FBF(\tilde g)\bigr)\in \pi_0^{\mathrm{Pin}(2)}(S^0)
\oplus
\pi_1^{\mathrm{Pin}(2)}(S^0)
\oplus
\pi_1^{\mathrm{Pin}(2)}(S^0).
\]
Consider the pinch map
\[
\gamma\colon S^1\longrightarrow S^1\vee S^1
\]
that traverse first the circle (with monodromy $\tilde g$) and then the circle (with monodromy $\tilde f$). Then we have
\[
\gamma^\ast E\cong\mathrm{Torus}(\tilde f\circ\tilde g)
\]
as $(X,\s_X)$-bundles over $S^1$. Under the given stable splittings, pullback along $\gamma$ preserves the $\pi_0$ component and adds the two $\pi_1$ components. Thus \Cref{lem: naturality of BF} gives
\[
BF_{\mathrm{Torus}(\tilde f\circ\tilde g)}=\gamma^\ast BF_E=\bigl(BF_{X,\s_X},FBF(\tilde f)+FBF(\tilde g)\bigr).
\]
as $\mathrm{Pin}(2)$-equivariant stable homotopy classes; note that, since $S^1$ and $S^1 \vee S^1$ are both 1-dimensional CW-complexes, the trivialization of quaternionic Dirac index bundles over them are unique up to homotopy, as discussed in the proof of \Cref{lem: isotopy invariance}. Taking the reduced component then proves the lemma.
\end{proof}

\begin{lem}\label{lem: -id is trivial}
Let $-\mathrm{id}_X$ denote the nontrivial spin lift of the diffeomorphism $\mathrm{id}_X$. Then $\mathfrak{LM}_X(-\mathrm{id}_X)=0$.
\end{lem}

\begin{proof}
Choose a stabilized $\Pin(2)$-equivariant finite-dimensional approximation $F\colon V^+\to V^+$ representing $BF_{X,\s_X}$. The spin automorphism $-\mathrm{id}_X$ acts on both source and target as the central element $-1\in\Pin(2)$: it acts by $-1$ on the quaternionic summands and trivially on the real summands. Using constant auxiliary data, the family over $\mathrm{Torus}(-\mathrm{id}_X)$ is therefore represented by gluing $F$ along this action.

The resulting representation bundles over $S^1$ are $\Pin(2)$-equivariantly trivial. Indeed, on the quaternionic summands, the clutching map $-I$ can be joined to $I$ within the commuting quaternionic frame group $\mathrm{Sp}(m)$, while the real summands have identity clutching. We may therefore give any further stabilizing summands the same central $-1$ action in the gluing.

By \Cref{lem: BF is identity}, after such a stabilization there is a $\Pin(2)$-equivariant homotopy from $F$ to $\mathrm{id}_{V^+}$. Every map in this homotopy commutes with $-1$, so the homotopy descends to the mapping-torus bundles. Its endpoint is the fiberwise identity. Using the same equivariant trivialization for source and target, with the real summands fixed and the quaternionic trivialization unique up to homotopy, this becomes the constant identity family. Its reduced component is zero, and hence $\mathfrak{LM}_X(-\mathrm{id}_X)=0$.
\end{proof}

We now show that, for an $\mathbb F_2$-homology $4$-sphere, the family Bauer--Furuta invariant of Lin--Mukherjee is independent of the choice of spin lift and defines a $\Z_2$-valued homomorphism on the mapping class group.

\begin{thm}\label{thm: defining Lin}
Let $X$ be a smooth, oriented $\mathbb F_2$-homology $4$-sphere. Then the Lin--Mukherjee invariant defines a group homomorphism
\[
\mathfrak{LM}_X\colon\pi_0\mathrm{Diff}^+(X)\longrightarrow\Z_2.
\]
\end{thm}

\begin{proof}
We first show that the invariant is independent of the choice of spin lift. Given $f\in\mathrm{Diff}^+(X)$, denote its two spin lifts by $\tilde f$ and $-\tilde f$. Since
\[
-\tilde f=(-\mathrm{id}_X)\circ\tilde f,
\]
\Cref{lem: pre-additivity of Lin--Mukherjee invariant,lem: -id is trivial} give
\[
\mathfrak{LM}_X(-\tilde f)
=\mathfrak{LM}_X(-\mathrm{id}_X)+\mathfrak{LM}_X(\tilde f)
=\mathfrak{LM}_X(\tilde f).
\]
We may therefore write $\mathfrak{LM}_X(f)$ instead of $\mathfrak{LM}_X(\tilde f)$.

The map $\mathrm{Diff}^+(X,\s_X)\to\mathrm{Diff}^+(X)$ is a double covering, so every isotopy of diffeomorphisms lifts to an isotopy of spin diffeomorphisms. By \Cref{lem: isotopy invariance} and independence of the lift, $\mathfrak{LM}_X(f)$ depends only on the isotopy class of $f$.

Finally, choose $f,g\in\mathrm{Diff}^+(X)$ and spin lifts $\tilde f,\tilde g\in\mathrm{Diff}^+(X,\s_X)$. Their composition $\tilde f\circ\tilde g$ is a spin lift of $f\circ g$, so \Cref{lem: pre-additivity of Lin--Mukherjee invariant} gives
\[
\begin{aligned}
\mathfrak{LM}_X(f\circ g)
&=\mathfrak{LM}_X(\tilde f\circ\tilde g)\\
&=\mathfrak{LM}_X(\tilde f)+\mathfrak{LM}_X(\tilde g)\\
&=\mathfrak{LM}_X(f)+\mathfrak{LM}_X(g).
\end{aligned}
\]
Thus $\mathfrak{LM}_X$ is a well-defined group homomorphism.
\end{proof}

Finally, we prove some properties of $\mathfrak{LM}_X$.

\begin{lem} \label{lem: space of psc metrics and lin invariant}
    Suppose that $X$ admits an orientation-preserving self-diffeomorphism $f$ satisfying $\mathfrak{LM}_X(f)=1$. Then, for any psc metric $g$ on $X$, $g$ and $f^\ast g$ are contained in different connected components of the space $\cR_+(X)$ of psc metrics on $X$. Consequently, $\cR_+(X)$ is either empty or disconnected.
\end{lem}

\begin{proof}
    Suppose that $g$ and $f^\ast g$ are contained in the same connected component. Since $\cR_+(X)$ is a Fr\'{e}chet manifold, it is locally path-connected. Hence there exists a path between $g$ and $f^\ast g$. This path, after smoothing and making it constant near its endpoints, defines a fiberwise psc metric on the spin mapping torus $\mathrm{Torus}(\tilde{f})$ where $\tilde{f}$ denotes a spin lift of $f$. By \cite[Lemma 2.8]{baraglia2026exotic}, we have a $\mathrm{Pin}(2)$-equivariant parametrized (over the base $S^1$) homotopy
    \[
    BF_{\mathrm{Torus}(\tilde{f})} \sim_{\mathrm{Pin}(2)} \left[ \underline{S^0}_{S^1} \xrightarrow{\mathrm{id}} \underline{S^0}_{S^1} \right].
    \]
    Hence we can represent this family Bauer--Furuta invariant as
    \[
    \left[BF_{\mathrm{Torus}(\tilde{f})}\right] = ([\mathrm{id}],0) \in \pi_0^{\mathrm{Pin}(2)}(S^0)\oplus \pi_1^{\mathrm{Pin}(2)}(S^0),
    \]
    which gives $FBF(\tilde{f})=0$ and thus $\mathfrak{LM}_X(f)=0$, a contradiction. Therefore $g$ and $f^\ast g$ lie in different connected components of $\cR_+(X)$. In particular, $\cR_+(X)$ is either empty or disconnected.
\end{proof}

\begin{cor} \label{cor: fiberwise psc metric implies vanishing of LM}
    Let $f$ be an orientation-preserving self-diffeomorphism of $X$. If the mapping torus family $\mathrm{Torus}(f)\rightarrow S^1$ admits a fiberwise psc metric, then $\mathfrak{LM}_X(f)=0$.
\end{cor}
\begin{proof}
    This was already proven in the proof of \Cref{lem: space of psc metrics and lin invariant}.
\end{proof}

\begin{cor}\label{cor: psc vanishing}
Suppose that $f$ is an orientation-preserving diffeomorphism satisfying $f^*g=g$ for some Riemannian metric $g$ on $X$ that has positive scalar curvature at every point. Then we have $\mathfrak{LM}_X(f)=0$.
\end{cor}
\begin{proof}
    This follows directly from \Cref{lem: space of psc metrics and lin invariant}.
\end{proof}

\begin{lem} \label{lem: connected sum FBF}
Let $Z$ be any closed smooth spin 4-manifold with $H^1(Z;\mathbb{F}_2)=0$ and choose $f\in \mathrm{Diff}^+(X)$ fixing an open ball $U\subset X$ pointwise, so that the connected sum $f\# \mathrm{id}_Z$ is again an orientation-preserving self-diffeomorphism of $X\# Z = (X\smallsetminus U) \cup (Z\smallsetminus B^4)$. Then there exists a spin lift $F$ of $f\# \mathrm{id}_Z$ such that the family Bauer--Furuta invariant $BF_{\mathrm{Torus}(F),\s_F}$ of the spin mapping torus family $\mathrm{Torus}(F)\rightarrow S^1$, together with the fiberwise spin structure $\s_F$ induced by $F$, satisfies
\[
BF_{\mathrm{Torus}(F),\s_F} \sim_{\mathrm{Pin}(2)} BF_{Z,\s_Z}\circ (\mathrm{id},\mathfrak{LM}_X(f)\cdot \alpha),
\]
where $\s_Z$ denotes the unique spin structure on $Z$ and $(\mathrm{id},\alpha) \in \pi_0^{\mathrm{Pin}(2)}(S^0) \oplus \pi_1^{\mathrm{Pin}(2)}(S^0)$ is canonically understood as an element of $[\Sigma^\infty S^1 _+,S^0]^{\mathrm{Pin}(2)}$.
\end{lem}
\begin{proof}
    Choose the unique spin lift $\tilde{f}$ of $f$ that acts by identity on $U$; this specifies the spin lift $\tilde{f}\# \mathrm{id}$ of $f\# \mathrm{id}$. Since $BF_{X,\s_X} \sim_{\mathrm{Pin}(2)} \mathrm{id}$ by \Cref{lem: BF is identity}, we get
    \[
    BF_{\mathrm{Torus}(\tilde{f})} = (\mathrm{id},\mathfrak{LM}_X(f)\cdot \alpha) \in \pi_0^{\mathrm{Pin}(2)}(S^0) \oplus \pi_1^{\mathrm{Pin}(2)}(S^0) \cong [\Sigma^\infty S^1_+,S^0]^{\mathrm{Pin}(2)}.
    \]
    The chosen lift identifies the spin mapping torus of $\tilde f\#\mathrm{id}$ with the fiberwise connected sum 
    \[
    \mathrm{Torus}(\tilde f) \#_{S^1}(Z\times S^1,\operatorname{pr}_Z^*\s_Z) \rightarrow S^1
    \]
    formed using the product neighborhoods of the chosen balls, where $\mathrm{pr}_Z:Z\times S^1\rightarrow Z$ denotes the projection. The lemma then follows from the $\mathrm{Pin}(2)$-equivariant family connected-sum argument of \cite[Proof of Proposition 2.26]{Lin:2023-1}, applied to the full parametrized monopole maps before passing to the reduced invariant. Note that the same gluing argument applies in the present $b_1(Z)=0$ (instead of requiring $Z$ to be simply-connected) setting with the spin structures and product neighborhoods fixed.
\end{proof}

\begin{lem} \label{lem: connected sum and Lin invariant}
    Let $X,Y$ be $\mathbb{F}_2$-homology 4-spheres and $f,g$ be their orientation-preserving self-diffeomorphisms, respectively. Suppose that there exist open balls $U_X \subset X$ and $U_Y \subset Y$ such that $f\vert_{U_X}=\mathrm{id}_{U_X}$ and $g\vert_{U_Y} = \mathrm{id}_{U_Y}$. Consider the connected sum diffeomorphism
    \[
    f\# g \in \mathrm{Diff}^+(X\# Y)
    \]
    constructed by taking $X\# Y = (X\smallsetminus U_X) \cup (Y\smallsetminus U_Y)$. Then we have
    \[
    \mathfrak{LM}_{X\# Y}(f\# g) = \mathfrak{LM}_X(f)+\mathfrak{LM}_Y(g).
    \]
\end{lem}
\begin{proof}
    This follows from adapting the proof of \Cref{lem: connected sum FBF} to the $\mathbb{F}_2$-homology 4-sphere setting, using
    \[
    f\# g = (f\# \mathrm{id}_Y) \circ (\mathrm{id}_X \# g).
    \]
\end{proof}

\section{Relation to Miyazawa invariant} \label{sec: relation to miyazawa}

Throughout this section, we fix an $\mathbb{F}_2$-homology 4-sphere $X$ and a smoothly embedded sphere $F\subset X$. Since $H_2(X;\Z)$ is 2-torsion-free and finite, and thus $[F]$ is divisible by 2, it admits a uniquely defined branched double cover $\Sigma(F)$, which is again a $\mathbb{F}_2$-homology 4-sphere. Clearly, $\Sigma(F)$ has a unique spin structure, which we denote by $\s$.

Consider the deck transformation $\tau$ of the branched double covering map $\Sigma(F)\rightarrow X$. Since its fixed point set is nonempty and has codimension 2, we know from \cite[Proposition 8.46]{atiyah1968lefschetz} that it is odd, i.e. lifts to an order 4 spin diffeomorphism $\tilde\tau$ of $(\Sigma(F),\s)$ satisfying $\tilde\tau^2 = -\mathrm{id}_{\Sigma(F)}$, which gives rise to the following symmetry group:
\[
G_\mathrm{odd} = \mathrm{Pin}(2)\times \Z_4 \langle \tilde\tau \rangle / (-1,\tilde\tau^2).
\]
 There are exactly two such structures coming from $\tilde\tau$ and $-\tilde\tau$.  
Choosing $\tilde\tau$-equivariant auxiliary data and performing an equivariant finite-dimensional approximation of the monopole map gives the $G_\mathrm{odd}$-equivariant Bauer--Furuta invariant
\[
BF_{\Sigma(F),\s,\tilde\tau} : S^0 \rightarrow S^0;
\]
see \cite[Section~6]{Ba21} and \cite[Theorem~2.6]{BF04}. As a result, we obtain its equivariant homotopy class $\left[ BF_{\Sigma(F),\s,\tilde\tau} \right] \in \pi_0^{G_\mathrm{odd}} (S^0)$. 
\begin{rem}
For the universes, we may take
$\mathcal{U}_{\mathrm{Pin}(2)}$
for the $\mathrm{Pin}(2)$-equivariant theory and
$\mathcal{U}_{G_{\mathrm{odd}}}$
for the $G_{\mathrm{odd}}$-equivariant theory, as shown in the
introduction.

For the real part of the monopole map, one can check that the required representation
types appear easily. Let us explain the quaternionic part. Recall that the
original $\mathrm{Pin}(2)$-action on $\H$ is given by left
quaternionic multiplication. The lift of $\tilde\tau$ is required to
commute with this $\mathrm{Pin}(2)$-action. Once the lift of the
underlying involution is fixed, there are only two possible lifts,
which differ by multiplication by the central element $-1$. In our
quaternionic coordinates, they are therefore given by
\[
w\longmapsto wj
\qquad\text{or}\qquad
w\longmapsto -wj.
\]
These two $G_{\mathrm{odd}}$-representations are isomorphic via $w\mapsto wi$.

In either case, the square of the action of $\tilde\tau$ is
multiplication by $-1$. Thus the element
$(-1,[2])\in \mathrm{Pin}(2)\times\Z_4$
acts trivially on $\H$, and the action descends to $ 
G_{\mathrm{odd}}
=
\mathrm{Pin}(2)\times\Z_4/
\langle(-1,[2])\rangle.$
\end{rem}

We use the following compatible index identifications in the
equivariant and family constructions. Let $D$ be the spin Dirac
operator on $(\Sigma(F),\s)$, with $\tau$-invariant auxiliary data.
On the spinor spaces, $\tilde\tau$ is quaternionic-linear and
satisfies $\tilde\tau^2=-1$. Since every orthogonal
quaternionic-linear complex structure on $\H^N$ is conjugate
to the standard one under $\mathrm{Sp}(N)$, the spinor
representations in a finite-dimensional approximation are sums
of the representation $\H$ in \Cref{table: universe summands}.
Consequently,
\[
\operatorname{Ind}_{G_{\mathrm{odd}}}(D)
=
-\frac{\sigma(\Sigma(F))}{16}[\H]
=
0.
\]
The equivariant identifications of the stabilized spinor source
and target are unique up to homotopy, since the group of orthogonal $G_{\mathrm{odd}}$-equivariant automorphisms of $\H^N$
is $\mathrm{U}(N)$\footnote{Indeed, a real-linear endomorphism of $\H$ commuting with
the left $\mathrm{Pin}(2)$-action is right multiplication
by some $a\in\H$.
It commutes with $\tilde\tau(w)=-wj$ precisely when
$aj=ja$, or equivalently when
$a\in\C_j:=\R\oplus\R \cdot j\cong\C$.
Hence
\[
\operatorname{End}_{G_{\mathrm{odd}}}(\H^N)
\cong M_N(\C_j),
\]
with the adjoint corresponding to conjugate transpose.
The orthogonal automorphisms therefore form $\mathrm{U}(N)$.}, which is connected. For the real summands,
we use the isomorphism $d^*\oplus d^+$ with the constant functions
removed from its target, as in the preceding section.
This isomorphism is $G_{\mathrm{odd}}$-equivariant for the chosen auxiliary data.

Taking mapping-torus bundles along $\tilde\tau$ carries these
identifications to the index identifications for
$\mathrm{Torus}(\tilde\tau)\to S^1$.
On the spinor summands, the induced identification is quaternionic
and therefore agrees up to homotopy with the stable
trivialization used in the family construction, by
\cite[Remark~12]{lin2025family}.
On the real summands, the induced identification is precisely
the family of normalized AHS operators used there.
The associated real bundles need not be trivial; they are
identified by this family of operators and cancelled by
fiberwise desuspension. Thus the equivariant construction
and its mapping-torus family use the same index normalizations.

Replacing $\tilde\tau$ by $-\tilde\tau$ replaces the
complex structures $J_U$ and $J_V$ on the spinor source
and target by their negatives.
This does not change the space of equivariant index
identifications, since
\[
aJ_U=J_Va
\quad\Longleftrightarrow\quad
a(-J_U)=(-J_V)a.
\]
In particular, the centralizers of $J$ and $-J$
in $\mathrm{Sp}(N)$ coincide.
The family construction uses stable quaternionic
trivializations, so the comparison with the family
index identifications applies to either choice of lift.
The resulting reduced family class is independent of
the lift by
\Cref{lem: pre-additivity of Lin--Mukherjee invariant,%
lem: -id is trivial}, since
\[
[FBF(-\tilde\tau)]
=
[FBF(\tilde\tau)]
+
[FBF(-\mathrm{id}_{\Sigma(F)})]
=
[FBF(\tilde\tau)].
\]

\begin{rem} \label{rem: S1 fixed point identity}
Applying the arguments outlined in \cite[Sketch of proof of Theorem 1.14(2)]{hopkins2022intersection} to the $G_\mathrm{odd}$-equivariant setting gives 
\[
BF_{\Sigma(F),\s,\tilde\tau}^{S^1} \sim_{G_\mathrm{odd}/S^1} \mathrm{id}.
\]
This argument is originally from \cite{Fu01}.
\end{rem}

Furthermore, if we consider the $I$-fixed point map
\[
BF_{\Sigma(F),\s,\tilde\tau}^{\langle I \rangle}:S^0 \rightarrow S^0,
\]
where $I=j\tilde\tau \in G_\mathrm{odd}$ (which has order 2), then $\left\vert\deg BF_{\Sigma(F),\s,\tilde\tau}^{\langle I \rangle}\right\vert\in \Z_{>0}$ is well-defined and independent of the choice of the lift $\tilde\tau$ of $\tau$. This is because replacing $\tilde\tau$ with $-\tilde\tau$ replaces $I$ with $-I$, which is conjugate to $I$ as $iIi^{-1} = -I$. Here, we are taking absolute value to follow the definition of Miyazawa invariant in \cite{miyazawa2023gauge}, where it was taken to mitigate the absence of a canonical orientation of virtual real index bundles.

\begin{defn} \label{defn: generalized Miyazawa}
    We define the \emph{(generalized) Miyazawa invariant} $\mathfrak{M}(X,F)$ of $(X,F)$ as
    \[
    \mathfrak{M}(X,F) := \left\vert\deg BF_{\Sigma(F),\s,\tilde\tau}^{\langle I \rangle}\right\vert.
    \]
\end{defn}

\begin{rem}
    If $X=S^4$, i.e. $F$ is a 2-knot, then $\mathfrak{M}(S^4,F)$ is the mapping invariant of $F$ defined in \cite{miyazawa2023gauge}.
\end{rem}

We now consider the equivariant Becker--Gottlieb transfer. Given a compact Lie group $K$ and a closed smooth manifolds $Z,Z'$ equipped with smooth $K$-actions, as well as a proper $K$-equivariant submersion $f:Z\rightarrow Z'$. Suppose that $Z$ can be $K$-equivariantly and smoothly embedded in some finite-dimensional $K$-representation $V$ contained in a prescribed universe (potentially incomplete) for stable $K$-equivariant homotopy theory; choose such an embedding $i:Z\rightarrow V$. Then the map
\[
(f,i):Z\rightarrow Z'\times V
\]
is also an embedding; consider its normal bundle $\nu$. Since we have $\nu \oplus \ker df = \underline{V}$, we can consider the following composite map:
\[
\mathrm{Tr}^K_f:V^+ \wedge Z'_+ \xrightarrow{\text{collapse}} \mathrm{Th}(\nu)\xrightarrow{\text{zero-section of }\ker (df)}\mathrm{Th}(\nu\oplus \ker (df))\simeq V^+ \wedge Z_+.
\]
After desuspending by $V$, we get a $K$-equivariant stable map
\[
\mathrm{Tr}^K_f:\Sigma^\infty Z'_+ \rightarrow \Sigma^\infty Z_+;
\]
we call this the \emph{$K$-equivariant Becker--Gottlieb transfer} of the given map $f:Z\rightarrow Z'$. For any closed subgroup $L\subset K$, we have $\left( \mathrm{Tr}^K_f \right)^L \sim_{N_K(L)/L} \mathrm{Tr}^{N_K(L)/L}_{f^L}$; this follows from \cite[Equations 6.2 and 6.3]{may1996equivariant} and \cite[Proposition 15.2.4(i)]{may2006parametrized}.

There is also a related construction. Instead of using the zero-section map, we can consider the following composition, using again the isomorphism $\nu \cong V\oplus (-\ker df)$:
\[
f^!: \Sigma^V \Sigma^\infty Z'_+ \xrightarrow{\text{collapse}} \Sigma^\infty\mathrm{Th}(\nu) \xrightarrow\simeq \Sigma^V \Sigma^\infty \mathrm{Th}(-\ker df),
\]
which induces a $K$-equivariant stable map
\[
f^!:\Sigma^\infty Z'_+ \rightarrow \Sigma^\infty\mathrm{Th}(-\ker df),
\]
which is defined as the $K$-equivariant \emph{Pontryagin--Thom umkehr map} of $f$; see \cite[Section 3]{BECKER19751}. Furthermore, if $\ker df$ is equivariantly stably trivial, then choosing its equivariant stable framing $\mu:\underline{W}\simeq \ker df$ gives the $K$-equivariant \emph{framed umkehr map}
\[
f^!_\mu:\Sigma^W \Sigma^\infty Z'_+ \rightarrow \Sigma^\infty Z_+.
\]
Moreover, this construction can be generalized to the setting where $f$ is a proper smooth map, not necessarily a submersion, with a framing of its virtual normal bundle $\nu_f = f^\ast TZ'-TZ$.

\begin{rem} \label{rem: rephrasing of definition of alpha}
The definition of $\alpha$ in \cite[Equation~32]{lin2025family} can be rephrased using umkehr maps as follows:
\[
\alpha = \left[ c_{\mathrm{Pin}(2)}\circ  (\mathrm{Pin}(2)\rightarrow \ast)^!_\mu \right],
\]
where $c_{\mathrm{Pin}(2)}: \Sigma^\infty\mathrm{Pin}(2)_+\rightarrow \Sigma^\infty S^0$ denotes the $\mathrm{Pin}(2)$-equivariant collapse map and $\mu$ denotes the unique (up to sign) left-invariant parallelization of $\mathrm{Pin}(2)$. Here the above $V$ can be taken as $\H$ with the natural embedding $\mathrm{Pin}(2) \to \H$.
\end{rem}

\begin{rem}
    In general, the equivariant Becker--Gottlieb transfer can be defined in a much general setting. Given $K$-CW-complexes $Z,Z'$ and a $K$-equivariant map $f:Z\rightarrow Z'$, for each $b\in Z'$, consider
    \[
    K_b = \{g\in K\,\vert\,gb=b\} \subset K,\quad F_b = \mathrm{hofib}_b(f);
    \]
    note that $F_b$ has a natural $K_b$-action. Suppose that $F_b$ is dualizable in the universe for $K_b$ (induced by the one for $K$). Then we have a well-defined $K$-equivariant Becker--Gottlieb transfer map
    \[
    \mathrm{Tr}_f^K:\Sigma^\infty Z'_+ \rightarrow \Sigma^\infty Z_+.
    \]
    See \cite[Section 6.1]{carmeli2025characters} for more detailed explanations.
\end{rem}

Now consider the composition of the equivariant Becker--Gottlieb transfer map with the collapsing map:
\[
\chi_K(f) : \Sigma^\infty Z'_+ \rightarrow \Sigma^\infty Z_+ \xrightarrow{\text{collapse}} \Sigma^\infty S^0.
\]
This can be seen as a stable $K$-equivariant cohomotopy class:
\[
[\chi_K(f)] \in \pi_K^0 (Z'),
\]
which is the \emph{$K$-equivariant family Euler characteristic} of $f$ (with respect to the given universe, where the admissible $K$-representation $V$ lies in). For any $b\in Z'$, pulling this class back along the inclusion $\{b\}\hookrightarrow Z'$, with stabilizer $K_b$, gives the $K_b$-equivariant Euler characteristic of the homotopy fiber of $f$ at $b$.

For any closed smooth $K$-manifold $Z$ (that can be $K$-equivariantly embedded in the relevant universe), consider the collapsing map $c_Z:Z\rightarrow \ast$. Then we have
\[
\chi_K(Z) := \chi_K(c_Z) :\Sigma^\infty S^0 \xrightarrow{\mathrm{Tr}^K_{c_Z}} \Sigma^\infty Z_+ \xrightarrow{\text{collapse}}\Sigma^\infty S^0,
\]
which can be considered as an element of $\pi_0^K(S^0)$. This is called the \emph{$K$-equivariant Euler characteristic} of $Z$. For any closed subgroup $L\subset K$, we clearly have
\[
\chi_K(Z)^L = \chi_{N_K(L)/L}(Z^L).
\]
Furthermore, if $K=1$, then $\chi_{\{1\}}(Z) = \chi(Z)$, i.e. the Euler characteristic of $Z$; this is \cite[Theorem 2.4]{BECKER19751}.

The following gives a classification of equivariant maps between spheres, using tom Dieck spilittings and Adams isomorphisms. 
\begin{thm}
\label{thm: tom dieck and eqv euler characteristic}
Let $K$ be a compact Lie group and let $\mathcal U$ be
a $K$-universe, not necessarily complete.
Denote by $\mathcal C(K)$ the set of conjugacy classes
of subgroups of $K$ which occur as isotropy groups
of points in $\mathcal U$. For each $H\in\mathcal C(K)$, let
\[
\xi_H:=EW_K H\times_{W_K H}\operatorname{Lie}(W_K H)
\longrightarrow BW_K H
\]
be the vector bundle associated with the adjoint
representation of $W_K H$. Then the following statements hold.

\begin{enumerate}[label=\textup{(\roman*)}]

\item
For every $H\in\mathcal C(K)$, the Adams
isomorphism gives
\[
\operatorname{Ad}_H:
\pi_0^s\bigl(\operatorname{Th}(\xi_H)\bigr)
\xrightarrow{\cong}
\pi_{0,\mathcal U^H}^{W_G H}(EW_K H{}_+).
\]

\item
The tom Dieck splitting gives an isomorphism
\[
\operatorname{TD}_K
:=\sum_{H\in\mathcal C(K)}\zeta_H^K:
\bigoplus_{H\in\mathcal C(K)}
\pi_{0,\mathcal U^H}^{W_K H}(EW_K H{}_+)
\xrightarrow{\cong}
\pi_{0,\mathcal U}^K(S^0),
\]
where $\zeta_H^K$ denotes the map associated
with the $H$-summand of the tom Dieck splitting.

\item
Consequently, their composition
\[
\begin{aligned}
\Theta_K
&:=\operatorname{TD}_K\circ
\left(\bigoplus_{H\in\mathcal C(K)}
\operatorname{Ad}_H\right):
\bigoplus_{H\in\mathcal C(K)}
\pi_0^s\bigl(\operatorname{Th}(\xi_H)\bigr)
\xrightarrow{\cong}
\pi_{0,\mathcal U}^K(S^0)
\end{aligned}
\]
is an isomorphism.

\item We have
\[
\pi_0^s\bigl(\operatorname{Th}(\xi_H)\bigr)
\cong
\begin{cases}
\mathbb Z &\text{if } W_K H \text{ is finite},\\
0 &\text{if } \dim W_K H>0.
\end{cases}
\]
For each $H$ with finite $W_K H$, let $1_H$
denote the generator of the corresponding
$\mathbb Z$-summand. Then
\[
\Theta_K(1_H)=\chi_K(K/H).
\]
In particular, the elements $\chi_K(K/H)$,
with $H\in\mathcal C(K)$ and $W_K H$ finite,
form a free abelian basis of
$\pi_{0,\mathcal U}^K(S^0)$.
\end{enumerate}
\end{thm}

\begin{proof}
Part \textup{(ii)} is the tom Dieck splitting theorem.
For the classical complete-universe version, see
\cite[Satz~2]{tomDieck1975Orbittypen} and
\cite[Chapter~V, Section~9]{LewisMaySteinberger}.
For the formulation for incomplete universes used here, see
\cite[Section~4.1]{lin2025family}.
The summand indexed by $H$ is $\pi_{0,\mathcal U^H}^{W_KH}(EW_KH{}_+)$,
since $(S^0)^H=S^0$.

For part \textup{(i)}, choose a complete $W_K H$-universe
$\widehat{\mathcal U}_H$ containing $\mathcal U^H$.
For every based $W_K H$-CW complex $Y$ that is free away from
the basepoint, change of universe induces an isomorphism
\begin{equation} \label{eqn: change of universe}
\pi_{q,\mathcal U^H}^{W_K H}(Y)
\rightarrow
\pi_{q,\widehat{\mathcal U}_H}^{W_K H}(Y).
\end{equation}
Indeed, the assumed equivariant embedding of $K/H$ to $\mathcal{U}$, and thus of $W_K H$ into
$\mathcal U^H$ by taking $H$-fixed points, allows us to apply the Wirthm\"uller isomorphism
for $\{1\}\subset W_K H$ in both universes; see
\cite[Theorem~1.1]{may2003wirthmuller}.
For a free cell $Y=W_K H{}_+\wedge S^n$, both groups identify
naturally with
\[
\pi_{q-n-\dim W_K H}^s(S^0),
\]
and the change-of-universe map is the identity under these
identifications. The assertion follows by cellular induction
and passage to filtered colimits.

Applying this comparison to $EW_K H_+$ and then using the
Adams isomorphism (See \cite[Theorem~2.2]{may2003wirthmuller} for complete universe version, and \cite[Theorem 12.2]{lewis2000splitting} for the incomplete universe version), we obtain
\[
\begin{aligned}
\pi_{0,\mathcal U^H}^{W_K H}(EW_K H_+)
&\cong
\pi_{0,\widehat{\mathcal U}_H}^{W_K H}(EW_K H_+)\\
&\cong
\pi_0^s\left(
  (\mathrm{Lie}(W_K H)^+\wedge EW_K H_+)/W_K H
\right)\\
&\cong
\pi_0^s\bigl(\operatorname{Th}(\xi_H)\bigr).
\end{aligned}
\]
The inverse of this composite is the stated isomorphism
$\operatorname{Ad}_H$. Note that this did not require us to assume that $\mathfrak{w}_H$ can be $W_K H$-equivariantly embedded in $\mathcal{U}^H$, as suspension by the representation sphere $\mathrm{Lie}(W_K H)^+$ happens in the completed universe $\widehat{\mathcal U}_H$.

Part \textup{(iii)} follows by composing
the isomorphisms in \textup{(i)} and
\textup{(ii)}.

For part \textup{(iv)}, Since $\xi_H$ has rank $\dim W_K H$, its
Thom spectrum is connective, and the
stable Hurewicz theorem gives
\[
\pi_0^s(\operatorname{Th}(\xi_H))
\cong
\widetilde H_0(\operatorname{Th}(\xi_H);
\mathbb Z).
\]
This group is $\mathbb Z$ when $W_K H$ is
finite, and zero when $W_K H$ has positive
dimension.

Finally, suppose that $W_K H$ is finite.
Under the Adams isomorphism, the generator
of $\pi_0^s(BW_K H_+)\cong\mathbb Z$
corresponds to the stable class of a free
$W_K H$-orbit in $EW_K H$.
We compute the following normalizer transfer in a complete $K$-universe, since change of universe is injective in degree zero and preserves both the standard generator of the $H$-summand and the equivariant Euler characteristic of $K/H$ (see \cite[Example~6.1]{lewis1995change}).
The construction of the tom Dieck
summand map $\zeta_H^K$, using the
transfer from $N_K(H)$ to $K$, identifies
its image with the equivariant
Euler characteristic of $K/H$ \emph{in the complete universe of $K$}\footnote{This is because the transfer from $N_K(H)$ to $K$ requires $K/N_K(H)$ to be $K$-equivariantly dualizable in the given $K$-universe $\mathcal{U}$; in general, this need not be true.}, i.e. the change-of-universe map
\[
\iota_\ast:\pi_{0,\mathcal{U}}^K (S^0)\rightarrow \pi_{0,\hat{\mathcal{U}}}^K(S^0)
\]
to our (possibly incomplete) $K$-universe $\mathcal{U}$ to the complete universe $\hat{\mathcal{U}}$ maps $\Theta_K(1_H)$ to $\chi_K(K/H)$. Observe that part (ii) of the theorem, together with the fact that change-of-universe maps are compatible with tom Dieck splittings and change-of-universe maps as in \Cref{eqn: change of universe} (i.e. for spaces with free actions outside the basepoint) induce isomorphisms, implies that $\iota_\ast$ is injective, and in fact, the inclusion map of a direct summand of $\pi_{0,\hat{\mathcal{U}}}^K(S^0)$. Since $H$ is an isotropy subgroup in $\mathcal{U}$, we see that $K/H$ embeds $K$-equivariantly in an admissible representation, and hence $\chi_K(K/H)$ is already well-defined in $\pi_{0,\mathcal{U}}^K(S^0)$ and is mapped to the same element in $\pi_{0,\hat{\mathcal{U}}}^K(S^0)$. Therefore we deduce that
\[
\Theta_K(1_H)=\chi_K(K/H).
\]
See also \cite[XIX, Section~1]{may1996equivariant} for the description in terms of equivariant Euler characteristics.
\end{proof}

\begin{lem} \label{lem: kernel of fixed point map by identity component}
    Let $K$ be a compact Lie group and denote its identity component by $K_0$. Suppose that we use the universes $\mathcal{U}_K$ and $\mathcal{U}_{\pi_0 K}$ for $K$- and $\pi_0 K$-equivariant stable homotopy theories, where the natural map $\mathcal{U}_K^{K_0}\rightarrow \mathcal{U}_{\pi_0 K}$ is an embedding. Then the kernel of the fixed point map
    \[
    \Phi^{K_0}:\pi_0 ^G(S^0) \rightarrow \pi_0 ^{\pi_0 G}(S^0)
    \]
    is generated by equivariant Euler characteristics $\chi_K(K/H)$ where $H$ runs over all conjugacy classes of closed subgroups of $K$ that does not contain $K_0$, has finite Weyl group in $K$, and arise as isotropy subgroups of finite-dimensional $K$-subrepresentations of $\mathcal{U}_K$. 

    Moreover, for every such subgroup $H$ with finite Weyl group,
    \[
    \deg \Phi^{K_0}\bigl(\chi_K(K/H)\bigr)
    =
    \chi\bigl((K/H)^{K_0}\bigr)
    =
    \begin{cases}
    [K:H], & K_0\subset H,\\
    0, & K_0\not\subset H,
    \end{cases}
    \]
    where $\deg$ denotes the ordinary degree after forgetting the residual $\pi_0K$-action.
\end{lem}
\begin{proof}
   Observe that, for any closed subgroup $H\subset K$, we have
\[
[\chi_K(K/H)^{K_0}]
=
[\chi_{K/K_0}((K/H)^{K_0})]
=
\begin{cases}
[\chi_{\pi_0 K}(\pi_0 K/\pi_0 H)]
& \text{if } K_0\subset H,\\
0
& \text{if } K_0\not\subset H.
\end{cases}
\]
    Since the map
    \[
    H\mapsto \pi_0 H:\left\{ H\in \mathcal{C}(K)\,\vert\,K_0 \subset H \right\} \rightarrow \mathcal{C}(\pi_0 K)
    \]
    is clearly injective, the lemma follows from \Cref{thm: tom dieck and eqv euler characteristic}.

    Moreover, more generally, for any closed subgroup $L\subset K$, taking
$L$-fixed points of an equivariant Euler characteristic gives
\[
\Phi^L\bigl(\chi_K(K/H)\bigr)
=
\chi_{W_KL}\bigl((K/H)^L\bigr).
\]
Hence, after forgetting the residual $W_KL$-action,
\[
\deg \Phi^L\bigl(\chi_K(K/H)\bigr)
=
\chi\bigl((K/H)^L\bigr).
\]
In particular, taking $L=K_0$, we obtain
\[
\deg \Phi^{K_0}\bigl(\chi_K(K/H)\bigr)
=
\begin{cases}
[K:H], & K_0\subset H,\\
0, & K_0\not\subset H.
\end{cases}
\]
\end{proof}

\begin{lem} \label{lem: classifying integer}
    Consider the equivariant Euler characteristic
    \[
    e :=\chi_{G_\mathrm{odd}}(G_\mathrm{odd}/\langle I \rangle) \in \pi_0^{G_\mathrm{odd}}(S^0).
    \]
    Then, for any $[f]\in \pi_0^{G_\mathrm{odd}}(S^0)$ satisfying $f^{S^1}\sim _{G_\mathrm{odd}/S^1} \mathrm{id}$, there exists a unique integer $n$ such that
    \[
    f \sim_{G_\mathrm{odd}} \mathrm{id}+ne.
    \]
    Furthermore, we have $\deg f^{\langle I \rangle} = 4n+1$.
\end{lem}
\begin{proof}
    Observe from \Cref{rem: S1 fixed point identity} that $[f]-[\mathrm{id}]$ is contained in the kernel of the fixed point map
    \[
    \Phi^{S^1}:\pi_0 ^{G_\mathrm{odd}}(S^0)\rightarrow \pi_0 ^{\pi_0 G_\mathrm{odd}}(S^0).
    \]
    Up to conjugacy, the list of all possible isotropy groups in the $G_\mathrm{odd}$-universe $\mathcal{U}_{G_\mathrm{odd}}$ is given by
    \[
    \{1\},\quad \langle i\tilde\tau\rangle,\quad \langle I \rangle,\quad S^1,\quad \langle S^1,I\rangle,\quad \langle S^1,\tilde\tau\rangle,\quad \mathrm{Pin}(2),\quad G_\mathrm{odd}.
    \]
    Among them, $\langle I \rangle$ is the only one that does not contain $S^1$ and has finite Weyl group. Hence, by \Cref{lem: kernel of fixed point map by identity component}, we see that $[f]-[\mathrm{id}]$ is a multiple of $e$, so we may write
    \[
    f \sim_{G_\mathrm{odd}} \mathrm{id}+ne,\quad n\in\Z.
    \]
Put $H=\langle I\rangle=\langle j\tilde\tau\rangle$.
Since $H$ has order two, its normalizer is the centralizer of $I$.
As $\tilde\tau$ is central and
$C_{\mathrm{Pin}(2)}(j)=\langle j\rangle$, we have
\[
N_{G_{\mathrm{odd}}}(H)
=\langle j,\tilde\tau\rangle
=\langle j\rangle \sqcup \langle j\rangle\tilde\tau.
\]
Thus $(G_{\mathrm{odd}}/H)^H=N_{G_{\mathrm{odd}}}(H)/H$.
Recall that taking $H$-fixed points of an equivariant Euler
characteristic gives the equivariant Euler characteristic of the
$H$-fixed set; after forgetting the remaining Weyl-group action, its
degree is the ordinary Euler characteristic. Hence
\[
\deg\Phi^H(e)
=\chi\bigl((G_{\mathrm{odd}}/H)^H\bigr)
=4.
\]
    Therefore we get $\deg f^{\langle I \rangle} = 1+n\cdot \deg \Phi^{\langle I \rangle}(e) = 4n+1$, as desired.
\end{proof}

\begin{rem} \label{rem: signed refinement of miyazawa}
    In the present spin setting, the Miyazawa invariant admits a canonical signed refinement. Indeed, $j$ commutes with $I=j\tilde\tau$ and satisfies $j^2=-1$. It therefore gives a complex structure on the $I$-fixed spinorial index, and hence a canonical orientation of its real determinant line. These orientations define the signed degree
\[
\mathfrak M_{\mathrm{sign}}(X,F)
:=
\deg BF_{\Sigma(F),\s,\tilde\tau}^{\langle I\rangle},
\qquad
|\mathfrak M_{\mathrm{sign}}(X,F)|=\mathfrak M(X,F).
\]
This sign is natural; it is gauge theoretically determined. Moreover, the condition
\[
\mathfrak M_{\mathrm{sign}}(X,F)\equiv1\pmod4
\]
shown in \Cref{lem: classifying integer} gives a simple numerical way to determine the sign. This shows that $\mathfrak{M}_\mathrm{sign}(X,F)$ is also independent of the spin lift $\tilde\tau$, as $\mathfrak{M}(X,F)$ already uniquely determines the sign of $\mathfrak{M}_\mathrm{sign}(X,F)$. In particular, the \emph{signed Miyazawa invariant} $\mathfrak{M}_\mathrm{sign}(X,F)$ is the unique integer congruent to $1$ modulo $4$ whose absolute value is $\mathfrak M(X,F)$.
\end{rem}

We now prove \Cref{thm: Miyazawa determines equivariant BF for S4} more generally for smoothly embedded $2$-spheres in smooth $\mathbb F_2$-homology $4$-spheres.

\begin{thm}\label{thm: Miyazawa determines equivariant BF}
Let $X$ be a smooth $\mathbb F_2$-homology $4$-sphere, and let $F$ be a smoothly embedded $2$-sphere in $X$. Equip the double cover $\Sigma(F)$ branched along $F$ with its unique spin structure $\s$, and let $\tilde\tau$ be a spin lift of the deck involution. Then the Miyazawa invariant $\mathfrak M(X,F)$ determines the $G_{\mathrm{odd}}$-equivariant stable homotopy class
\[
[BF_{\Sigma(F),\s,\tilde\tau}]\in\pi_0^{G_{\mathrm{odd}}}(S^0).
\]
\end{thm}


\begin{proof}
As discussed above, the $S^1$-fixed-point part of
$BF_{\Sigma(F),\s,\tilde\tau}$ is
$G_\mathrm{odd}/S^1$-equivariantly homotopic to
the identity. Hence, by
\Cref{lem: classifying integer}, there exists a
unique integer $n$ such that
\[
[BF_{\Sigma(F),\s,\tilde\tau}]
=[\mathrm{id}]+ne, \text{ and }
\deg BF_{\Sigma(F),\s,\tilde\tau}^{\langle I\rangle}
=4n+1.
\]
By the definition of the (generalized) Miyazawa
invariant,
\[
\mathfrak{M}(X,F)=|4n+1|.
\]
Since $4n+1\equiv 1\pmod 4$, its sign is uniquely
determined by its absolute value. Consequently,
\[
4n+1=
\begin{cases}
\mathfrak{M}(X,F),
&\mathfrak{M}(X,F)\equiv 1\pmod 4,\\
-\mathfrak{M}(X,F),
&\mathfrak{M}(X,F)\equiv 3\pmod 4.
\end{cases}
\]
This completes the proof. 
\end{proof}

Now we extract the (reduced) family Bauer--Furuta invariant $FBF(\tilde\tau)$ from equivariant Bauer--Furuta invariant $BF_{\Sigma(F),\s,\tilde\tau}$. We first give an algebraic model, i.e. a map
\[
T:\pi_0^{G_\mathrm{odd}}(S^0)\rightarrow \pi_1^{\mathrm{Pin}(2)}(S^0),
\]
and then show that it intertwines the two Bauer--Furuta invariants. 

Recall that the universes fixed above satisfy
\[
\operatorname{Res}^{G_{\mathrm{odd}}}_{\mathrm{Pin}(2)}
\mathcal{U}_{G_{\mathrm{odd}}}
\cong
\mathcal{U}_{\mathrm{Pin}(2)}
=
\R^\infty\oplus\wt\R^\infty\oplus\H^\infty.
\]
We fix such a $\mathrm{Pin}(2)$-equivariant identification
and use it for all fiberwise stabilizations below.
The construction uses the same associated sphere bundle
on the source and target, which is cancelled by fiberwise
desuspension; no trivialization of this bundle is required.
Compatibility with the index identifications used in the
family Bauer--Furuta invariant will be verified in the proof
of \Cref{lem: eqv to family BF}.

Given a $G_\mathrm{odd}$-equivariant stable map $f:S^0 \rightarrow S^0$, realize it as a concrete $G_\mathrm{odd}$-equivariant map $F:V^+ \rightarrow V^+$ for some $V\subset \mathcal{U}_{G_\mathrm{odd}}$, and form the $\mathrm{Pin}(2)$-equivariant vector bundle
\[
V_\tau := V \times [0,1] / (v,1) \sim (\tilde\tau (v),0) \rightarrow S^1.
\]
Since $F$ commutes with $\tilde\tau$, it induces a fiberwise endomorphism
\[
F_\tau : (V_\tau)^+ \rightarrow (V_\tau)^+,\quad (v,t) \mapsto (F(v),t),
\]
of $V^+$-bundle $(V_\tau)^+ \rightarrow S^1$. Thus, after fiberwise desuspending by the sphere bundle $V_\tau^+$, we get a parametrized endomorphism
\[
F_\tau:\underline{S^0}_{S^1} \rightarrow \underline{S^0}_{S^1},
\]
which, by composing the projection
\[
\Sigma ^\infty \underline{S^0}_{S^1} \to \Sigma ^\infty S^0, 
\]
defines an element 
\[
[F_\tau]\in [\Sigma^\infty S^1_+,\Sigma^\infty S^0]^{\mathrm{Pin}(2)} = \pi_0 ^{\mathrm{Pin}(2)}(S^0) \oplus \pi_1 ^{\mathrm{Pin}(2)}(S^0).
\]
We denote its $\pi_1 ^{\mathrm{Pin}(2)}(S^0)$ component by $T([f])$, i.e.
\[
[F_\tau] = ([F],T([f])) \in \pi_0 ^{\mathrm{Pin}(2)}(S^0) \oplus \pi_1 ^{\mathrm{Pin}(2)}(S^0).
\]

\begin{rem}
    Equivariant homotopies descend to fiberwise homotopies of the associated bundles. Stabilization by a representation $W$ adds the same associated sphere bundle to source and target, which cancels under fiberwise desuspension. The construction preserves addition of stable maps. Consequently, it defines a homomorphism 
    \[ 
    T:\pi_0^{G_{\mathrm{odd}}}(S^0) \longrightarrow\pi_1^{\operatorname{Pin}(2)}(S^0).
    \]
\end{rem}

\begin{lem} \label{lem: eqv to family BF}
    We have $T([BF_{\Sigma(F),\s,\tilde\tau}]) = [FBF(\tilde\tau)]$.
\end{lem}
\begin{proof}
Choose $\tau$-invariant auxiliary data and a
$G_{\mathrm{odd}}$-equivariant finite-dimensional approximation
to the monopole map on $(\Sigma(F),\s)$.
Using the equivariant index identifications fixed above
and stabilizing if necessary, represent the equivariant
Bauer--Furuta invariant by a map
\[
\Phi:V^+\longrightarrow V^+,
\qquad V\subset\mathcal{U}_{G_{\mathrm{odd}}}.
\]
Keep these data and this approximation constant over $[0,1]$.
Gluing the endpoints by $\tilde\tau$, with the convention
$(v,1)\sim(\tilde\tau v,0)$, gives a finite-dimensional
approximation for $\mathrm{Torus}(\tilde\tau)\to S^1$:
\[
\Phi_\tau:(V_\tau)^+\longrightarrow(V_\tau)^+.
\]

We check that the index identifications also agree.
The equivariant identification of the spinor summands
descends to a quaternionic bundle identification over $S^1$.
After stabilization, this agrees up to homotopy with the
index trivialization used in the family construction,
since quaternionic trivializations over $S^1$ are unique
up to homotopy; see \cite[Remark~12]{lin2025family}.
For the real summands, the identification is induced by
the normalized AHS operator $d^*\oplus d^+$, with
mean-zero functions in the target.
This operator is equivariant for the chosen auxiliary data,
so its mapping-torus construction is precisely the family
of operators used to identify the real source and target
in the family invariant.

Consequently, using the universe identification above,
fiberwise desuspension of $\Phi_\tau$ by the common sphere
bundle $(V_\tau)^+$ represents the full family
Bauer--Furuta invariant of $\mathrm{Torus}(\tilde\tau)$.
Under the splitting
\[
[\Sigma^\infty S^1_+,\Sigma^\infty S^0]^{\mathrm{Pin}(2)}
\cong
\pi_0^{\mathrm{Pin}(2)}(S^0)
\oplus
\pi_1^{\mathrm{Pin}(2)}(S^0),
\]
we therefore have
\[
[\Phi_\tau]
=
\bigl([BF_{\Sigma(F),\s}],[FBF(\tilde\tau)]\bigr).
\]
By the definition of $T$, the second component is also
$T([BF_{\Sigma(F),\s,\tilde\tau}])$, proving the assertion.
\end{proof}

We recall a relation between Becker--Gottlieb transfer maps and Pontryagin--Thom umkehr maps.

\begin{prop} \label{prop: transfer and umkehr}
    Let $K$ be a compact Lie group, $f:M\rightarrow B$ be a $K$-equivariant smooth submersion between closed $K$-manifolds, and $s$ be a $K$-equivariant smooth section of $T_v M:=\ker df$ that is transverse to the zero-section. Consider the inclusion map $i:s^{-1}(0)\hookrightarrow M$ and note that $s$ gives an isomorphism
    \[
    ds:\nu_{s^{-1}(0)} \xrightarrow\cong T_v M\vert_{s^{-1}(0)},
    \]
    which in turn induces a stable normal trivialization $ds:\nu_{s^{-1}(0)}-T_v M\vert_{s^{-1}(0)}\xrightarrow\simeq 0$. Suppose that, for each $b\in B$, the submanifold $f^{-1}(b) \subset M$ can be $\mathrm{Stab}_K(b)$-equivariantly embedded in some admissible $K$-representation. Then the $K$-equivariant Becker--Gottlieb transfer $\mathrm{Tr}^K_f$ is $K$-equivariantly stably homotopic to the following composition:
    \[
    \Sigma^\infty B_+ \xrightarrow{(f\vert_{s^{-1}(0)})_{ds}^!} \Sigma^\infty s^{-1}(0)_+ \xrightarrow{\Sigma^\infty i_+} \Sigma^\infty M_+.
    \]
\end{prop}
\begin{proof}
    This follows from the equivariant versions of the functoriality and transverse-pullback formulas for Pontryagin--Thom umkehr maps \cite[Part II, Propositions 12.7 and 12.11]{CrabbJames1998}, together with the Euler class description of the transfer \cite[Part II, Proposition 12.30]{CrabbJames1998}. Note that, since we are working with incomplete universes, we have to make sure that $M$ can be $K$-equivariantly embedded in $B\times V$ for some admissible $K$-representation $V$; this is simply a consequence of the given fiberwise admissibility.\footnote{This can be done via equivariant handle induction as $B$ admits a finite $K$-equivariant handle decomposition.}
\end{proof}

We also need a relation between Becker--Gottlieb transfer maps and associated bundle. Let $K$ and $H$ be compact Lie groups, $B$ be a closed smooth manifold with trivial $K$-action, and $Q\rightarrow B$ be a smooth principal right $H$-bundle. Let $X$ be a closed smooth manifold with commuting smooth $K$- and $H$-actions, and consider the bundle
\[
p:E = Q\times_H X \rightarrow B.
\]
Also, choose a finite-dimensional $K\times H$-representation $V$ and consider the \emph{associated (vector) bundle}
\[
\xi : Q\times_H V \rightarrow B;
\]
denote its fiberwise one-point compactification by $\xi^+$.

\begin{prop} \label{prop: associated bundle and transfer}
    $F:V^+ \rightarrow V^+$ represents $\chi_{K\times H}(X)$ and let $F_Q:\xi^+ \rightarrow \xi^+$ be the parametrized endomorphism (over $B$) defined by applying $F$ fiberwise. Fiberwise desuspension by $\xi$ then gives a parametrized stable endomorphism class
    \[
    [F_Q]\in [\underline{S^0}_B,\underline{S^0}_B]^K.
    \]
    Suppose that we are working with $K\times H$-universe $\mathcal{V}$ and $K$-universe $\mathcal{W}$ such that:
    \begin{itemize}
        \item For each $x\in X$, its orbit $(K\times H)/\mathrm{Stab}(x)$ embeds equivariantly in a finite-dimensional subrepresentation of $\mathcal{V}$;
        \item $\mathcal{W}$ contains $\mathrm{Res}_K^{K\times H}\mathcal{V}$.
    \end{itemize}
    Then, under the natural identification $[\underline{S^0}_B,\underline{S^0}_B]^K \cong [\Sigma^\infty B_+,\Sigma^\infty S^0]^K$, we have
    \[
    [F_Q] = c_E \circ \mathrm{Tr}_p^K,
    \]
    where $c_E:\Sigma^\infty E_+ \rightarrow S^0$ is the collapsing map.
\end{prop}
\begin{proof}
Write $\Gamma=K\times H$ and let $q\colon\Gamma\to K$ be the
projection. The orbit hypothesis and a finite $\Gamma$-CW
decomposition of $X$ imply that $\Sigma^\infty X_+$ is dualizable
in $\mathcal V$. Set
\[
\mathcal V'=\mathcal V\oplus q^*\mathcal W.
\]
Change of universe along $\mathcal V\hookrightarrow\mathcal V'$
preserves the equivariant Euler characteristic, since its left
derived functor is symmetric monoidal
\cite[Theorem~14.2.4]{may2006parametrized}. Thus $F$ also represents
$\chi_\Gamma(X)$ in $\mathcal V'$.

Give $Q$ the trivial $K$-action. Since $H$ acts freely on $Q$,
the stabilizer of each point of $Q$ is $K\times\{1\}$. Moreover,
the second hypothesis gives
\[
\mathrm{Res}^{\Gamma}_K\mathcal V'
=\mathrm{Res}^{\Gamma}_K\mathcal V\oplus\mathcal W
\cong_K\mathcal W.
\]
Hence the inclusion $q^*\mathcal W\hookrightarrow\mathcal V'$
satisfies the condition of
\Cref{rem: may sigurdsson lemma on incomplete universe} after
pullback to $Q$, as required in the proof of
\cite[Proposition~15.6.1]{may2006parametrized}.
The associated-bundle construction sends $F$ to $F_Q$.
Applying the transfer formula
\cite[Theorem~15.6.3]{may2006parametrized} to $\mathrm{id}_X$
and composing with collapse gives
$[F_Q]=c_E\circ\mathrm{Tr}_p^K$.
\end{proof}

\begin{rem} \label{rem: may sigurdsson lemma on incomplete universe}
Let $\Pi$ be a closed normal subgroup of $\Gamma$, let
$q\colon\Gamma\to G=\Gamma/\Pi$ be the quotient map, and let $Y$
be a $\Gamma$-space. Fix a $\Gamma$-equivariant isometric embedding
\[
i\colon q^*\mathcal U_G\hookrightarrow\mathcal U_\Gamma.
\]
Suppose that, for every $y\in Y$ and every closed subgroup
$\Lambda\subset\mathrm{Stab}_\Gamma(y)$ with
$\Lambda\cap\Pi=\{1\}$, the restrictions of $q^*\mathcal U_G$
and $\mathcal U_\Gamma$ to $\Lambda$ contain the same irreducible
real $\Lambda$-representations. Then change of universe along $i$
induces a symmetric monoidal equivalence on the homotopy categories
of $\Pi$-free $\Gamma$-spectra over $Y$. This follows by applying
the proof of \cite[Proposition~15.5.5]{may2006parametrized}
fiberwise.
\end{rem}

\begin{lem} \label{lem: computation of T}
    We have $T(\mathrm{id})=0$ and $T(e)=\alpha$.
\end{lem}
\begin{proof}
    The first statement is obvious: by tracking the definition of $T$, we see that it maps $\mathrm{id}$ to the homotopy class of the constant loop based at $\mathrm{id}$, and hence $T(\mathrm{id})=0$. We thus only have to compute $T(e)$.

    Recall that $e=\chi_{G_\mathrm{odd}}(G_\mathrm{odd}/\langle I \rangle)$; the quotient $G_\mathrm{odd}/\langle I \rangle$ can be identified with $G_\mathrm{odd}$-orbit $\mathrm{Pin}(2)\subset \H$:
    \[
    [g]\mapsto g\cdot 1:G_\mathrm{odd}/\langle I \rangle \hookrightarrow \mathrm{Pin}(2)\subset \H.
    \]
    Here, $I=j\tilde\tau$ acts on $\H$ by $I\cdot z = jzj^{-1}$, so we see that the stabilizer of $1\in \H$ is exactly $\langle I \rangle \subset G_\mathrm{odd}$. Consider the mapping torus family
    \[
    E := \mathrm{Pin}(2)\times [0,1]/(g,1)\sim (-gj,0) \xrightarrow{f} S^1.
    \]
    Since $e$ is the $G_\mathrm{odd}$-equivariant (and thus $\mathrm{Pin}(2)\times \langle \tilde\tau \rangle$-equivariant) Euler characteristic of this orbit, we see from \Cref{prop: associated bundle and transfer} under the setting $K=\mathrm{Pin}(2)$, $H=\langle \tilde\tau\rangle$, and $X=G_\mathrm{odd}/\langle I \rangle$, that $T(e)$ is the $\pi_1^{\mathrm{Pin}(2)}(S^0)$-component of the following family:
    \[
    \left[ \Sigma^\infty S^1_+ \xrightarrow{\mathrm{Tr}_f^{\mathrm{Pin}(2)}} \Sigma^\infty E_+ \xrightarrow{c_E} \Sigma^\infty S^0 \right] \in [\Sigma^\infty S^1_+,\Sigma^\infty S^0]^{\mathrm{Pin}(2)},
    \]
    where $c_E$ denotes the collapsing map for $E$.

    Consider the vertical tangent bundle $T_vE :=\ker df$ of $E\rightarrow S^1$. Then the real line bundle $T_v E/\mathrm{Pin}(2)\rightarrow S^1$ is nonorientable, so we may choose its smooth section $s$ which intersects the zero-section transversely at one point, say $p_0$. Then $s$ corresponds to a $\mathrm{Pin}(2)$-invariant vertical vector field, which we denote again by $s$, that is also transverse to the zero-section. The zero-locus of $s$ is given by
    \[
    s^{-1}(0) = f^{-1}(p_0) \cong \mathrm{Pin}(2).
    \]
    Observe that $f:E\rightarrow S^1$ satisfies the fiberwise admissibility condition in \Cref{prop: transfer and umkehr}, as $\mathrm{Pin}(2)$ acts trivially on $S^1$ and can be left-equivariantly embedded in $\H (\subset \mathcal{U}_{\mathrm{Pin}(2)})$. Hence we see that $T(e)$ is the $\pi_1^{\mathrm{Pin}(2)}(S^0)$-component of 
    \[
    c_E \circ \mathrm{Tr}^{\mathrm{Pin}(2)}_f \sim_{\mathrm{Pin}(2)} c_E \circ \Sigma^\infty i_+ \circ (f\vert_{f^{-1}(p_0)})^!_{ds} \sim_{\mathrm{Pin}(2)} c_{f^{-1}(p_0)}\circ (f\vert_{f^{-1}(p_0)})_{ds}^!,
    \]
    where $c_{f^{-1}(p_0)}:\Sigma^\infty f^{-1}(p_0)_+\rightarrow \Sigma^\infty S^0$ denotes the stable collapsing map and $i:f^{-1}(p_0)\hookrightarrow E$ denotes the inclusion. Since the map
    \[
    f\vert_{f^{-1}(p_0)}:f^{-1}(p_0)\rightarrow \{p_0\} \hookrightarrow S^1
    \]
    is clearly the collapsing map to $p_0 \in S^1$ and $ds$ is the unique (up to sign) left-invariant parallelization of $\mathrm{Pin}(2)$, we deduce from \Cref{rem: rephrasing of definition of alpha} that
    \[
    \begin{split}
    T(e) &\sim_{\mathrm{Pin}(2)} \pi_1^{\mathrm{Pin}(2)}(S^0)\text{-component of } c_{f^{-1}(p_0)}\circ (f\vert_{f^{-1}(p_0)})_{ds}^! \\
    &\sim_{\mathrm{Pin}(2)} \pi_1^{\mathrm{Pin}(2)}(S^0)\text{-component of } c_{\mathrm{Pin}(2)} \circ (\mathrm{Pin}(2)\rightarrow \ast)^!_{ds} \circ (\ast\hookrightarrow S^1)^! \\
    &\sim_{\mathrm{Pin}(2)} c_{\mathrm{Pin}(2)} \circ (\mathrm{Pin}(2)\rightarrow \ast)^!_{ds} \\
    &\sim_{\mathrm{Pin}(2)} \alpha,
    \end{split}
    \]
    as the map $\pi_1^{\mathrm{Pin}(2)}(S^0)\rightarrow  [\Sigma^\infty S^1_+,S^0]^{\mathrm{Pin}(2)}$ induced by precomposing with $\iota^!:\Sigma^\infty S^1_+\rightarrow \R^+$ is identified with the inclusion to the second summand, i.e.
    \[
    x\mapsto (0,x):\pi_1^{\mathrm{Pin}(2)}(S^0) \hookrightarrow \pi_0^{\mathrm{Pin}(2)}(S^0) \oplus \pi_1^{\mathrm{Pin}(2)}(S^0),
    \]
    under the canonical splitting $[\Sigma^\infty S^1_+,S^0]^{\mathrm{Pin}(2)}\cong \pi_0^{\mathrm{Pin}(2)}(S^0) \oplus \pi_1^{\mathrm{Pin}(2)}(S^0)$. The lemma follows.
\end{proof}

We are now ready to prove \Cref{thm: Lin to Miyazawa}, whose statement we recall.

\LinToMiyazawa*

\begin{proof}
Let $\tilde\tau$ be a spin lift of $\tau_F$. As discussed above, the $G_{\mathrm{odd}}$-equivariant Bauer--Furuta map satisfies $BF_{\Sigma(F),\s,\tilde\tau}^{S^1}\sim_{G_{\mathrm{odd}}/S^1}\mathrm{id}$. Hence, by \Cref{lem: classifying integer}, there is a unique integer $n$ such that
\[
BF_{\Sigma(F),\s,\tilde\tau}\sim_{G_{\mathrm{odd}}}\mathrm{id}+ne,
\qquad
\mathfrak M(X,F)=\left|\deg BF_{\Sigma(F),\s,\tilde\tau}^{\langle I\rangle}\right|=|4n+1|.
\]
By \Cref{lem: eqv to family BF,lem: computation of T}, we have
\[
[FBF(\tilde\tau)]
=T([BF_{\Sigma(F),\s,\tilde\tau}])
=T(\mathrm{id}+ne)
=n\alpha.
\]
Thus $\mathfrak{LM}_{\Sigma(F)}(\tau_F)=n\bmod 2$. Since
\[
\frac{\mathfrak M(X,F)^2-1}{8}
=\frac{(4n+1)^2-1}{8}
=2n^2+n
\equiv n\pmod 2,
\]
the theorem follows.
\end{proof}

\section{Real Bauer--Furuta invariants of Cappell--Shaneson homotopy $4$-spheres}\label{sec: CS degree computation}

In this section, we compute the Miyazawa invariant of special Cappell--Shaneson $2$-knots. After recalling their construction in \Cref{sec: CS example construction}, we use a surgery formula to reduce the computation to the real Bauer--Furuta degree of associated three-torus mapping tori. The explicit degree computation in \Cref{prop: CS negative degree} then shows that their Miyazawa invariant equals $3$, providing the examples needed to prove our main results.

\subsection{Cappell--Shaneson homotopy $4$-spheres and $2$-knots}\label{sec: CS example construction}

Let $T^3=\R^3/\Z^3$ be the standard three-dimensional torus, and choose a matrix $B\in SL_3(\Z)$ such that $\det(I-B)=\pm1$. Consider the smooth oriented $4$-manifold
\[
X_B=(\R\times T^3)/\langle\varphi_B\rangle,
\qquad
\varphi_B(t,x)=(t+1,Bx).
\]
Thus $X_B$ is the mapping torus of the diffeomorphism $x\mapsto Bx$ of $T^3$. It is straightforward to check that $X_B$ is a homology $S^1\times S^3$ and that $\pi_1(X_B)\cong\Z^3\rtimes_B\Z$. Since the torus diffeomorphism fixes $p=[0]\in T^3$, it determines an embedded loop
\[
\gamma_B=(\R\times\{p\})/\langle\varphi_B\rangle\subset X_B.
\]

There are two homotopy classes of framings of the oriented normal bundle of $\gamma_B$. Choose a framing $\mu$, which identifies a tubular neighborhood $\nu(\gamma_B)$ with $S^1\times D^3$. Surgery along $\gamma_B$ with this framing gives
\[
\mathcal{S}_{B,\mu}
=
\bigl(X_B\smallsetminus\operatorname{int}\nu(\gamma_B)\bigr)
\cup_\mu(D^2\times S^2),
\]
where the gluing along $S^1\times S^2$ is determined by $\mu$. It is straightforward to check that $\mathcal{S}_{B,\mu}$ is a homology $4$-sphere and that
\[
\pi_1(\mathcal{S}_{B,\mu})
\cong\operatorname{coker}(I-B)=0.
\]
Thus $\mathcal{S}_{B,\mu}$ is a smooth homotopy $4$-sphere. We call it the \emph{Cappell--Shaneson homotopy $4$-sphere} associated with $B$ and the framing $\mu$; see \cite{cappellshaneson76} for the original construction.

Consider the matrices
\[
A_m=
\begin{pmatrix}
0 & 1 & 0\\
0 & 1 & 1\\
1 & 0 & m+1
\end{pmatrix}
\in SL_3(\Z),
\qquad m\in\Z.
\]
For these matrices, Akbulut proved standardness for one framing
\cite[Theorem~1]{akbulut2010cappell}, and Gompf's argument also establishes standardness for the other framing
\cite[Example~3.1(a)]{gompf2010more}. Thus we have the following.

\begin{thm}[{\cite{akbulut2010cappell,gompf2010more}}]\label{thm: akbulut standardness}
The Cappell--Shaneson homotopy $4$-sphere $\mathcal{S}_{A_m,\mu}$ is diffeomorphic to $S^4$ for every integer $m$ and every framing $\mu$ of $\gamma_{A_m}\subset X_{A_m}$.
\end{thm}

Cappell and Shaneson introduced a construction of fibered $2$-knots in smooth homotopy $4$-spheres whose fibers are punctured $3$-tori \cite{cs1976inequivalent}. We recall this construction, following \cite[p.~2]{Aitchison-Rubinstein:1984}, and describe the deck involution of the branched double cover explicitly. Let $B\in SL_3(\Z)$ satisfy $\det(I-B^2)=\pm1$. Since $\det(I-B^2)=\det(I-B)\det(I+B)$, both $B$ and $B^2$ satisfy the Cappell--Shaneson condition. Consider the three-torus bundle $X_{B^2}$ and the embedded loop $\gamma_{B^2}\subset X_{B^2}$. After rescaling the $\R$-coordinate, we identify
\[
X_{B^2}=(\R\times T^3)/\langle\varphi_B^2\rangle.
\]
The natural double cover $X_{B^2}\to X_B$ has deck involution
\[
\tau_B\colon X_{B^2}\longrightarrow X_{B^2};\qquad [t,x]\longmapsto[t+1,Bx].
\]

Among the two homotopy classes of framings of the oriented normal bundle of $\gamma_{B^2}$, exactly one admits a $\Z_2$-equivariant parametrization
\[
\nu(\gamma_{B^2})\cong S^1\times D^3,
\]
where the $\Z_2$-action on the right is given by $(z,v)\mapsto(-z,v)$. This framing class has two equivariant refinements up to equivariant homotopy, corresponding to the two homotopy classes of framings of the quotient loop $\gamma_B$. We choose either refinement and denote it by $\mu_0$.

Since $\mu_0$ is $\Z_2$-equivariant, the deck involution extends over the surgery to a $\Z_2$-action on the Cappell--Shaneson homotopy $4$-sphere $\mathcal{S}_{B^2,\mu_0}$. Let $\hat\mu_0$ denote the induced framing of the quotient loop $\gamma_B=\gamma_{B^2}/\Z_2$ in $X_B=X_{B^2}/\Z_2$. Then
\[
\mathcal{S}_{B^2,\mu_0}^{\Z_2}\cong S^2
\qquad\text{ and }\qquad
\mathcal{S}_{B^2,\mu_0}/\Z_2\cong\mathcal{S}_{B,\hat\mu_0}.
\]
The quotient map
\[
q_B\colon\mathcal{S}_{B^2,\mu_0}\longrightarrow\mathcal{S}_{B,\hat\mu_0}
\]
is a branched double covering whose branch locus is the smoothly embedded $2$-sphere given by the image of $\mathcal{S}_{B^2,\mu_0}^{\Z_2}$. We denote this branch locus by $K_B$, which is a $2$-knot in the homotopy $4$-sphere $\mathcal{S}_{B,\hat\mu_0}$.

\begin{defn}
Suppose that $\mathcal{S}_{B^2,\mu_0}$ and $\mathcal{S}_{B,\hat\mu_0}$ are both diffeomorphic to $S^4$. An orientation-preserving identification of $\mathcal{S}_{B,\hat\mu_0}$ with $S^4$ then determines a $2$-knot $K_B\subset S^4$, well-defined up to orientation-preserving self-diffeomorphism of $S^4$. In this paper, we call such a knot a \emph{Cappell--Shaneson $2$-knot}. We call it \emph{special} if $\det(I+B)=-1$.
\end{defn}

The significance of the condition $\det(I+B)=-1$ for our purposes is that the degree computation in \Cref{prop: CS negative degree} gives $\mathfrak M(S^4,K_B)=3$ for every special Cappell--Shaneson $2$-knot, so \Cref{thm: Lin to Miyazawa} detects its deck involution.

This condition also has a geometric interpretation. The two choices of framing on $\gamma_B$ give rise to two knot pairs with the same exterior, related by Gluck reconstruction. The knot $K_B$ is called \emph{reflexive} if these pairs are homeomorphic, allowing changes of orientation. Hillman and Wilson showed that, under the normalization $\det(I-B)=-1$, the knot $K_B$ is reflexive if and only if $\operatorname{tr}(B)=-1$ \cite{Hillman-Wilson:1989}. To compare this criterion with ours, replace $B$ by $B^{-1}$ if necessary. This reverses the direction of the base circle and changes the sign of $\det(I-B)$ while preserving $\det(I+B)$, so we may assume that $\det(I-B)=-1$. Write $a=\operatorname{tr}(B)$. Since the characteristic polynomial $p_B(t)=t^3-at^2+bt-1$ satisfies $p_B(1)=-1$, we have $b=a-1$. It follows that
\[
\det(I+B)=-p_B(-1)=2a+1.
\]
Thus $K_B$ is special if and only if it is reflexive.

We now show that $K_{A_{-3}}$ is a special Cappell--Shaneson $2$-knot. For a description of the corresponding branched covering in terms of involutions, see \cite[p.~23]{Aitchison-Rubinstein:1984}.

\begin{lem}\label{lem: cappell-shaneson 2-knots exist}
The knot $K_{A_{-3}}$ is a special Cappell--Shaneson $2$-knot.
\end{lem}

\begin{proof}
For every $m\in\Z$, we have
\[
\det(I-A_m)=-1,\qquad \det(I+A_m)=2m+5.
\]
Since $\det(I-A_m^2)=\det(I-A_m)\det(I+A_m)$, the condition $\det(I-A_m^2)=\pm1$ holds precisely when $m=-2$ or $m=-3$. Of these two values, only $m=-3$ satisfies $\det(I+A_m)=-1$.

It remains to show that the quotient and the branched double cover are both diffeomorphic to $S^4$. For the quotient $\mathcal{S}_{A_{-3},\hat\mu_0}$, this follows directly from \Cref{thm: akbulut standardness}. For the cover, a direct calculation gives
\[
A_{-3}^{-2}=Q A_4 Q^{-1},\qquad
Q=\begin{pmatrix}
0&1&2\\
0&0&-1\\
-1&2&1
\end{pmatrix}\in SL_3(\Z).
\]
Reversing the base circle replaces the monodromy by its inverse, and conjugating the monodromy amounts to changing coordinates on $T^3$. Both operations preserve the mapping torus together with its surgery circle up to diffeomorphism. Thus surgery on $\gamma_{A_{-3}^2}$ with framing $\mu_0$ is diffeomorphic to surgery on $\gamma_{A_4}$ with the corresponding framing. Since \Cref{thm: akbulut standardness} applies to both framings, the branched double cover $\mathcal{S}_{A_{-3}^2,\mu_0}$ is also diffeomorphic to $S^4$. This proves the lemma.
\end{proof}

\begin{rem}
To obtain a self-diffeomorphism of the standard $S^4$ from the deck involution of $q_B$, we need only the branched double cover $\mathcal{S}_{B^2,\mu_0}$ to be diffeomorphic to $S^4$. We require the quotient $\mathcal{S}_{B,\hat\mu_0}$ to be diffeomorphic to $S^4$ in the definition above so that the branch locus $K_B$ is a $2$-knot in the standard $S^4$. Our computation does not require this assumption on the quotient.
\end{rem}

\subsection{Invariants of homology $S^1\times S^3$}

Since $K_B$ is the belt sphere of the surgery on $\gamma_B$, surgery on $K_B$ recovers the mapping torus $X_B$. The surgery formula of Miyazawa--Park--Taniguchi \cite[Proposition~7.3]{miyazawa2025satellite} therefore identifies the Miyazawa invariant of $K_B$ with the absolute degree of the real Bauer--Furuta invariant associated with the double cover $X_{B^2}\to X_B$. We recall the definition of this degree and explain how to compute it by counting framed real monopoles with signs.

Let $X$ be a closed oriented integral homology $S^1\times S^3$, and let $\widetilde X\to X$ be its unique nontrivial double cover. By \cite[Lemma~7.1]{miyazawa2025satellite}, $\widetilde X$ admits a real spin structure. With respect to the induced real $\mathrm{spin}^c$ structure, its real Bauer--Furuta invariant is represented by a stable map
\[
BF^R_{\widetilde X}\colon S^{-\sigma(\widetilde X)/16}=S^0\longrightarrow S^{b^+(\widetilde X)-b^+(X)}=S^0.
\]
Following \cite[Section~7]{miyazawa2025satellite}, we set
\[
|\deg(X)|:=\bigl|\deg BF^R_{\widetilde X}\bigr|\in\mathbb Z_{\geq0}.
\]
When the framed real moduli space is cut out transversely, $|\deg(X)|$ is the absolute value of its signed count. Here we quotient only by the identity component of the real gauge group, so the residual $\{\pm1\}$-action is retained. By \cite[Lemma~7.2]{miyazawa2025satellite}, $|\deg(X)|$ is invariant under orientation-preserving diffeomorphisms and equals $1$ if $X$ admits a metric of positive scalar curvature.

We describe the framed real monopole moduli space. Let $\tau$ denote the deck involution of $\widetilde X\to X$ and $I$ the real structure, which also acts on connections and sections. Fix a $\tau$-invariant metric on $\widetilde X$. The real gauge group and its identity component are
\[
\mathcal G^\R_{\wt X}=\{u\colon\widetilde X\to S^1\mid u\circ\tau=\overline u\}
\qquad\text{ and }\qquad
 (\mathcal G^\R_{\wt X})_0=\{e^{i\chi}\mid\chi\in C^\infty(\widetilde X;\R),\ \chi\circ\tau=-\chi\}
\]
with the appropriate Sobolev completions.  Define the framed real monopole moduli space by
\[
\widetilde{\mathcal M}^I
=
\left\{
(A,\psi)\ \middle|\
\begin{gathered}
I(A,\psi)=(A,\psi),\\
D_A^+\psi=0,\quad F_A^+=q(\psi)
\end{gathered}
\right\}\big/\mathcal (\mathcal{G}^\R_{\wt X})_0,
\]
where $A$ is a $\mathrm{spin}^c$ connection and $\psi$ is a positive spinor; note that this includes reducible solutions as well. Following the computations in \cite[Lemma 2.10]{nakamura2015mathrm}, we get
\[
\mathcal G^\R_{\wt X}/( \mathcal G^\R_{\wt X})_0 \cong \pi_0 \mathcal G^\R_{\wt X}\simeq H^1({X};l)\cong \Z_2
\]
where $l$ denotes the $\Z$-coefficient local system induced by the double covering $\wt{X}\rightarrow X$, which implies that the framed moduli space carries a residual $\{\pm1\}$-action. Since the mod 2 Gysin sequence shows that $\widetilde X$ is a rational homology $S^1\times S^3$, and transfer gives $H^1(\widetilde X;\R)^-=0$, we see that its formal dimension is 0.

\subsection{The geometry of the three-torus and explicit spinor conventions}
\label{subsec:cs-explicit-geometry}

Throughout the computation, we assume that
\[
B\in SL_3(\Z),\qquad \det(I-B^2)=\pm1,\qquad \det(I+B)=-1.
\]
The characteristic polynomial $p_B(t)=\det(tI-B)$ satisfies $p_B(-1)=1$ and $p_B(0)=-1$. Since $p_B$ is a monic cubic, it has one root in each of the intervals $(-\infty,-1)$, $(-1,0)$, and $(0,\infty)$. The positive root is not $1$, since $\det(I-B^2)\neq0$. Thus we may write the eigenvalues of $B$ as
\[
-e^{-k_1},\qquad -e^{-k_2},\qquad e^{-k_3},
\]
where $k_1<0<k_2$, $k_3\neq0$, and $k_1+k_2+k_3=0$ since $\det B=1$.

Choose an eigenbasis in this order, normalized to have determinant one, and let $y=(y_1,y_2,y_3)$ be the corresponding coordinates on $\R^3$. In these coordinates, the original lattice becomes a $B$-invariant lattice $\Lambda$ of covolume one, and we write
\[
Y=\R^3/\Lambda,\qquad K=\operatorname{diag}(k_1,k_2,k_3),\qquad R=\operatorname{diag}(-1,-1,1),\qquad B=e^{-K}R.
\]
For $s\in\R$, define a coframe, its dual frame, and a flat metric on $Y$ by
\[
\theta_s^i=e^{k_i s}dy_i,\qquad E_i(s)=e^{-k_i s}\partial_{y_i},\qquad g_s=\sum_i\theta_s^i\otimes\theta_s^i.
\]
The frame $(E_1(s),E_2(s),E_3(s))$ is oriented, orthonormal, and parallel with respect to $g_s$. Since $k_1+k_2+k_3=0$, we have
\[
d\operatorname{vol}_{g_s}=dy_1\wedge dy_2\wedge dy_3
\qquad\text{ and }\qquad
\operatorname{vol}(Y,g_s)=1.
\]
Finally, $B^\ast\theta_{s+1}^i=R_{ii}\theta_s^i$, so $B^\ast g_{s+1}=g_s$.

The parallel frame $(E_1(s),E_2(s),E_3(s))$ is invariant under translations on $Y$ and thus trivializes the oriented orthonormal frame bundle:
\[
\mathrm{Fr}(Y,g_s)\cong Y\times SO(3).
\]
The corresponding spin structure is given by the double cover
\[
Y\times\mathrm{Sp}(1)\longrightarrow Y\times SO(3);\qquad (y,q)\longmapsto(y,\operatorname{Ad}_q),
\]
where $\operatorname{Ad}_q(v)=qvq^{-1}$ on $\operatorname{Im}\H$. For the Clifford conventions below, we identify the oriented frame $$(E_1(s),E_2(s),E_3(s))$$ with the ordered orthonormal basis $(\mathbf k,\mathbf j,-\mathbf i)$ of $\operatorname{Im}\H$. This is the \emph{Lie spin structure} on $Y$, the unique spin structure to which the translation action of $Y$ lifts.

In this trivialization, the spinor bundle is
\[
S_{Y,s}=Y\times\H;\qquad (\psi_1,\psi_2)\leftrightarrow\psi_1+\mathbf j\psi_2,
\]
where $\psi_1,\psi_2\in\C$ and complex scalars act on the right. In particular, the complex structure is $J(\psi)=\psi\mathbf i$, so
\[
J(\psi_1+\mathbf j\psi_2)=\mathbf i\psi_1-\mathbf k\psi_2.
\]
Thus a spinor is a $\Lambda$-periodic quaternion-valued function on $\R^3$. Clifford multiplication, in contrast, acts on the left. Setting
\[
\varepsilon_1=-\mathbf k,\qquad \varepsilon_2=-\mathbf j,\qquad \varepsilon_3=\mathbf i,
\]
we have
\[
\rho_s(E_i(s))\psi=\varepsilon_i\psi
\qquad\text{ and }\qquad
\rho_s\left(\sum_i v_i\theta_s^i\right)\psi=\sum_i v_i\varepsilon_i\psi.
\]
The complex-antilinear quaternionic structure is given by right multiplication by $\mathbf j$:
\[
j(\psi)=\psi\mathbf j,\qquad j^2=-1,\qquad j\rho_s=\rho_sj.
\]
Since the frame is parallel, the spin Dirac operator for $g_s$ is
\[
D_s\psi=\sum_i\varepsilon_i e^{-k_i s}\partial_{y_i}\psi.
\]
Although the coordinates $y_i$ on the universal cover $\R^3$ do not descend to real-valued functions on $Y$, the translation-invariant vector fields $\partial_{y_i}$ do.

\begin{lem}\label{lem: constant C}
For every $s\in\R$, we have
\[
\ker D_s=\{\text{constant spinors}\}\cong\H.
\]
Moreover, there is a constant $C>0$ such that
\begin{equation}\label{eq:cs-detailed-fiber-gap}
\|\psi\|_{L^2(g_s)}\leq C\|D_s\psi\|_{L^2(g_s)}
\end{equation}
for every $s\in[0,2]$ and every smooth spinor $\psi$ that is $L^2(g_s)$-orthogonal to $\ker D_s$.
\end{lem}

\begin{proof}
The Clifford relations give
\[
D_s^2=-\sum_i e^{-2k_i s}\partial_{y_i}^2.
\]
Let $\Lambda^\ast=\{\xi\in\R^3\mid\langle\xi,\Lambda\rangle\subset\Z\}$ be the dual lattice. The functions $e_\xi(y)=e^{2\pi i\langle\xi,y\rangle}$, for $\xi\in\Lambda^\ast$, form an orthonormal basis of $L^2(Y;\C)$. Thus every spinor has a Fourier expansion $\psi=\sum_{\xi\in\Lambda^\ast}e_\xi v_\xi$ with $v_\xi\in\H$. Writing $\langle\xi,y\rangle=\sum_i\xi_i y_i$, we obtain
\[
D_s^2 e_\xi=\lambda_s(\xi)e_\xi,\qquad \lambda_s(\xi)=4\pi^2\sum_i e^{-2k_i s}\xi_i^2.
\]
Since $\lambda_s(\xi)=0$ precisely when $\xi=0$, and $D_s$ annihilates constant spinors, the first statement follows.

For $s\in[0,2]$, the weights $e^{-2k_i s}$ have a common positive lower bound. Since the nonzero vectors of $\Lambda^\ast$ have lengths bounded away from zero, there is a constant $\lambda_0>0$ such that $\lambda_s(\xi)\geq\lambda_0$ for every $s\in[0,2]$ and every nonzero $\xi\in\Lambda^\ast$. If $\psi$ is orthogonal to the constant spinors, then $v_0=0$. By formal self-adjointness of $D_s$ and Parseval's identity,
\[
\|D_s\psi\|_{L^2(g_s)}^2
=\langle D_s^2\psi,\psi\rangle_{L^2(g_s)}
=\sum_{\xi\ne0}\lambda_s(\xi)|v_\xi|^2
\geq\lambda_0\|\psi\|_{L^2(g_s)}^2.
\]
Taking $C=\lambda_0^{-1/2}$ proves the estimate.
\end{proof}

\subsection{Computations for mapping tori}

Consider the mapping torus
\[
X:=X_{B^2}=(\R\times Y)/\langle\Phi^2\rangle,
\qquad \Phi(s,y)=(s+1,By).
\]
Then $\Phi^2(s,y)=(s+2,B^2y)$, and $\Phi$ induces the deck involution $\tau_B$ of the double cover $X\to X_B$, given by $\tau_B([s,y])=[s+1,By]$. Fix $T>0$ and equip $\R\times Y$ with the metric
\[
g_T=T^2ds^2+g_s.
\]
Set
\[
\theta^0=T\,ds,\qquad \theta^i=\theta_s^i;
\qquad e_0=T^{-1}\partial_s,\qquad e_i=e^{-k_i s}\partial_{y_i}
\qquad (1\leq i\leq 3).
\]
We use the orientation for which $\theta^0\wedge\theta^1\wedge\theta^2\wedge\theta^3$ is positive. Since $B^\ast g_{s+1}=g_s$, the metric $g_T$ is $\Phi$-invariant and hence descends to a $\tau_B$-invariant metric on $X$. Moreover, $\Phi^\ast\theta^0=\theta^0$ and $\Phi^\ast\theta^i=R_{ii}\theta^i$, so $\Phi^2$ preserves the orthonormal coframe $(\theta^0,\theta^1,\theta^2,\theta^3)$. Thus this coframe and its dual frame $(e_0,e_1,e_2,e_3)$ both descend to $X$.

\begin{lem}\label{lem: scal}
The scalar curvature of $\R\times Y$ with respect to $g_T$ is
\[
\operatorname{scal}(g_T)=-\frac{1}{T^2}\sum_i k_i^2.
\]
\end{lem}

\begin{proof}
A direct computation gives $[e_0,e_i]=-\frac{k_i}{T}e_i$ and $[e_i,e_j]=0$ for $1\leq i,j\leq 3$. The Koszul formula then gives
\[
\nabla_{e_i}e_0=\frac{k_i}{T}e_i,\qquad
\nabla_{e_i}e_i=-\frac{k_i}{T}e_0
\qquad (1\leq i\leq 3),
\]
while $\nabla_{e_0}e_a=0$ for $0\leq a\leq 3$ and $\nabla_{e_i}e_j=0$ for distinct $i,j\in\{1,2,3\}$. Thus the sectional curvatures of the planes spanned by $e_0,e_i$ and by $e_i,e_j$ ($i<j$) are $-k_i^2/T^2$ and $-k_i k_j/T^2$, respectively. Since $\sum_i k_i=0$, we obtain
\[
\operatorname{scal}(g_T)
=-\frac{2}{T^2}\left(\sum_i k_i^2+\sum_{i<j}k_i k_j\right)
=-\frac{1}{T^2}\left(\left(\sum_i k_i\right)^2+\sum_i k_i^2\right)
=-\frac{1}{T^2}\sum_i k_i^2,
\]
which completes the proof.
\end{proof}

The global oriented orthonormal frame $(e_0,e_1,e_2,e_3)$ induces a spin structure $\s_B$ on $X$, given in this trivialization by the double cover $X\times\mathrm{Spin}(4)\to X\times SO(4)$. The corresponding half-spin bundles are
\[
S_B^+=X\times\H,\qquad S_B^-=X\times\H.
\]
We identify $\mathrm{Spin}(4)$ with $\mathrm{Sp}(1)_+\times\mathrm{Sp}(1)_-$, with the two factors acting by left multiplication on the fibers of $S_B^+$ and $S_B^-$, respectively. The subgroup fixing $e_0$ is the diagonal $\mathrm{Sp}(1)$, and both half-spin representations restrict to the spin representation used on the torus fibers. On both bundles, we retain the complex structure $J(\psi)=\psi\mathbf i$ and the quaternionic structure $j(\psi)=\psi\mathbf j$.

For a $1$-form $v=v_0\theta^0+\sum_{i=1}^3 v_i\theta^i$, define Clifford multiplication $c^\pm(v)\colon S_B^\pm\to S_B^\mp$ by
\[
c^\pm(v)\psi=\left(\pm v_0+\sum_{i=1}^3 v_i\varepsilon_i\right)\psi.
\]
The induced Clifford multiplication on $2$-forms is given by
\[
\begin{aligned}
\rho^+(v\wedge w)&=\frac{1}{2}\left(c^-(v)c^+(w)-c^-(w)c^+(v)\right),\\
\rho^-(v\wedge w)&=\frac{1}{2}\left(c^+(v)c^-(w)-c^+(w)c^-(v)\right).
\end{aligned}
\]
Let $A$ be a $U(1)$-connection on the determinant line bundle $\det_{\C}S_B^+$. Together with the Levi--Civita connection on $X$, it induces a connection $\nabla^A$ on $S_B^+$. The corresponding positive Dirac operator is
\[
D_A^+\psi=\sum_{i=0}^3 c^+(\theta^i)\nabla_{e_i}^A\psi.
\]
Let $\Lambda^+$ denote the real rank-three bundle of self-dual $2$-forms. The quadratic map $q\colon S_B^+\to i\Lambda^+$ is characterized by
\[
\rho^+(q(\psi))=(\psi\psi^\ast)_0
=\psi\psi^\ast-\frac{|\psi|^2}{2}\,\mathrm{id}_{S_B^+}.
\]
With these conventions, the Seiberg--Witten equations on $X$ take the following form.

\begin{lem}\label{lem:cs-quaternionic-SW}
Let $A_0$ be the trivial connection on $\det_{\C}S_B^+$, and write
\[
A=A_0+2ia,\qquad a=f\theta^0+\sum_{i=1}^3 u_i\theta^i.
\]
Then the Seiberg--Witten equations
\[
D_A^+\psi=0,\qquad F_A^+=q(\psi)
\]
are equivalent to
\begin{equation}\label{eq:cs-quaternionic-SW}
\begin{aligned}
\frac{1}{T}\partial_s\psi+D_s\psi
+\left(f+\sum_{i=1}^3 u_i\varepsilon_i\right)\psi\mathbf i &=0,\\
\sum_{\mathrm{cyc}}
\left(
\frac{\partial_s u_i+k_i u_i}{T}
-e_i(f)+e_j(u_k)-e_k(u_j)
\right)\varepsilon_i
&=\frac{1}{4}\psi\mathbf i\overline{\psi}.
\end{aligned}
\end{equation}
Here $\sum_{\mathrm{cyc}}$ denotes the sum over the cyclic permutations $(i,j,k)$ of $(1,2,3)$, and $\overline{\psi}$ denotes quaternionic conjugation.
\end{lem}

\begin{proof}
Since $A=A_0+2ia$ is a connection on the determinant line, the induced spinor connections satisfy
\[
\nabla^A_Z\psi=\nabla^{A_0}_Z\psi+a(Z)\psi\mathbf i.
\]
In the trivialization induced by the orthonormal frame $(e_0,e_1,e_2,e_3)$, the positive spin connection is given by
\[
\nabla^{A_0}_Z\psi
=Z(\psi)+\frac{1}{4}\sum_{r,t=0}^3
\langle\nabla_Z e_r,e_t\rangle
c^-(\theta^r)c^+(\theta^t)\psi.
\]
Using the computation in the proof of \Cref{lem: scal}, we obtain
\[
\nabla_{e_0}^{A_0}\psi=\frac{1}{T}\partial_s\psi,\qquad
\nabla_{e_i}^{A_0}\psi=e_i(\psi)-\frac{k_i}{2T}\varepsilon_i\psi
\qquad (1\leq i\leq 3).
\]
It follows that
\[
\begin{aligned}
D_{A_0}^+\psi
&=\sum_{r=0}^3 c^+(\theta^r)\nabla_{e_r}^{A_0}\psi\\
&=\frac{1}{T}\partial_s\psi
+\sum_{i=1}^3\varepsilon_i
\left(e_i(\psi)-\frac{k_i}{2T}\varepsilon_i\psi\right)
=\frac{1}{T}\partial_s\psi+D_s\psi,
\end{aligned}
\]
where the last equality uses $\varepsilon_i^2=-1$ and $\sum_i k_i=0$. Hence
\[
D_A^+\psi
=D_{A_0}^+\psi+c^+(a)\psi\mathbf i
=\frac{1}{T}\partial_s\psi+D_s\psi
+\left(f+\sum_{i=1}^3 u_i\varepsilon_i\right)\psi\mathbf i,
\]
which gives the first equation.

For the curvature equation, set
\[
\omega_i=\theta^0\wedge\theta^i+\theta^j\wedge\theta^k
\]
for each cyclic permutation $(i,j,k)$ of $(1,2,3)$. Since
\[
d\theta^i=\frac{k_i}{T}\theta^0\wedge\theta^i,\qquad F_A=2i\,da,
\]
we obtain
\[
F_A^+=i\sum_{\mathrm{cyc}}
\left(
\frac{\partial_s u_i+k_i u_i}{T}
-e_i(f)+e_j(u_k)-e_k(u_j)
\right)\omega_i.
\]
The Clifford action satisfies $\rho^+(\omega_i)\varphi=-2\varepsilon_i\varphi$. Since complex scalars act on the right, this gives
\[
\rho^+(F_A^+)\varphi
=-2\sum_{\mathrm{cyc}}
\left(
\frac{\partial_s u_i+k_i u_i}{T}
-e_i(f)+e_j(u_k)-e_k(u_j)
\right)\varepsilon_i\varphi\mathbf i.
\]
On the other hand, the Hermitian pairing on $S_B^+$ is given by
\[
\langle\psi,\varphi\rangle
=\frac{1}{2}\left(\overline{\psi}\varphi
-\mathbf i\overline{\psi}\varphi\mathbf i\right).
\]
Since $\psi\overline{\psi}=|\psi|^2$, we obtain
\[
(\psi\psi^\ast)_0\varphi
=\psi\langle\psi,\varphi\rangle-\frac{|\psi|^2}{2}\varphi
=-\frac{1}{2}\psi\mathbf i\overline{\psi}\varphi\mathbf i.
\]
Comparing the two expressions and using $\rho^+(q(\psi))=(\psi\psi^\ast)_0$ gives the second equation, which completes the proof.
\end{proof}

We now describe the real structure and the conditions it imposes on solutions to the monopole equations.

\begin{lem}\label{lem:cs-real-conditions}
On each half-spin bundle, the map
\[
I_B\colon S_B^\pm\longrightarrow S_B^\pm;
\qquad [s,y,\psi]\longmapsto[s+1,By,\mathbf i\psi\mathbf j]
\]
is a complex-antilinear involution that covers $\tau_B$ and is compatible with Clifford multiplication. With the notation of \Cref{lem:cs-quaternionic-SW}, a configuration $(A,\psi)$ is \emph{Real} (i.e., $I_B$-invariant) if and only if
\begin{equation}\label{eq:cs-detailed-spinor-boundary}
\psi(s+1,By)=\mathbf i\psi(s,y)\mathbf j
\end{equation}
and
\begin{equation}\label{eq:cs-detailed-form-boundary}
f(s+1,By)=-f(s,y),\qquad
u_i(s+1,By)=-R_{ii}u_i(s,y)
\qquad (1\leq i\leq 3).
\end{equation}
\end{lem}

\begin{proof}
In the product spin trivializations, the identification over $(s,y)\sim(s+2,B^2y)$ acts trivially on the spinor fibers. Thus $I_B$ is well-defined and covers $\tau_B$. Left multiplication by $\mathbf i$ is complex-linear, while right multiplication by $\mathbf j$ is complex-antilinear. Moreover, $I_B^2=\mathrm{id}$ since $\mathbf i^2=\mathbf j^2=-1$.

To check compatibility with Clifford multiplication, it suffices to consider the coframe $(\theta^0,\theta^1,\theta^2,\theta^3)$. Since $\tau_B^\ast\theta^0=\theta^0$ and $c^\pm(\theta^0)=\pm\mathrm{id}$ in the chosen trivializations, the claim holds for $\theta^0$. For $1\leq i\leq 3$, we have $\tau_B^\ast\theta^i=R_{ii}\theta^i$ and $\mathbf i\varepsilon_i\mathbf i^{-1}=R_{ii}\varepsilon_i$. Hence
\[
I_B(c^\pm(\theta^i)\psi)
=\mathbf i\varepsilon_i\psi\mathbf j
=(\mathbf i\varepsilon_i\mathbf i^{-1})I_B(\psi)
=R_{ii}c^\pm(\theta^i)I_B(\psi).
\]
Thus $I_B$ is compatible with Clifford multiplication by $1$-forms and hence with the induced Clifford multiplication by $2$-forms.

The spinor condition \eqref{eq:cs-detailed-spinor-boundary} follows directly from the definition of $I_B$. The induced real structure preserves the trivial determinant connection $A_0$, so $A=A_0+2ia$ is Real precisely when $\tau_B^\ast a=-a$. Since $\tau_B^\ast\theta^0=\theta^0$ and $\tau_B^\ast\theta^i=R_{ii}\theta^i$, this is equivalent to \eqref{eq:cs-detailed-form-boundary}, which completes the proof.
\end{proof}

As in \Cref{lem:cs-real-conditions}, a configuration $(A,\psi)$ is \emph{Real} if it is $I_B$-invariant. The real gauge group
\[
\mathcal G_X^\R=\left\{g\colon X\to U(1)\mid g\circ\tau_B=\overline g\right\}
\]
acts naturally on the space of Real configurations. Its identity component $(\mathcal G_X^\R)_0$ acts freely.

\begin{prop}\label{prop:cs-constant-real-solutions}
For all sufficiently large $T$, every Real solution to the monopole equations is mapped by a gauge transformation in $(\mathcal G_X^\R)_0$ to a solution satisfying the following conditions:
\begin{itemize}
    \item $f=0$;
    \item the $u_i$ and $\psi$ are independent of both $s$ and $y$.
\end{itemize}
Such a solution has the form
\[
a=u_1\theta^1+u_2\theta^2,\qquad
\psi=z+\mathbf k\overline z,\qquad
u_1,u_2\in\R,\qquad z\in\C,
\]
where
\begin{equation}\label{eq:cs-detailed-finite-equations}
k_1u_1=-\frac{T}{2}\operatorname{Im}(z^2),\qquad
k_2u_2=-\frac{T}{2}\operatorname{Re}(z^2),\qquad
(u_2+iu_1)\overline z=0.
\end{equation}
\end{prop}

We prove the proposition using the following two lemmas.

\begin{lem}\label{lem:cs-fiberwise-constant-solutions}
For all sufficiently large $T$, every Real solution to the monopole equations is carried by a gauge transformation in $(\mathcal G_X^\R)_0$ to a solution satisfying the following conditions:
\begin{itemize}
    \item $f=0$;
    \item the $u_i$ and $\psi$ are independent of $y$.
\end{itemize}
\end{lem}

\begin{proof}
Let $(A,\psi)=(A_0+2ia,\psi)$ be a Real solution. On $X$, all pointwise norms on forms and spinors are induced by $g_T$, and the $L^p(g_T)$ norms are taken with respect to these norms and $d\operatorname{vol}_{g_T}$. We write $d_{g_T}^\ast$ for the formal adjoint of $d$ and $D_{A_0,g_T}^\pm$ for the Dirac operators. For a $1$-form $\alpha$, we write $d_{g_T}^+\alpha=\frac{1}{2}(d\alpha+\ast_{g_T}d\alpha)$ for the self-dual part of $d\alpha$. The maps $c^\pm$, $\rho^+$, and $q$ are also defined using $g_T$. Norms of bundle endomorphisms are operator norms. With the $k_i$ and $\Lambda$ fixed, the constants $C$ below may change from line to line but are independent of the solution and of all sufficiently large $T$.

The Weitzenb\"ock formula and the maximum principle, together with \Cref{lem: scal}, give
\[
\|\psi\|_{L^\infty(g_T)}^2
\leq -\operatorname{scal}(g_T)
=\frac{\sum_i k_i^2}{T^2}.
\]
Since
\[
\langle d_{g_T}^\ast a,1\rangle_{L^2(g_T)}
=\langle a,d(1)\rangle_{L^2(g_T)}=0,
\]
there exists a solution $\chi$ to $d_{g_T}^\ast d\chi=d_{g_T}^\ast a$. Since $X$ is closed and connected, the kernel of $d_{g_T}^\ast d$ consists of the constant functions, so the normalization $\int_X\chi\,d\operatorname{vol}_{g_T}=0$ determines $\chi$ uniquely. As $\tau_B$ is an isometry of $g_T$ and $\tau_B^\ast a=-a$, uniqueness gives $\chi\circ\tau_B=-\chi$. Thus $e^{i\chi}$ puts $a$ in Coulomb gauge and is joined to the identity by the path $e^{it\chi}$, $0\leq t\leq1$, in $\mathcal G_X^\R$. We may therefore assume that $d_{g_T}^\ast a=0$.

Since $(\theta_s^1,\theta_s^2,\theta_s^3)$ is an oriented orthonormal coframe on $Y$ and $k_1+k_2+k_3=0$, the fiber volume form $dy=d\operatorname{vol}_{g_s}$ is given by
\[
d\operatorname{vol}_{g_s}
=\theta_s^1\wedge\theta_s^2\wedge\theta_s^3
=e^{(k_1+k_2+k_3)s}dy_1\wedge dy_2\wedge dy_3
=dy_1\wedge dy_2\wedge dy_3.
\]
In particular, it is independent of $s$ and satisfies
\[
\int_Y d\operatorname{vol}_{g_s}
=\int_{\R^3/\Lambda}dy_1\wedge dy_2\wedge dy_3
=\operatorname{covol}(\Lambda)=1.
\]
Let $P_0$ denote fiberwise averaging, applied coefficientwise to forms and spinors in the chosen frames, and let $P_\perp$ denote its complement:
\[
P_0h(s)=\int_Y h(s,y)\,dy,\qquad P_\perp=1-P_0.
\]
Since $B$ preserves $dy$, these projections are well-defined on $X$ and preserve the Real conditions. They are orthogonal projections in $L^2(g_T)$. In the chosen frames, the coefficients of $d_{g_T}^\ast$, $d_{g_T}^+$, and $D_{A_0,g_T}^+$ depend only on $s$. Since $dy$ is independent of $s$ and
\[
\int_Y e_i(h)\,dy
=e^{-k_i s}\int_Y\partial_{y_i}h\,dy=0,
\]
all three operators commute with $P_0$.

Write $a=a_0+a_\perp$ and $\psi=\psi_0+\psi_\perp$ for the corresponding decompositions. Averaging $d_{g_T}^\ast a=0$ gives $\frac{d}{ds}(P_0f)=0$, while the Real condition gives $(P_0f)(s+1)=-(P_0f)(s)$. Hence $P_0f=0$. Moreover,
\[
d_{g_T}^\ast a_\perp
=-d_{g_T}^\ast a_0
=-P_0d_{g_T}^\ast a=0.
\]
Since $\varepsilon_1,\varepsilon_2,\varepsilon_3$ are orthonormal, averaging the curvature equation in \eqref{eq:cs-quaternionic-SW} gives
\[
\sup_{s\in[0,2]}
\left|\frac{d}{ds}(P_0u_i)(s)+k_i(P_0u_i)(s)\right|
\leq \frac{T}{4}\|\psi\|_{L^\infty(g_T)}^2
\leq \frac{C}{T}
\qquad (1\leq i\leq 3).
\]
Since $P_0u_i$ is $2$-periodic, evaluating at a point where $|P_0u_i|$ attains its maximum gives
\[
|k_i|\sup_{s\in[0,2]}|(P_0u_i)(s)|\leq\frac{C}{T}.
\]
Each $k_i$ is nonzero. Since $P_0f=0$ and $\|\psi\|_{L^\infty(g_T)}^2\leq \sum_i k_i^2/T^2$, we obtain
\[
\|a_0\|_{L^\infty(g_T)}+\|\psi\|_{L^\infty(g_T)}
\leq\frac{C}{T}.
\]

We next claim that, for all sufficiently large $T$, every positive spinor $\varphi$ with $P_0\varphi=0$ and every real $1$-form $b$ with $P_0b=0$ satisfy
\begin{equation}\label{eq:cs-zero-average-estimates}
\begin{aligned}
\|\varphi\|_{L^2(g_T)}
&\leq C\|D_{A_0,g_T}^+\varphi\|_{L^2(g_T)},\\
\|b\|_{L^2(g_T)}
&\leq C\bigl(\|d_{g_T}^\ast b\|_{L^2(g_T)}
             +\|d_{g_T}^+b\|_{L^2(g_T)}\bigr).
\end{aligned}
\end{equation}
The fiberwise Poincar\'e inequality used in the proof of \eqref{eq:cs-detailed-fiber-gap} gives
\[
\sum_{i=1}^3\|E_i(s)h\|_{L^2(g_s)}^2
\geq\lambda\|h\|_{L^2(g_s)}^2
\]
for every coefficient function $h$ with $P_0h=0$, where $\lambda>0$ is independent of $s\in[0,2]$ and $T$. In the chosen frames, the connection matrices of $\nabla^{A_0}$ and of the Levi--Civita connection $\nabla$ of $g_T$ are $O(T^{-1})$. Applying the inequality coefficientwise and integrating against $T\,ds$ gives
\begin{equation}\label{cs: ineq}
\begin{aligned}
\|\nabla^{A_0}_{g_T}\varphi\|_{L^2(g_T)}^2
&\geq\left(\frac{\lambda}{2}-\frac{C}{T^2}\right)
       \|\varphi\|_{L^2(g_T)}^2,\\
\|\nabla b\|_{L^2(g_T)}^2
&\geq\left(\frac{\lambda}{2}-\frac{C}{T^2}\right)
       \|b\|_{L^2(g_T)}^2.
\end{aligned}
\end{equation}
Since $A_0$ is flat, the Weitzenb\"ock formula and \Cref{lem: scal} give
\[
\|D_{A_0,g_T}^+\varphi\|_{L^2(g_T)}^2
=\|\nabla^{A_0}_{g_T}\varphi\|_{L^2(g_T)}^2
-\frac{1}{4}\int_X\frac{\sum_{i=1}^3 k_i^2}{T^2}
|\varphi|^2\,d\operatorname{vol}_{g_T}.
\]
It follows from \eqref{cs: ineq} that
\begin{equation}\label{cs: ineq2}
\|D_{A_0,g_T}^+\varphi\|_{L^2(g_T)}^2
\geq\left(\frac{\lambda}{2}-\frac{C}{T^2}\right)
\|\varphi\|_{L^2(g_T)}^2.
\end{equation}
Applying the Weitzenb\"ock formula to $b$ gives
\[
\|d_{g_T}^\ast b\|_{L^2(g_T)}^2+\|db\|_{L^2(g_T)}^2
=\|\nabla b\|_{L^2(g_T)}^2
+\int_X\operatorname{Ric}_{g_T}(b,b)\,d\operatorname{vol}_{g_T}.
\]
Since $\|\operatorname{Ric}_{g_T}\|_{L^\infty(g_T)}=O(T^{-2})$ and $\|db\|_{L^2(g_T)}^2=2\|d_{g_T}^+b\|_{L^2(g_T)}^2$ by $\int_X db\wedge db=0$, we deduce from \eqref{cs: ineq} that
\begin{equation}\label{cs: ineq3}
\left(\frac{\lambda}{2}-\frac{C}{T^2}\right)\|b\|_{L^2(g_T)}^2
\leq\|d_{g_T}^\ast b\|_{L^2(g_T)}^2
+2\|d_{g_T}^+b\|_{L^2(g_T)}^2.
\end{equation}
The estimates in \eqref{eq:cs-zero-average-estimates} now follow from \eqref{cs: ineq2} and \eqref{cs: ineq3} for all sufficiently large $T$.

Finally, since $q(\psi_0)$ is independent of $y$, the $P_\perp$-projections of the monopole equations are
\[
\begin{aligned}
2i\,d_{g_T}^+a_\perp
&=P_\perp\bigl(q(\psi)-q(\psi_0)\bigr),\\
D_{A_0,g_T}^+\psi_\perp
&=-c^+(a_0)\psi_\perp\mathbf i
  -P_\perp\bigl(c^+(a_\perp)\psi\mathbf i\bigr).
\end{aligned}
\]
Using $|q(v)-q(w)|_{g_T}\leq C(|v|_{g_T}+|w|_{g_T})|v-w|_{g_T}$, the inequality $\|\psi_0\|_{L^\infty(g_T)}\leq\|\psi\|_{L^\infty(g_T)}$, and the bounds above, we obtain
\[
\begin{aligned}
\|a_\perp\|_{L^2(g_T)}+\|\psi_\perp\|_{L^2(g_T)}
&\leq C\bigl(\|d_{g_T}^+a_\perp\|_{L^2(g_T)}
            +\|D_{A_0,g_T}^+\psi_\perp\|_{L^2(g_T)}\bigr)\\
&\leq\frac{C}{T}
\bigl(\|a_\perp\|_{L^2(g_T)}+\|\psi_\perp\|_{L^2(g_T)}\bigr),
\end{aligned}
\]
where the first inequality follows from \eqref{eq:cs-zero-average-estimates} and $d_{g_T}^\ast a_\perp=0$. Taking $T$ sufficiently large that $C/T<1$ gives $a_\perp=\psi_\perp=0$. Thus the coefficients of $a$ and $\psi$ are independent of $y$. Since $P_0f=0$, we also have $f=P_0f=0$, which completes the proof.
\end{proof}

\begin{lem}\label{lem:cs-s-independent-solutions}
Let $T>0$. Suppose that a solution of \eqref{eq:cs-quaternionic-SW} on $X$ satisfies $f=0$ and that the $u_i$ and $\psi$ are independent of $y$. Then the $u_i$ and $\psi$ are also independent of $s$.
\end{lem}


\begin{proof}
Write
\[
E=\sum_{i=1}^3 u_i\varepsilon_i,\qquad
Q=\psi\mathbf i\bar\psi=\sum_{i=1}^3 Q_i\varepsilon_i.
\]
Since $f=0$ and all vertical derivatives vanish, \eqref{eq:cs-quaternionic-SW} reduces to
\begin{equation}\label{eqn: ode}
\frac{d\psi}{ds}=-TE\psi\mathbf i,\qquad
\frac{du_i}{ds}=\frac{T}{4}Q_i-k_i u_i.
\end{equation}
Since $E$ is purely imaginary, the first equation gives $\frac{d\bar\psi}{ds}=-T\mathbf i\bar\psi E$. Hence
\[
\frac{dQ}{ds}
=\frac{d\psi}{ds}\mathbf i\bar\psi+\psi\mathbf i\frac{d\bar\psi}{ds}
=T(E\psi\bar\psi+\psi\bar\psi E)
=2T|\psi|^2 E.
\]
Comparing coefficients gives $\frac{dQ_i}{ds}=2T|\psi|^2u_i$. Moreover,
\[
\left|\frac{d\psi}{ds}\right|^2
=|-TE\psi\mathbf i|^2
=T^2|E|^2|\psi|^2
=T^2|\psi|^2\sum_{i=1}^3 u_i^2.
\]
Consider the quantity
\[
F=\frac{T}{4}\sum_{i=1}^3 u_iQ_i-\frac12\sum_{i=1}^3 k_i u_i^2.
\]
Differentiating with respect to $s$ and using \eqref{eqn: ode}, we obtain
\[
\begin{aligned}
\frac{dF}{ds}
&=\sum_{i=1}^3\left[
\frac{du_i}{ds}\left(\frac{T}{4}Q_i-k_i u_i\right)
+\frac{T}{4}u_i\frac{dQ_i}{ds}\right]\\
&=\sum_{i=1}^3\left|\frac{du_i}{ds}\right|^2
+\frac{T^2}{2}|\psi|^2\sum_{i=1}^3u_i^2\\
&=\left|\frac{dE}{ds}\right|^2+\frac12\left|\frac{d\psi}{ds}\right|^2.
\end{aligned}
\]
Since the $u_i$ and $\psi$ are $2$-periodic in $s$, so is $F$. Thus
\[
0=F(2)-F(0)
=\int_0^2\left(\left|\frac{dE}{ds}\right|^2
+\frac12\left|\frac{d\psi}{ds}\right|^2\right)\,ds.
\]
Both terms in the integrand are nonnegative, so $\frac{dE}{ds}=\frac{d\psi}{ds}=0$ everywhere. Therefore the $u_i$ and $\psi$ are independent of $s$.
\end{proof}

\begin{rem}
In the flat product case $k_1=k_2=k_3=0$, the change of variables $t=Ts$, $a=-E$, and $\phi=\psi/\sqrt{2}$ identifies \eqref{eqn: ode} with the equations on the center manifold at the trivial reducible on $T^3$ computed by Morgan, Mrowka, and Szab\'o \cite[p.~922, after Lemma~2.5]{MorganMrowkaSzabo1997}.
\end{rem}

\begin{proof}[Proof of \Cref{prop:cs-constant-real-solutions}]
For all sufficiently large $T$, \Cref{lem:cs-fiberwise-constant-solutions} gives a Real gauge transformation in the identity component after which $f=0$ and the coefficients $u_i$ and $\psi$ are independent of $y$. By \Cref{lem:cs-s-independent-solutions}, they are also independent of $s$. The Real condition therefore reduces to
\[
u_3=-u_3,\qquad \psi=\mathbf i\psi\mathbf j.
\]
The first equation gives $u_3=0$. For the second, write $\psi=x_0+x_1\mathbf i+x_2\mathbf j+x_3\mathbf k$ with real coefficients. Then
\[
\mathbf i\psi\mathbf j=x_3-x_2\mathbf i-x_1\mathbf j+x_0\mathbf k.
\]
Hence $\psi=\mathbf i\psi\mathbf j$ is equivalent to $x_3=x_0$ and $x_2=-x_1$. Setting $z=x_0+ix_1$ gives
\[
\psi=z+\mathbf k\overline z.
\]
We then have
\[
\psi\mathbf i\bar\psi
=(z+\mathbf k\overline z)\mathbf i(\overline z-z\mathbf k)
=2z^2\mathbf j
=-2\operatorname{Im}(z^2)\varepsilon_1-2\operatorname{Re}(z^2)\varepsilon_2.
\]
Thus the second equation in \eqref{eq:cs-quaternionic-SW} reduces to
\[
k_1u_1=-\frac T2\operatorname{Im}(z^2),\qquad
k_2u_2=-\frac T2\operatorname{Re}(z^2).
\]
Finally, the first equation in \eqref{eq:cs-quaternionic-SW} reads
\[
-(u_1\mathbf k+u_2\mathbf j)(z+\mathbf k\overline z)\mathbf i=0.
\]
Setting $w=(u_2+iu_1)\overline z$, we can write this equation as $w+\mathbf k\overline w=0$. Since the map
\[
\C\longrightarrow\H;\qquad w\longmapsto w+\mathbf k\overline w
\]
is injective, the equation is equivalent to $(u_2+iu_1)\overline z=0$. Thus the monopole equations reduce to \eqref{eq:cs-detailed-finite-equations}, as desired.
\end{proof}

\begin{rem}\label{rem:cs-degenerate-real-reducible}
On $X=X_{B^2}$, the unperturbed reducible solution $(A_0,0)$ is degenerate even after imposing the Real condition. Indeed, $D_{A_0}^+=T^{-1}\partial_s+D_s$ annihilates constant spinors, which are Real precisely when $\psi=\mathbf i\psi\mathbf j$. Hence, for every $T>0$, we have
\[
\{z+\mathbf k\overline z\mid z\in\C\}\subset(\ker D_{A_0}^+)^{I_B},
\]
so the Real Dirac kernel has real dimension at least two. Thus, although the unperturbed Real moduli space consists of a single reducible solution for all sufficiently large $T$, its absolute degree need not be $1$.
\end{rem}

\subsection{The degree computation}
\label{subsec:cs-adapted-real-return}

We work in the global Coulomb slice for the identity component of the Real gauge group,
\[
(\mathcal G_X^\R)_0
=\{e^{i\chi}\mid \chi\in C^\infty(X;\R),\ \chi\circ\tau_B=-\chi\},
\]
which acts freely on the space of Real configurations $(A,\psi)$. Let $\Omega^k(X)_-$ denote the space of real $k$-forms on $X$ that are anti-invariant under $\tau_B^\ast$. By Hodge theory, the gauge-fixed operator
\[
d^\ast\oplus d^+\colon
\Omega^1(X)_-\longrightarrow\Omega^0(X)_-\oplus\Omega^+(X)_-
\]
has kernel and cokernel given by
\[
\begin{aligned}
\ker(d^\ast\oplus d^+)
&\cong H^1(X;\R)^{-\tau_B^\ast}=0,\\
\operatorname{coker}(d^\ast\oplus d^+)
&\cong H^0(X;\R)^{-\tau_B^\ast}\oplus H^2_+(X;\R)^{-\tau_B^\ast}=0.
\end{aligned}
\]
Thus $d^\ast\oplus d^+$ is invertible, so its determinant line has a canonical orientation. We orient the Real kernel and cokernel of $D_{A_0}^+$ using the complex structures induced by $j$. These choices determine the reference orientation of the Real monopole map.

\begin{prop}\label{prop:cs-local-degree-model}
For all sufficiently large $T$, the unperturbed Real monopole equations admit a Kuranishi chart at $(A_0,0)$ whose obstruction map, defined on a neighborhood $U_0$ of $0$ in $\C$, is
\begin{equation}\label{eq:cs-detailed-h}
h_T(z)
=-\frac T2\overline z
\left(
\frac{\operatorname{Re}(z^2)}{k_2}
+i\frac{\operatorname{Im}(z^2)}{k_1}
\right),
\qquad z\in U_0.
\end{equation}
Its local degree at the origin, computed with the usual complex orientations, equals the local oriented degree of the Real monopole map with the reference orientation fixed above.
\end{prop}

\begin{proof}
Fix $T>0$ sufficiently large as in \Cref{prop:cs-constant-real-solutions}. We use the local Kuranishi model for the Real monopole map; see \cite[Theorem~3.5(1)]{Konno-Taniguchi:2020-PSC}. Throughout the proof, we regard all Sobolev spaces as real Hilbert spaces.

Fix an integer $k\geq3$. Writing $\Lambda^j$ for the bundle of real $j$-forms on $X$, consider the spaces
\[
\begin{aligned}
V&=L^2_k(\Lambda^1)_-\oplus L^2_k(S_B^+)^{I_B},\\
V'&=L^2_{k-1}(\Lambda^0\oplus\Lambda^+)_-
\oplus L^2_{k-1}(S_B^-)^{I_B}.
\end{aligned}
\]
We represent the monopole equations in Coulomb gauge by the map
\[
\mu\colon V\longrightarrow V';\qquad
\mu(a,\psi)
=\left(d^\ast a,\,
d^+a-\frac{q(\psi)}{2i},\,
D_{A_0+2ia}^+\psi\right).
\]
Since $F_{A_0+2ia}=2i\,da$, the second component is the curvature equation divided by $2i$, written in terms of real self-dual forms. The first component imposes the Coulomb condition $d^\ast a=0$, which is not imposed on the domain $V$.

As shown in the proof of \Cref{lem:cs-fiberwise-constant-solutions}, the normalized Poisson equation
\[
d^\ast d\chi=d^\ast a,\qquad
\int_X\chi\,d\operatorname{vol}_{g_T}=0
\]
has a unique solution $\chi\in L^2_{k+1}(X;\R)$, and $\chi\circ\tau_B=-\chi$. Thus $e^{i\chi}\in(\mathcal G_X^\R)_0$ sends $a$ to $a-d\chi$, which satisfies
\[
d^\ast(a-d\chi)=d^\ast a-d^\ast d\chi=0.
\]
Conversely, if $e^{i\chi}\in(\mathcal G_X^\R)_0$ carries one Coulomb representative to another, then $d^\ast d\chi=0$. Hence $\chi$ is constant, and the condition $\chi\circ\tau_B=-\chi$ gives $\chi=0$. Thus each $(\mathcal G_X^\R)_0$-orbit has a unique Coulomb representative, and $\mu^{-1}(0)$ is identified with the framed Real monopole moduli space.

Since $k>\frac12\dim X$, Sobolev multiplication $L^2_k\times L^2_k\to L^2_k$ is continuous. Thus $\mu$, being the sum of a bounded linear map and continuous quadratic terms, is smooth. Its derivative $L\coloneqq d\mu_{(0,0)}$ is given by
\[
L(a,\psi)=(d^\ast a,d^+a,D_{A_0}^+\psi).
\]
The operator $L$ is Fredholm, and its form component $d^\ast\oplus d^+$ is invertible, as shown above.

We next identify the kernel and cokernel of the spinor component $D_{A_0}^+$. The zero-average estimate \eqref{eq:cs-zero-average-estimates} also holds for $D_{A_0}^-$ by the same Weitzenb\"ock argument, since $A_0$ is flat. Both operators commute with the fiberwise averaging operator $P_0$. Thus, if $\psi$ lies in the kernel of either operator, its zero-average part $\psi-P_0\psi$ lies in the same kernel and vanishes by the estimate. On fiberwise-constant spinors, $D_{A_0}^\pm$ act by $\pm T^{-1}\partial_s$, so their kernels consist of constant spinors. By the proof of \Cref{prop:cs-constant-real-solutions}, the constant Real spinors are precisely
\[
\psi_z=z+\mathbf k\overline z,\qquad z\in\C.
\]
The operators
\[
D_{A_0}^\pm\colon L^2_k(S_B^\pm)^{I_B}
\longrightarrow L^2_{k-1}(S_B^\mp)^{I_B}
\]
are formal $L^2$-adjoints of one another, so elliptic regularity identifies the cokernel of $D_{A_0}^+$ with the kernel of $D_{A_0}^-$. We therefore have real-linear isomorphisms
\[
\begin{aligned}
\iota_+\colon\C&\longrightarrow\ker L;
&\qquad \iota_+(z)&=(0,\psi_z),\\
\iota_-\colon\C&\longrightarrow(\operatorname{Im}L)^\perp;
&\qquad \iota_-(z)&=(0,0,\psi_z),
\end{aligned}
\]
where orthogonality is taken with respect to $L^2$. In the second line, $\psi_z$ is regarded as a section of $S_B^-$. In particular, $L$ has index zero.

We now describe the Kuranishi coordinates. With respect to the $L^2$-orthogonal splittings
\[
V=\ker L\oplus(\ker L)^\perp,\qquad
V'=\operatorname{Im}L\oplus(\operatorname{Im}L)^\perp,
\]
the restriction $D\coloneqq L|_{(\ker L)^\perp}\colon(\ker L)^\perp\to\operatorname{Im}L$ is an isomorphism with bounded inverse. Define the Kuranishi coordinate map $\Psi\colon V\to V$ by
\[
\Psi(x)
=x+D^{-1}\operatorname{pr}_{\operatorname{Im}L}\bigl(\mu(x)-Lx\bigr)
=\operatorname{pr}_{\ker L}x
+D^{-1}\operatorname{pr}_{\operatorname{Im}L}\mu(x).
\]
Since $\Psi(0)=0$ and $d\Psi_0=\mathrm{id}$, the inverse function theorem gives a smooth local inverse $\mathcal T=\Psi^{-1}$ on a neighborhood $U$ of $0$ in $V$, with $d\mathcal T_0=\mathrm{id}$. Writing
\[
\wt\mu(v,w)
=\operatorname{pr}_{(\operatorname{Im}L)^\perp}
\mu(\mathcal T(v,w)),
\]
we obtain the Kuranishi normal form
\[
\mu\circ\mathcal T(v,w)=\bigl(Dw,\wt\mu(v,w)\bigr).
\]
For a sufficiently small neighborhood $U_0$ of $0$ in $\C$, define $\mathcal T_0\colon U_0\to V$ and the obstruction map $h_T\colon U_0\to\C$ by
\[
\mathcal T_0(z)=\mathcal T(\iota_+(z),0),\qquad
h_T(z)=\iota_-^{-1}\wt\mu(\iota_+(z),0).
\]
Then the following diagram commutes:
\[
\begin{tikzcd}[column sep=large,row sep=large]
U_0 \arrow[r,"h_T"] \arrow[d,"\mathcal T_0"']
  & \C \arrow[d,"\iota_-"]\\
V \arrow[r,"\mu"']
  & V'.
\end{tikzcd}
\]
Since $D$ is injective, every zero $(v,w)$ of $\mu\circ\mathcal T$ has $w=0$. Hence $\mathcal T_0$ identifies $h_T^{-1}(0)$ with $\mu^{-1}(0)$ near the origin, with inverse $\iota_+^{-1}\operatorname{pr}_{\ker L}$.

We compute this chart explicitly. For $z$ sufficiently close to $0$, set
\[
\begin{aligned}
\psi_z&=z+\mathbf k\overline z,\\
a(z)&=-\frac{T}{2k_1}\operatorname{Im}(z^2)\theta^1
      -\frac{T}{2k_2}\operatorname{Re}(z^2)\theta^2,\qquad
x_z=(a(z),\psi_z)\in V.
\end{aligned}
\]
Here $\psi_z$ is a constant section of $S_B^+=X\times\H$. The form $a(z)$ is anti-invariant under $\tau_B^\ast$ and satisfies $d^\ast a(z)=0$, since its coefficients are constant and its $\theta^0$-coefficient vanishes. Write
\[
a(z)=u_1(z)\theta^1+u_2(z)\theta^2.
\]
The identity
\[
\begin{aligned}
\frac14\psi_z\mathbf i\overline{\psi_z}
&=-\frac12\operatorname{Im}(z^2)\varepsilon_1
  -\frac12\operatorname{Re}(z^2)\varepsilon_2\\
&=\frac{k_1u_1(z)}{T}\varepsilon_1
  +\frac{k_2u_2(z)}{T}\varepsilon_2,
\end{aligned}
\]
together with \eqref{eq:cs-quaternionic-SW}, gives $d^+a(z)=q(\psi_z)/(2i)$. The spinor component is
\[
D_{A_0+2ia(z)}^+\psi_z
=\bigl(u_1(z)\varepsilon_1+u_2(z)\varepsilon_2\bigr)\psi_z\mathbf i
=\psi_{(u_2(z)+iu_1(z))\overline z}.
\]
This is a constant Real negative spinor, so it lies in $\ker D_{A_0}^-$. Consequently,
\[
\operatorname{pr}_{\operatorname{Im}L}\mu(x_z)=0,\qquad
\operatorname{pr}_{\ker L}x_z=\iota_+(z).
\]
It follows from the definition of $\Psi$ that $\Psi(x_z)=(\iota_+(z),0)$. Since $x_z\to0$ as $z\to0$, we may shrink $U_0$ so that $x_z$ lies in the neighborhood on which $\Psi$ is invertible. Uniqueness of its local inverse gives
\[
\mathcal T_0(z)=\mathcal T(\iota_+(z),0)=x_z=(a(z),\psi_z).
\]
Since $\mu(x_z)$ already lies in $(\operatorname{Im}L)^\perp$, we obtain
\[
\wt\mu(\iota_+(z),0)
=\mu(x_z)
=\bigl(0,0,\psi_{(u_2(z)+iu_1(z))\overline z}\bigr).
\]
Thus
\[
h_T(z)=\iota_-^{-1}\wt\mu(\iota_+(z),0)
=(u_2(z)+iu_1(z))\overline z.
\]
Substituting the coefficients of $a(z)$ gives \eqref{eq:cs-detailed-h}.

It remains to compare orientations. In Coulomb gauge, the local oriented degree of the Real monopole map is the local degree of $\mu$ at the origin. The coordinate change $\mathcal T$ has derivative $\mathrm{id}$ at the origin, so it does not change this degree. Since $k_1$ and $k_2$ are nonzero, \eqref{eq:cs-detailed-h} shows that $h_T(z)=0$ only at $z=0$. After shrinking $U$ to a product neighborhood, consider the homotopy
\[
H_t(v,w)=\bigl(Dw,\wt\mu(v,tw)\bigr),\qquad 0\leq t\leq1.
\]
Every zero of $H_t$ has $w=0$, and its remaining equation is $\wt\mu(v,0)=0$. Thus the origin remains the only zero throughout the homotopy, which reduces the degree computation to the product map
\[
H_0(v,w)=\bigl(Dw,\iota_-h_T(\iota_+^{-1}(v))\bigr).
\]
With the reference orientation fixed above, the invertible block $D$ contributes no sign: its form part has the canonical determinant-line orientation, and its spinor part is complex-linear with respect to $j$. Finally, $j(\psi_z)=\psi_z\mathbf j=\psi_{-iz}$, so $\iota_+$ and $\iota_-$ both reverse the complex orientations induced by $j$. The two signs cancel, and therefore
\[
\deg_0\mu=\deg_0 h_T,
\]
where the degree on the right is computed with the usual complex orientations.
\end{proof}

Since $d^\ast\oplus d^+$ is invertible and the Real Dirac index is zero, the Real monopole map has virtual dimension zero. Let $BF^R_{X_{B^2}}$ denote the real Bauer--Furuta map for $I_B=j\widetilde\tau_B$, constructed using $A_0$ as the reference connection. With the reference orientation fixed above, we define
\[
e(B)\coloneqq\deg BF^R_{X_{B^2}},
\]
where $\deg$ denotes the ordinary signed degree.

\begin{prop}\label{prop: CS negative degree}
Let $B\in SL_3(\Z)$ satisfy $\det(I-B^2)=\pm1$ and $\det(I+B)=-1$. Then
\[
e(B)=-3.
\]
In particular, the absolute real Bauer--Furuta degree is $3$.
\end{prop}

\begin{proof}
Fix $T>0$ sufficiently large that the conclusions of \Cref{prop:cs-constant-real-solutions,prop:cs-local-degree-model} hold. If $z\neq0$, then $(u_2+iu_1)\overline z=0$ gives $u_1=u_2=0$. The remaining equations give
\[
-\frac T2\operatorname{Im}(z^2)=k_1u_1=0,\qquad
-\frac T2\operatorname{Re}(z^2)=k_2u_2=0,
\]
so $z^2=0$, a contradiction. Thus $z=0$, and since $k_1$ and $k_2$ are nonzero, we also have $u_1=u_2=0$. Hence $(A_0,0)$ is the unique solution to the unperturbed Real monopole equations modulo $(\mathcal G_X^\R)_0$.

To compute the local degree of $h_T$, note that the winding numbers of its two factors add. Since $k_1<0<k_2$, the real-linear map $x+iy\mapsto x/k_2+iy/k_1$ reverses orientation. Thus \eqref{eq:cs-detailed-h} gives
\[
\begin{aligned}
\deg_0 h_T
&=\deg_0(z\mapsto\overline z)
+\deg_0\left(z\mapsto
\frac{\operatorname{Re}(z^2)}{k_2}
+i\frac{\operatorname{Im}(z^2)}{k_1}\right)\\
&=(-1)+(-2)=-3.
\end{aligned}
\]
By \Cref{prop:cs-local-degree-model}, the local oriented degree of the Real monopole map at $(A_0,0)$ is therefore $-3$.

It remains to identify this local degree with $e(B)$. Choose a closed disk $\Delta\subset U_0$ centered at the origin and a sufficiently small regular value $\eta\in\C$ of $h_T$ such that
\[
0<|\eta|<\min_{\partial\Delta}|h_T|.
\]
Using the cokernel inclusion $\iota_-\colon\C\to(\operatorname{Im}L)^\perp$ from the proof of \Cref{prop:cs-local-degree-model}, define the finite-rank perturbation
\[
\mu_\eta\colon V\longrightarrow V';\qquad
\mu_\eta(x)=\mu(x)-\iota_-(\eta).
\]
In Kuranishi coordinates, we have
\[
\mu_\eta(\mathcal T(v,w))
=\bigl(Dw,\wt\mu(v,w)-\iota_-(\eta)\bigr).
\]
Since $D$ is invertible, the zeros in this chart are precisely the points $\mathcal T_0(z)$ with $h_T(z)=\eta$. Each such zero is regular, and its sign agrees with $\operatorname{sgn}\det dh_T(z)$ by the orientation comparison in \Cref{prop:cs-local-degree-model}.

By properness of the monopole map on the Coulomb slice and uniqueness of its unperturbed zero, for all sufficiently small $\eta$ and every $t\in[0,1]$, all zeros of $\mu_{t\eta}$ lie in the Kuranishi chart over $\operatorname{int}\Delta$; see \cite[Lemma~2.2 and Proposition~3.1]{BF04} for properness.\footnote{Otherwise, there would be sequences $\eta_n\to0$, $t_n\in[0,1]$, and zeros $x_n$ outside this neighborhood with $\mu(x_n)=t_n\iota_-(\eta_n)\to0$. Properness gives a subsequence converging to the unique unperturbed zero, a contradiction.}
Thus $0$ is a regular value of $\mu_\eta$. Since $h_T(z)\neq t\eta$ on $\partial\Delta$ for every $t\in[0,1]$, the homotopy
\[
H_t=\mu-t\iota_-(\eta),\qquad 0\leq t\leq1,
\]
realizes the ordinary degree of the fixed-point map as the signed count of these zeros. Therefore
\[
e(B)
=\sum_{\substack{z\in\operatorname{int}\Delta\\ h_T(z)=\eta}}
  \operatorname{sgn}\det dh_T(z)
=\deg_0 h_T
=-3,
\]
as desired.
\end{proof}

\begin{cor}\label{cor: Miyazawa invariant of special CS knots}
For every special Cappell--Shaneson $2$-knot $K_B$, we have
\[
\mathfrak M(S^4,K_B)=|e(B)|=3.
\]
\end{cor}

\begin{proof}
Let $l$ be the integral sign local system on $X_B$ associated with the double cover $X_{B^2}\to X_B$. Since the monodromy acts by $-1$ on the coefficients, the twisted Wang sequence gives a short exact sequence
\[
0\longrightarrow\operatorname{coker}(I+B^\ast)
\longrightarrow H^2(X_B;l)
\longrightarrow\ker(I+\Lambda^2B^\ast)
\longrightarrow0,
\]
where $B^\ast$ acts on $H^1(T^3;\Z)$ and $\Lambda^2B^\ast$ acts on $H^2(T^3;\Z)$. Since $\det(I+B^\ast)=\det(I+B)=-1$, we have $\operatorname{coker}(I+B^\ast)=0$. Moreover, since $\det B=1$, Poincar\'e duality identifies $\Lambda^2B^\ast$ with the action of $B^{-1}$ on $H_1(T^3;\Z)$. Thus
\[
\det(I+\Lambda^2B^\ast)
=\det(I+B^{-1})
=(\det B)^{-1}\det(I+B)
=-1.
\]
It follows that $\ker(I+\Lambda^2B^\ast)=0$, and hence $H^2(X_B;l)=0$. By \Cref{lem:cs-real-conditions}, $I_B$ defines a real $\mathrm{spin}^c$ structure on $X_{B^2}$. Since $H^2(X_B;l)=0$, the classification in \cite[Lemma~2.8]{miyazawa2025satellite} identifies this structure with the one used in the real surgery formula.

By construction, $K_B$ is the belt sphere of the circle surgery on $\gamma_B$, which replaces $S^1\times D^3$ by $D^2\times S^2$. On the attached piece in the double cover, the involution acts by $(z,u)\mapsto(-z,u)$, so the branch locus is the sphere $\{0\}\times S^2$ in the quotient piece. Surgery on $K_B$ reverses this construction. Thus the manifold $X(K_B)$ obtained by surgery on $K_B$ is diffeomorphic to $X_B$, and its nontrivial double cover is $X_{B^2}$. The surgery formula \cite[Proposition~7.3]{miyazawa2025satellite} and \Cref{prop: CS negative degree} now give
\[
\mathfrak M(S^4,K_B)=|\deg(X_B)|=|e(B)|=3,
\]
as desired.
\end{proof}

\section{Proofs of the main theorems}\label{sec: top construction}

\subsection{Exotic diffeomorphisms of $\mathbb F_2$-homology $4$-spheres}

In this subsection, we prove \Cref{thm: main}. We first recall a standard fact about double covers branched along smoothly embedded $2$-spheres.

\begin{lem}\label{lem: branched double cover is F2 homology sphere}
Let $X$ be a smooth $\mathbb F_2$-homology $4$-sphere, and let $F\subset X$ be a smoothly embedded $2$-sphere. Then $\Sigma(F)$ is also an $\mathbb F_2$-homology $4$-sphere.
\end{lem}

\begin{proof}
For the branched double cover $p\colon\Sigma(F)\to X$, the Smith exact sequence \cite[Theorem~1]{lee1995homology} gives
\[
\cdots\longrightarrow H_i(X,F;\mathbb F_2)
\longrightarrow H_i(\Sigma(F);\mathbb F_2)
\overset{p_\ast}{\longrightarrow} H_i(X;\mathbb F_2)
\longrightarrow H_{i-1}(X,F;\mathbb F_2)
\longrightarrow\cdots.
\]
The long exact sequence of the pair $(X,F)$ gives $H_i(X,F;\mathbb F_2)=0$ for $i=0,1,2$. It follows that
\[
H_0(\Sigma(F);\mathbb F_2)\cong\mathbb F_2,\qquad
H_1(\Sigma(F);\mathbb F_2)=H_2(\Sigma(F);\mathbb F_2)=0.
\]
Poincar\'e duality and the universal coefficient theorem then give
\[
H_3(\Sigma(F);\mathbb F_2)=0,\qquad
H_4(\Sigma(F);\mathbb F_2)\cong\mathbb F_2.
\]
Thus $\Sigma(F)$ is an $\mathbb F_2$-homology $4$-sphere.
\end{proof}

We are now ready to prove \Cref{thm: main}, whose statement we recall.

\MainSphereTheorem*

\begin{proof}
Choose a special Cappell--Shaneson $2$-knot $K$ as in \Cref{lem: cappell-shaneson 2-knots exist}. Since $\Sigma(K)\cong S^4$, we may regard its deck involution $\tau_K$ as an orientation-preserving diffeomorphism of $S^4$. By \Cref{thm: Lin to Miyazawa,cor: Miyazawa invariant of special CS knots}, we have
\[
\mathfrak{LM}_{S^4}(\tau_K)
=\frac{\mathfrak M(S^4,K)^2-1}{8}
=\frac{3^2-1}{8}
=1
\qquad\text{in }\Z_2.
\]

By smooth isotopy extension, choose a diffeomorphism $f_0$ smoothly isotopic to $\tau_K$ and supported in the interior of a standard $4$-ball. Form the connected sum $f=f_0\#\id_X$ away from this ball. By \Cref{lem: connected sum and Lin invariant}, we have
\[
\mathfrak{LM}_X(f)
=\mathfrak{LM}_{S^4}(f_0)+\mathfrak{LM}_X(\id_X)
=1.
\]
Thus $f$ is not smoothly isotopic to the identity. Since $f$ is supported in a standard $4$-ball in $X$, the Alexander trick gives a topological isotopy from $f$ to the identity. Hence $f$ is exotic. Finally, \Cref{lem: space of psc metrics and lin invariant} implies that $\mathcal R_+(X)$ is either empty or disconnected.
\end{proof}


\subsection{Homotopy groups of $\mathcal R_+(S^4)$}

Let $X$ be an oriented smooth $\mathbb F_2$-homology $4$-sphere, and let $F$ be a smoothly embedded $2$-sphere in $X$. Let $\tau$ denote the deck involution of the branched double cover $\Sigma(F)\to X$, and choose a spin lift $\tilde\tau$ with respect to the unique spin structure $\s$ on $\Sigma(F)$. As discussed in \Cref{sec: relation to miyazawa}, this gives the $G_{\mathrm{odd}}$-equivariant Bauer--Furuta invariant $BF_{\Sigma(F),\s,\tilde\tau}$. By \Cref{lem: classifying integer,rem: signed refinement of miyazawa}, there is a unique integer $k$ such that
\[
BF_{\Sigma(F),\s,\tilde\tau}\sim_{G_{\mathrm{odd}}}\mathrm{id}+ke,
\qquad
\mathfrak M_{\mathrm{sign}}(X,F)=4k+1,
\]
where $e=\chi_{G_{\mathrm{odd}}}(G_{\mathrm{odd}}/\langle I\rangle)$ and $I=j\tilde\tau$.

For $n\geq1$, consider the smooth spin bundle
\[
E_n(X,F)=(\Sigma(F),\s)\times S^n\longrightarrow S^n.
\]
Let $\Z_2$ act diagonally by $\tau$ on $\Sigma(F)$ and by the antipodal map on $S^n$. The quotient is a smooth bundle over $\R P^n$. However, since $\tilde\tau^2=-1$, the chosen lift does not define a $\Z_2$-action on the vertical spin bundle and hence does not induce a fiberwise spin structure on the quotient.

We recall the parametrized version of the equivariant $KO$-degree. Let $K$ be a compact Lie group, and let $E\to B$ be a real $K$-equivariant vector bundle over a compact $K$-space. 
To define the $KO$-degree, we pass to a complete $K$-universe
before desuspending by $E$.
Write $E^+$ for its fiberwise one-point compactification. A fiberwise based $K$-equivariant map
\[
F\colon E^+\longrightarrow E^+
\]
determines, after fiberwise stabilization and desuspension by $E$, a map
\[
\Sigma_B^{-E}F\colon\underline{S^0}_B\longrightarrow\underline{S^0}_B
\]
in the genuine $K$-equivariant stable homotopy category parametrized over $B$. We define its equivariant $KO$-degree by
\[
\deg_{KO}(F):=
\left[
\underline{S^0}_B
\xrightarrow{\Sigma_B^{-E}F}
\underline{S^0}_B
\xrightarrow{\eta_B}
\underline{KO}_B
\right]
\in KO_K^0(B),
\]
where $\underline{S^0}_B$ and $\underline{KO}_B$ are the pullbacks to $B$ of the $K$-equivariant sphere spectrum and the $K$-equivariant real $K$-theory spectrum, respectively, and $\eta_B$ is the unit map.

For a $K$-equivariant vector bundle automorphism $f\colon E\to E$ covering the identity of $B$, we write $\deg_{KO}(f):=\deg_{KO}(f^+)$. In particular,
\[
\deg_{KO}(\mathrm{id}_E)=[\underline{\R}]=1\in KO_K^0(B).
\]
When $B$ is a point, this recovers the equivariant $KO$-degree map $\pi_0^K(S^0)\to KO_K^0(\mathrm{pt})=RO(K)$.

\begin{thm}\label{thm: euler characteristic to RO}
Let $K$ be a compact Lie group, and let $X$ be a finite $K$-CW complex whose suspension spectrum $\Sigma^\infty_K X_+$ is dualizable in the chosen $K$-universe. Under the $KO$-degree map
\[
\deg_{KO}\colon\pi_0^K(S^0)\longrightarrow RO(K),
\]
the equivariant Euler characteristic $\chi_K(X)$ maps to $\sum_{i\geq 0}(-1)^i[H^i(X;\R)]$.
\end{thm}
\begin{proof}
Change of universe preserves the equivariant Euler characteristic and its $KO$-degree, so it suffices to work in a complete $K$-universe. In this case, the formula for an orbit $K/H$ is given in \cite[Section~14]{fausk2008survey}. Since both sides are additive under equivariant cell attachments, the result follows by cellular induction.
\end{proof}

\begin{lem}\label{lem: invariance of KO-degree}
Let $B$ be a compact $G_{\mathrm{odd}}/\mathrm{Pin}(2)$-space, and let $E\to B$ be a quaternionic $G_{\mathrm{odd}}$-equivariant vector bundle (so that $\mathrm{Pin}(2)$ acts by left multiplication on each fiber). Then every $G_{\mathrm{odd}}$-equivariant automorphism of $E$ has $G_{\mathrm{odd}}$-equivariant $KO$-degree $1$.
\end{lem}

\begin{proof}
Let $f$ be a $G_{\mathrm{odd}}$-equivariant automorphism of $E$. For the induced parametrized automorphism $f^+\colon E^+\to E^+$, where $E^+$ denotes the fiberwise one-point compactification of $E$, we have
\[
\deg_{KO}(f^+)=[\operatorname{sgn}f]\in KO^0_{G_{\mathrm{odd}}}(B),
\]
where $\operatorname{sgn}f$ denotes the sign bundle of $f$, a $G_{\mathrm{odd}}$-equivariant real line bundle over $B$.\footnote{The sign bundle can be described as follows. Choose a $G_{\mathrm{odd}}$-invariant metric on $E$, and let $u=f(f^\ast f)^{-1/2}$ be the orthogonal part of $f$. Since $u_b$ commutes with the $\mathrm{Pin}(2)$-action, it lies in $\mathrm{Sp}(m)\subset SO(E_b)$, where $m=\operatorname{rank}_{\H}E$. Let $P_f\to B$ be the double cover whose fiber over $b$ consists of the two lifts of $u_b$ under $\mathrm{Spin}(E_b)\to SO(E_b)$. Then
\[
\operatorname{sgn}f=P_f\times_{\{\pm1\}}\R,
\]
where $-1$ acts on $\R$ by multiplication by $-1$. The action on $E$ induces the $G_{\mathrm{odd}}$-action on this line bundle.}
Since $\mathrm{Sp}(m)$ is simply connected for $m>0$, the sign bundle admits a canonical fiberwise trivialization normalized at the identity. This trivialization is $G_{\mathrm{odd}}$-equivariant and compatible with changes of coordinates, so the sign bundle is globally trivial. Hence
\[
\deg_{KO}(f^+)=[\underline{\R}]=1\in KO^0_{G_{\mathrm{odd}}}(B),
\]
as desired.
\end{proof}

\begin{lem}\label{lem: psc on n-skeleton and KO}
Let $k\in\Z$ be the integer determined by $\mathfrak M_{\mathrm{sign}}(X,F)=4k+1$, and let $L\in KO^0(\R P^n)$ denote the class of the tautological real line bundle. Suppose that $E_n(X,F)$ admits a $\Z_2$-equivariant fiberwise psc metric. Then
\[
k(L-1)=0\quad\text{in}\quad KO^0(\R P^n).
\]
\end{lem}

\begin{proof}
Let $p\colon S^n\to\ast$ be the collapse map. Choose fixed
$\tilde\tau$-equivariant geometric data on $(\Sigma(F),\s)$,
with the canonical flat determinant connection, and pull them back
to $S^n$. The family monopole map and its finite-dimensional
approximation are then pullbacks, and the equivariant quaternionic
index is the pullback of the zero representation class computed in
\Cref{sec: relation to miyazawa}. With the real normalization fixed
there, \Cref{lem: invariance of KO-degree} makes the $KO$-degree
independent of the stable quaternionic index identification. Hence
\[
\deg_{KO}(BF_{E_n(X,F)})
=p^*\deg_{KO}(BF_{\Sigma(F),\s,\tilde\tau})
\quad\text{in }KO^0_{G_{\mathrm{odd}}}(S^n).
\]
Since
\[
\begin{aligned}
H^0(G_{\mathrm{odd}}/\langle I\rangle;\R)
&\cong \R\oplus\wt\R_I,\quad 
H^1(G_{\mathrm{odd}}/\langle I\rangle;\R)\cong \wt\R_j\oplus\wt\R_\tau,
\end{aligned}
\]
the equality
$BF_{\Sigma(F),\s,\tilde\tau}\sim_{G_{\mathrm{odd}}}\mathrm{id}+ke$
and \Cref{thm: euler characteristic to RO} give
\[
\deg_{KO}(BF_{\Sigma(F),\s,\tilde\tau})
=[\R]+k\left([\R]+[\wt\R_I]-[\wt\R_j]-[\wt\R_\tau]\right)
\]
in $RO(G_{\mathrm{odd}})$.

Let $g_b$ be the assumed equivariant famly psc metrics . The equivariant
path $g_{t,b}=(1-t)g_0+t g_b$, where $g_0$ is the fixed metric above,
transports the index identification to this family; the path need
not remain psc. Along it, identify the real summands by
$d^*\oplus d^+$ with mean-zero functions in the target, which is
invertible since $b_1(\Sigma(F))=b^+(\Sigma(F))=0$.
At the psc endpoint, the Dirac operators are isomorphisms by the
Weitzenb\"ock formula and index zero. The homotopy and boundedness
argument of \cite[Lemma~2.8]{baraglia2026exotic} apply over the
compact base $S^n$ and respect both the deck lift and the
quaternionic symmetry. The linear endpoint has identity real
block; after a common quaternionic stabilization, its spinor block
in the transported identification is an equivariant quaternionic
automorphism. Thus \Cref{lem: invariance of KO-degree} gives
$\deg_{KO}(BF_{E_n(X,F)})=1$ with the same normalization. Hence
\[
k\,p^*\left([\R]+[\wt\R_I]-[\wt\R_j]-[\wt\R_\tau]\right)=0
\quad\text{in }KO^0_{G_{\mathrm{odd}}}(S^n).
\]

To conclude, consider the $\mathrm{Pin}(2)$-fixed-point map
\[
\begin{aligned}
\Phi^{\mathrm{Pin}(2)}\colon KO^0_{G_{\mathrm{odd}}}(S^n)
&\longrightarrow KO^0_{G_{\mathrm{odd}}/\mathrm{Pin}(2)}(S^n) \cong KO^0_{\Z_2}(S^n)\cong KO^0(\R P^n).
\end{aligned}
\]
This map is additive (but not multiplicative) and satisfies
\[
\begin{aligned}
\Phi^{\mathrm{Pin}(2)}(p^*[\R])=1,
\Phi^{\mathrm{Pin}(2)}(p^*[\wt\R_\tau])=L, \ \Phi^{\mathrm{Pin}(2)}(p^*[\wt\R_j])=0,\ \Phi^{\mathrm{Pin}(2)}(p^*[\wt\R_I])=0.
\end{aligned}
\]
Therefore
\[
\begin{aligned}
0&=\Phi^{\mathrm{Pin}(2)}\left(
k\,p^*\left([\R]+[\wt\R_I]-[\wt\R_j]-[\wt\R_\tau]\right)
\right) =k(1-L),
\end{aligned}
\]
as desired.
\end{proof}

To state the following corollary, we define, for each integer $n\geq1$,
\[
\psi(n)=\left|\left\{i\in\Z\mid 1\leq i\leq n,\ i\equiv 0,1,2,4\pmod{8}\right\}\right|.
\]

\begin{cor}\label{cor: miyazawa mod condition}
Suppose that $E_n(X,F)$ admits a $\Z_2$-equivariant fiberwise psc metric. Then
\[
\mathfrak M_{\mathrm{sign}}(X,F)\equiv 1\pmod{2^{\psi(n)+2}}.
\]
\end{cor}

\begin{proof}
By \Cref{lem: classifying integer,rem: signed refinement of miyazawa,lem: psc on n-skeleton and KO}, there is an integer $k$ such that
\[
\mathfrak M_{\mathrm{sign}}(X,F)=4k+1,
\qquad
k(L-1)=0\quad\text{in }KO^0(\R P^n).
\]
By \cite[Theorem~7.4]{adams1962vector}, the group $\widetilde{KO}^0(\R P^n)$ is cyclic of order $2^{\psi(n)}$, generated by $L-1$. Hence $k$ is divisible by $2^{\psi(n)}$, so
\[
\mathfrak M_{\mathrm{sign}}(X,F)=4k+1\equiv 1\pmod{2^{\psi(n)+2}},
\]
as desired.
\end{proof}

We next consider the deck involution associated with the connected sum of a $2$-knot with itself.

\begin{lem}\label{lem: deck transformation is smoothly trivial}
Let $K\subset S^4$ be a $2$-knot such that $\Sigma(K)$ is diffeomorphic to $S^4$. Then the deck involution of the branched double cover $\Sigma(K\# K)\to S^4$ is smoothly isotopic to the identity.
\end{lem}

\begin{proof}
Write $X=\Sigma(K)$ and let $\tau$ be its deck involution. Near a point of the branching sphere, choose coordinates in which the covering map is $(z,w)\mapsto(z,w^2)$, so that $\tau(z,w)=(z,-w)$. Let $B\subset X$ be the preimage of a small closed round ball centered at the origin in the base. This is a smooth $4$-ball. To form $X\# X$, glue the two copies of $X\smallsetminus\operatorname{int}B$ along $\partial B$ by the reflection
\[
\rho(z,w)=(\bar z,w).
\]
The resulting branched cover is $\Sigma(K\# K)$, and its deck involution is the equivariant connected sum $\tau\#\tau$.

Choose a smooth radial function $\chi$ supported in the coordinate chart and equal to $1$ on a neighborhood of $B$. For $0\leq t\leq 1$, define $f_t$ on the chart by
\[
f_t(z,w)=\bigl(z,e^{\pi i(1-t\chi(z,w))}w\bigr)
\]
and set $f_t=\tau$ outside the chart. Since the rotation angle depends only on the radius, each $f_t$ is a diffeomorphism. It preserves $B$ because it preserves both $|z|$ and $|w|$. We have $f_0=\tau$, and $f=f_1$ is the identity on a neighborhood of $B$. Near $\partial B$, the map $f_t$ is the rotation $(z,w)\mapsto(z,e^{\pi i(1-t)}w)$, which commutes with $\rho$. Hence the two copies of $f_t$ glue to a smooth isotopy from $\tau\#\tau$ to $f\#f$.

Set $D=X\smallsetminus\operatorname{int}B$. Since $X$ is diffeomorphic to $S^4$, the disk theorem allows us to carry $B$ to a standard hemisphere by an ambient isotopy. In particular, $D$ is diffeomorphic to $D^4$. Fix an orientation-preserving identification $X\# X\cong X$. Under this identification, the inclusions of the two copies of $D$ in $X\# X$ and the original inclusion $D\hookrightarrow X$ are three orientation-preserving embeddings of the same disk into $X$. By the disk theorem, these embeddings are ambiently isotopic \cite[Chapter~8, Theorem~3.1]{hirsch1976differential}. Conjugating by such an ambient isotopy carries the extension of $f|_D$ by the identity along one embedding to its extension along another. Thus $f\#\id$, $\id\#f$, and $f$ are smoothly isotopic. It follows that in $\pi_0\mathrm{Diff}^+(X)$ we have
\[
[f\#f]=[\id\#f]\,[f\#\id]=[f]^2=[\tau]^2=1,
\]
as desired.
\end{proof}

We can now prove the first assertion of \Cref{thm: pi 1 of psc metrics on S4}.

\begin{lem}\label{lem: pi 1 of psc components on S4}
Every connected component of $\mathcal R_+(S^4)$ has nontrivial fundamental group.
\end{lem}

\begin{proof}
Let $C$ be a connected component of $\mathcal R_+(S^4)$, and choose $g\in C$. Let $K=K_{A_{-3}}$ be the special Cappell--Shaneson $2$-knot from \Cref{lem: cappell-shaneson 2-knots exist}. By \Cref{lem: deck transformation is smoothly trivial}, the deck involution $\tau$ of the branched double cover $\Sigma(K\# K)\to S^4$ is smoothly isotopic to the identity. Pulling back $g$ along such an isotopy gives a path $\{g_t\}$ in $C$ from $g$ to $\tau^\ast g$. Concatenating this path with its pullback $\{\tau^\ast g_t\}$ gives a loop $\gamma\colon S^1\to C$ that is equivariant for the antipodal action on $S^1$ and the action of $\tau^\ast$ on $C$.

We claim that $\gamma$ is not nullhomotopic. Otherwise, it would extend over one hemisphere of $S^2$. Defining the map on the other hemisphere by equivariance would then give an equivariant map $S^2\to C$ and hence a $\tau$-equivariant fiberwise psc metric on $E_2(S^4,K\# K)$. Since $\psi(2)=2$, \Cref{cor: miyazawa mod condition} and \cite[Theorem~1.8]{miyazawa2023gauge} would imply
\[
9=\mathfrak{M}(S^4,K)^2
=\mathfrak{M}(S^4,K\# K)
\equiv\pm1\pmod{16},
\]
a contradiction. Thus $\gamma$ represents a nontrivial element of $\pi_1(C,g)$.
\end{proof}

We now turn to higher homotopy groups. Let $M$ and $N$ be closed smooth oriented manifolds of the same dimension $n\geq 3$. Recall from \cite[Section~1]{ebert2022positive} that the \emph{psc connected sum map}
\[
\#\colon \mathcal R_+(M)\times\mathcal R_+(N)\longrightarrow\mathcal R_+(M\#N)
\]
is well-defined up to homotopy. When $M=N=S^n$, the resulting multiplication
\[
\#\colon \mathcal R_+(S^n)\times\mathcal R_+(S^n)\longrightarrow\mathcal R_+(S^n)
\]
is homotopy associative and homotopy commutative, with the standard round metric as a homotopy unit. Thus $\pi_0\mathcal R_+(S^n)$ is a commutative monoid whose identity is the component containing the round metric.

\begin{lem}\label{lem: parametrized doubling}
Let $n\geq 1$, and let $C_0$ be the component of $\mathcal R_+(S^4)$ containing the standard round metric. Let $F$ be a smoothly embedded $2$-sphere in a smooth homology $4$-sphere $X$. Suppose that the branched double cover $\Sigma(F)$ is diffeomorphic to $S^4$, and use such an identification to regard its deck involution $\tau$ as a diffeomorphism of $S^4$. If there exists a map $f\colon S^n\to C_0$ satisfying $f(-x)=\tau^\ast f(x)$, then there exists a map $g\colon S^{n+1}\to C_0$ satisfying $g(-x)=(\tau\#\tau)^\ast g(x)$. Here $\tau\#\tau$ is regarded as an involution of $S^4$ via the connected-sum identification used in the multiplication above.
\end{lem}

\begin{proof}
Let $f\colon S^n\to C_0$ be an equivariant map. Choose a $\tau$-invariant ball centered at a point of the fixed sphere, with coordinates in which $\tau$ acts orthogonally. The subspace of $C_0$ consisting of metrics that agree with a fixed rotationally symmetric torpedo metric on this ball is $\tau^\ast$-invariant, and its inclusion is a weak homotopy equivalence by \cite[Theorem~1.1]{Chernysh2004}. Since $S^n$ is a free $\Z_2$-CW complex, we may homotope $f$ equivariantly into this subspace by working over one cell in each orbit and extending over its translate by equivariance. On these metrics, gluing along the standard collars by a reflection that commutes with $\tau$ makes the connected sum commute strictly with pullback by the involutions.

Choose a linear reflection $r\colon S^n\to S^n$ and define
\[
g\colon S^n\longrightarrow C_0;
\qquad x\longmapsto f(x)\# f(r(x)).
\]
Here $g$ takes values in $C_0$ because $C_0$ is the identity in the monoid $\pi_0\mathcal R_+(S^4)$. Since $r$ commutes with the antipodal map, we have
\[
g(-x)=\tau^\ast f(x)\#\tau^\ast f(r(x))=(\tau\#\tau)^\ast g(x).
\]
Thus $g$ is equivariant.

It remains to show that $g$ is nullhomotopic after forgetting equivariance. Since $C_0$ is a path-connected $H$-space, we may identify free homotopy classes of maps $S^n\to C_0$ with $\pi_n(C_0)$, with connected sum inducing addition. As $r$ has degree $-1$, it follows that
\[
[g]=[f]+[f\circ r]=[f]-[f]=0.
\]
Viewing $S^n$ as the equator of $S^{n+1}$, extend $g$ to a map $h$ on the northern hemisphere. On the southern hemisphere, define the map by $x\mapsto(\tau\#\tau)^\ast h(-x)$. The two maps agree on the equator by the equivariance of $g$, giving the desired $\Z_2$-equivariant map $S^{n+1}\to C_0$.
\end{proof}

\begin{lem}\label{lem: top construction for pi n}
For every integer $n\geq 1$, there exists a $2$-knot $F\subset S^4$ with the following properties:
\begin{itemize}
\item $\Sigma(F)$ is diffeomorphic to $S^4$ and $\mathfrak{M}(S^4,F)=3^{2^n}$;
\item $E_n(S^4,F)$ admits a $\Z_2$-equivariant fiberwise psc metric that takes values in $C_0$.
\end{itemize}
\end{lem}

\begin{proof}
For $n=1$, take $F=K\# K$, where $K=K_{A_{-3}}$. The proof of \Cref{lem: pi 1 of psc components on S4} gives $\mathfrak{M}(S^4,F)=9$ and, taking $C=C_0$, an equivariant loop $S^1\to C_0$ that gives the required fiberwise psc metric. Moreover, $\Sigma(F)\cong\Sigma(K)\#\Sigma(K)\cong S^4$.

Suppose that the lemma holds for some $n$ with $2$-knot $F$. Since the standard round metric is a homotopy unit for the psc connected sum map
\[
\#\colon\mathcal R_+(S^4)\times\mathcal R_+(S^4)\longrightarrow\mathcal R_+(S^4),
\]
we have $C_0\#C_0=C_0$. The branched double cover $\Sigma(F\#F)$ is diffeomorphic to $\Sigma(F)\#\Sigma(F)$ and hence to $S^4$, while \Cref{lem: parametrized doubling} shows that $E_{n+1}(S^4,F\#F)$ satisfies the second condition. By \cite[Theorem~1.8]{miyazawa2023gauge}, we have
\[
\mathfrak{M}(S^4,F\#F)=\mathfrak{M}(S^4,F)^2
=\bigl(3^{2^n}\bigr)^2=3^{2^{n+1}}.
\]
The lemma follows by induction on $n$.
\end{proof}

\begin{rem}
The proof of \Cref{lem: pi 1 of psc components on S4} is similar in form to Hitchin's argument for the nontriviality of homotopy groups of spaces of psc metrics; see \cite[Proposition~4.6 and Theorem~4.7]{hitchin1974harmonic}. In Hitchin's setting, a family of diffeomorphisms determines a bundle over a sphere by clutching. A nullhomotopy of the corresponding family of psc metrics would give a fiberwise psc metric on this bundle, forcing its family Dirac index to vanish.

In our setting, a nullhomotopy of the loop $\gamma$ would give, by equivariance on the other hemisphere, a $\Z_2$-equivariant family of psc metrics parametrized by $S^2$. The obstruction comes from the $KO$-degree of the $G_{\mathrm{odd}}$-equivariant Bauer--Furuta map rather than from a family Dirac index: after passing to $\mathrm{Pin}(2)$-fixed subbundles, such a family would force the image of this degree to equal the unit in $KO^0_{\Z_2}(S^2)\cong KO^0(\R P^2)$.
\end{rem}

We now prove a quantitative form of the second assertion of \Cref{thm: pi 1 of psc metrics on S4}. For each integer $r\geq 1$, set
\[
N(r)=\min\left\{m\geq 1\mid\psi(m)=r+1\right\}.
\]

\begin{lem}\label{thm: nontrivial pi n in unbounded degree}
Let $C_0$ be the component of $\mathcal R_+(S^4)$ containing the standard round metric. For every integer $k\geq 2$, there exists an integer $n$ with $k\leq n<N(k)$ such that $\pi_n(C_0)\neq 0$.
\end{lem}

\begin{proof}
Fix $k\geq 2$ and suppose that $\pi_n(C_0)=0$ whenever $k\leq n<N(k)$. By \Cref{lem: top construction for pi n}, there exists a $2$-knot $F\subset S^4$ such that $\Sigma(F)$ is diffeomorphic to $S^4$, $\mathfrak{M}(S^4,F)=3^{2^k}$, and $E_k(S^4,F)$ admits a $\Z_2$-equivariant fiberwise psc metric taking values in $C_0$. Starting with the corresponding equivariant map $S^k\to C_0$, we extend it one dimension at a time. For $k\leq n<N(k)$, the underlying map $S^n\to C_0$ is nullhomotopic because $\pi_n(C_0)=0$. It therefore extends over one hemisphere of $S^{n+1}$, and we define the extension on the other hemisphere by equivariance. Iterating gives a $\Z_2$-equivariant fiberwise psc metric on $E_{N(k)}(S^4,F)$. Since $\psi(N(k))=k+1$, \Cref{cor: miyazawa mod condition} gives
\[
3^{2^k}=\mathfrak{M}(S^4,F)
=\bigl|\mathfrak{M}_{\mathrm{sign}}(S^4,F)\bigr|
\equiv\pm 1\pmod{2^{k+3}}.
\]

On the other hand, writing $\nu_2$ for the $2$-adic valuation, we have
\[
\begin{aligned}
\nu_2(3^{2^k}-1)
&=\nu_2(3-1)+\nu_2(3+1)+\sum_{m=1}^{k-1}\nu_2(3^{2^m}+1)\\
&=1+2+(k-1)=k+2,
\end{aligned}
\]
since $3^{2^m}+1\equiv 2\pmod 4$ for $m\geq 1$. This rules out $+1$ modulo $2^{k+3}$, while $3^{2^k}\equiv 1\pmod 4$ rules out $-1$. Thus
\[
3^{2^k}\not\equiv\pm 1\pmod{2^{k+3}},
\]
a contradiction.
\end{proof}

Since $k$ is arbitrary and $n\geq k$, we have $\pi_n(C_0)\neq 0$ for infinitely many $n\geq 2$. Together with \Cref{lem: pi 1 of psc components on S4}, this completes the proof of \Cref{thm: pi 1 of psc metrics on S4}.

\subsection{Exotic diffeomorphisms of homotopy $K3$ surfaces}

In this subsection, we prove \Cref{thm: K3}. Here we take virtual $\mathrm{Pin}(2)$-representations to be
integral linear combinations of $\R$, $\wt\R$, and $\H$. 
 For a virtual $\mathrm{Pin}(2)$-representation $V$ of rank $n$, let
\[
\mathcal{K}_V\coloneqq
\ker\left[
\Phi^{S^1}\colon\pi_V^{\mathrm{Pin}(2)}(S^0)
\longrightarrow\pi_{V^{S^1}}^{\Z_2}(S^0)
\right],
\]
where $\Phi^{S^1}$ is the $S^1$-fixed-point map. We also use the homomorphism
\[
c_\ast\colon\pi_V^{\mathrm{Pin}(2)}(E\mathrm{Pin}(2)_+)
\longrightarrow\pi_V^{\mathrm{Pin}(2)}(S^0)
\]
induced by the collapse map $c\colon E\mathrm{Pin}(2)_+\to S^0$.
\begin{lem}\label{lem: injectivity condition}
Let $V=r\R+m\H-b\wt\R$, where $r,m,b$ are integers satisfying $r,m\geq 0$ and $b>r+1$. Suppose that there exists a $\mathrm{Pin}(2)$-equivariant stable map $f\colon(m\H)^+\to(b\wt\R)^+$ such that $\deg f^{\mathrm{Pin}(2)}=1$. Then the map
\[
c_\ast\colon\pi_V^{\mathrm{Pin}(2)}(E\mathrm{Pin}(2)_+)\longrightarrow\mathcal{K}_V
\]
is an isomorphism.
\end{lem}

\begin{proof}
We use the models
\[
E\mathrm{Pin}(2)=\operatorname{colim}_k S(k\H),\qquad
\wt{E}\mathrm{Pin}(2)=\operatorname{colim}_k(k\H)^+,
\]
where $S(k\H)$ denotes the unit sphere in $k\H$. The cofiber sequence
\[
E\mathrm{Pin}(2)_+\overset{c}{\longrightarrow}S^0
\overset{i}{\longrightarrow}\wt{E}\mathrm{Pin}(2)
\]
gives the exact sequence
\[
\begin{aligned}
\pi_{V+\R}^{\mathrm{Pin}(2)}(S^0)
&\overset{i_\ast}{\longrightarrow}
\pi_{V+\R}^{\mathrm{Pin}(2)}(\wt{E}\mathrm{Pin}(2))
\overset{\partial}{\longrightarrow}
\pi_V^{\mathrm{Pin}(2)}(E\mathrm{Pin}(2)_+)\\
&\overset{c_\ast}{\longrightarrow}
\pi_V^{\mathrm{Pin}(2)}(S^0)
\overset{i_\ast}{\longrightarrow}
\pi_V^{\mathrm{Pin}(2)}(\wt{E}\mathrm{Pin}(2)).
\end{aligned}
\]
For $k\geq m$, the isomorphism $\pi_{V+\R}^{\mathrm{Pin}(2)}((k\H)^+)\cong\pi_{r+1}^{\mathrm{Pin}(2)}((b\wt\R+(k-m)\H)^+)$ allows us to apply the integer-graded tom Dieck splitting to an actual representation sphere. In our chosen universe, the splitting has three summands:
\[
\pi_{V+\R}^{\mathrm{Pin}(2)}((k\H)^+)
\cong\pi_{r+1}^s(S^0)
\oplus\pi_{r+1}^s(\mathrm{Th}(b\lambda))
\oplus\pi_{r+1}^s(\mathrm{Th}((b+1)\xi_{\wt\R}+(k-m)\xi_\H)),
\]
where $\lambda$ denotes the nonorientable line bundle over $B\Z_2=B(\mathrm{Pin}(2)/S^1)$. The two Thom spaces are associated with real vector bundles of ranks $b$ and $b+1+4(k-m)$, respectively. Both ranks are greater than $r+1$, so both Thom spaces are $(r+1)$-connected. Hence their $(r+1)$st stable homotopy groups vanish by the stable Hurewicz theorem. It follows that taking $\mathrm{Pin}(2)$-fixed points gives an isomorphism
\[
\Phi^{\mathrm{Pin}(2)}\colon
\pi_{V+\R}^{\mathrm{Pin}(2)}((k\H)^+)
\xrightarrow{\cong}\pi_{r+1}^s(S^0)
\]
for every $k\geq m$. Passing to the colimit gives an isomorphism
\[
\Phi^{\mathrm{Pin}(2)}\colon
\pi_{V+\R}^{\mathrm{Pin}(2)}(\wt{E}\mathrm{Pin}(2))
\xrightarrow{\cong}\pi_{r+1}^s(S^0).
\]
Since $\Phi^{\mathrm{Pin}(2)}(f)$ is stably homotopic to the identity, the composition
\[
\pi_{r+1}^s(S^0)
\xrightarrow{\beta\mapsto f\beta}
\pi_{V+\R}^{\mathrm{Pin}(2)}(S^0)
\overset{i_\ast}{\longrightarrow}
\pi_{V+\R}^{\mathrm{Pin}(2)}(\wt{E}\mathrm{Pin}(2))
\xrightarrow{\Phi^{\mathrm{Pin}(2)}}
\pi_{r+1}^s(S^0)
\]
is the identity. Thus $i_\ast$ is surjective, so $\partial=0$ and $c_\ast$ is injective.

For surjectivity, the same connectivity argument in degree $r$ shows that the $S^1$-fixed-point map
\[
\Phi^{S^1}\colon
\pi_V^{\mathrm{Pin}(2)}(\wt{E}\mathrm{Pin}(2))
\longrightarrow\pi_{V^{S^1}}^{\Z_2}(S^0)
=\pi_{r\R-b\wt\R}^{\Z_2}(S^0)
\]
is an isomorphism. Since the composition
\[
\pi_V^{\mathrm{Pin}(2)}(S^0)
\overset{i_\ast}{\longrightarrow}
\pi_V^{\mathrm{Pin}(2)}(\wt{E}\mathrm{Pin}(2))
\xrightarrow{\Phi^{S^1}}
\pi_{V^{S^1}}^{\Z_2}(S^0)
\]
is the $S^1$-fixed-point map, exactness gives
\[
\operatorname{Im}(c_\ast)=\ker i_\ast
=\ker(\Phi^{S^1}\circ i_\ast)=\mathcal{K}_V.
\]
Therefore $c_\ast$ is surjective onto $\mathcal{K}_V$ and hence is an isomorphism.
\end{proof}

\begin{cor} \label{cor: our case injective}
    The assumptions of \Cref{lem: injectivity condition} is satisfied for $V=\R+\H-3\wt\R$. In particular, for this choice of $V$, the map 
    \[
    c_\ast:\pi_V^{\mathrm{Pin}(2)}(E\mathrm{Pin}(2)_+)\rightarrow \mathcal{K}_V
    \]
    is an isomorphism.
\end{cor}
\begin{proof}
    The inequalities in the assumptions are satisfied; it thus suffices to construct a $\mathrm{Pin}(2)$-equivariant map $f:\H^+\rightarrow (3\wt\R)^+$ such that $\deg f^{\mathrm{Pin}(2)}=1$. For that, we simply consider the following map:
    \[
    f(q)=\bar{q}iq:\H\rightarrow \wt\R^3;
    \]
    since $|\bar{q}iq| = |q|^2 \rightarrow \infty$ as $|q|\rightarrow\infty$, we see that $f$ extends to a map from $\H^+$ to $(3\wt\R)^+$. Taking $\mathrm{Pin}(2)$-fixed point then gives the identity map, so clearly we have
    \[
    \deg f^{\mathrm{Pin}(2)}=\deg \mathrm{id} = 1.
    \]
    The lemma follows.
\end{proof}

Suppose that $V$ satisfies the assumptions of \Cref{lem: injectivity condition}. Then, by tom Dieck splitting and Adams' isomorphism (see the parts (i) and (ii) of \Cref{thm: tom dieck and eqv euler characteristic} for the case $V=0$), we see that there exists a natural isomorphism
\[
    \mathcal{A}_V:\mathcal{K}_V \xrightarrow\simeq \pi_n^s (\mathrm{Th}(\xi_{\wt\R}-\xi_V+n\underline{\R})),
\]
where $\xi_W$ is defined as the associated virtual vector bundle (over $B\mathrm{Pin}(2)$) of a virtual $\mathrm{Pin}(2)$-representation. This isomorphism is explicitly given as the following composition of isomorphisms:
\[
    \mathcal{K}_V \xrightarrow{c_\ast^{-1}}\pi_V^{\mathrm{Pin}(2)}(E\mathrm{Pin}(2)_+) \xrightarrow{\Sigma^{n\R-V}} \pi_n^{\mathrm{Pin}(2)}(\Sigma^{n\R-V}E\mathrm{Pin}(2)_+) \xrightarrow{\text{Adams}}\pi_n^s (\mathrm{Th}(\xi_{\wt\R}-\xi_V+n\underline{\R})).
\]

\begin{lem} \label{lem: Adams isomorphism and alpha}
    Let $V$ be a virtual $\mathrm{Pin}(2)$-representation of rank $n$ and $f:V^+\rightarrow S^0$ be a $\mathrm{Pin}(2)$-equivariant stable map. Then $f\alpha \in \mathcal{K}_{V+\R}$. Furthermore, if $V+\R$ satisfies the assumptions of \Cref{lem: injectivity condition}, then $\mathcal{A}_{V+\R}(f\alpha)$ is stably homotopic to the following composition:
    \[
    S^{\dim V+1} \xrightarrow{\Sigma \mathrm{Res}^{\mathrm{Pin}(2)}_{\{1\}} f} S^1 \xrightarrow{e_V} \mathrm{Th}(\xi_{\wt\R}-\xi_V+n\underline{\R}).
    \]
    Here, since $\dim (\wt\R-V+n\R)=1$, we are choosing an identification $\mathrm{Res}^{\mathrm{Pin}(2)}_{\{1\}}(\wt\R-V+n\R)\simeq \R$ and using it to define the fiber-over-the-basepoint (nonequivariant) inclusion\footnote{
    This is determined up to sign (i.e. orientation reversal), which has no effect here as $\alpha$ is a 2-torsion class.}
    \[
    e_V:S^1 \simeq \R^+ \hookrightarrow\mathrm{Th}(\xi_{\wt\R} - \xi_V + n\underline{\R}).
    \]
\end{lem}
\begin{proof}
    Recall from \Cref{rem: rephrasing of definition of alpha} that $\alpha$ can be written as the following composition of $\mathrm{Pin}(2)$-equivariant stable maps:
    \[
    S^1 \xrightarrow{p^!_\mu} \mathrm{Pin}(2)_+ \xrightarrow{\text{collapse}} S^0.
    \]
    Hence $\alpha$ factors through $\mathrm{Pin}(2)_+$. Since $\Phi^{S^1}\mathrm{Pin}(2)_+=\ast$, we see that $\Phi^{S^1}\alpha=0$. Hence we have $\Phi^{S^1}(f\alpha)=0$ and thus $f\alpha\in \mathcal{K}_{V+\R}$. In order to compute $\mathcal{A}_{V+\R}(f\alpha)$, we first observe using \Cref{rem: rephrasing of definition of alpha} that the following composition is its lift to $\pi_{V+\R}^{\mathrm{Pin}(2)}(E\mathrm{Pin}(2)_+)$:
    \[
    (V+\R)^+ \xrightarrow{\mathrm{id}\wedge p^!_\mu} V^+ \wedge \mathrm{Pin}(2)_+ \xrightarrow{f\wedge \mathrm{id}}\mathrm{Pin}(2)_+ \xrightarrow{j_+} E\mathrm{Pin}(2)_+,
    \]
    where $p:\mathrm{Pin}(2)\rightarrow \ast$ is the collapsing map, $\mu$ is the left-invariant framing of $\mathrm{Pin}(2)$, and $j:\mathrm{Pin}(2)\rightarrow E\mathrm{Pin}(2)$ is the inclusion of a free orbit. Suspending by $n\R-V$ then gives the following composition, which we denote by $\widehat{f\alpha}$:
    \[
    S^{n+1}\xrightarrow{\mathrm{id}\wedge p^!_\mu} S^n \wedge \mathrm{Pin}(2)_+ \xrightarrow{\Sigma^{n\R-V}f\wedge \mathrm{id}} \Sigma^{n\R-V}\mathrm{Pin}(2)_+ \xrightarrow{\Sigma^{n\R-V} j_+} \Sigma^{n\R-V}E\mathrm{Pin}(2)_+.
    \]
    Hence, in order to prove the lemma, it suffices to show that
    \[
    \tilde\tau_V\circ e_V\circ \Sigma\mathrm{Res}^{\mathrm{Pin}(2)}_{\{1\}}(f) \sim \widehat{f\alpha},
    \]
    where $\tilde\tau_V:\mathrm{Th}(\xi_{\wt\R}-\xi_V+n\underline{\R})\rightarrow \Sigma^{n\R-V}E\mathrm{Pin}(2)_+$ is the ($\mathrm{Pin}(2)$-equivariant) Adams map for $\Sigma^{n\R-V}E\mathrm{Pin}(2)_+$.

    Observe that the $\mathrm{Pin}(2)$-equivariant Adams map for the free $\mathrm{Pin}(2)$-spectrum $\Sigma^{n\R-V}\mathrm{Pin}(2)_+$ takes the form
    \[
    \tilde\tau_V^0:S^{\mathrm{Res}^{\mathrm{Pin}(2)}_{\{1\}} (n\R-V)}\wedge S^1 \rightarrow \Sigma^{n\R-V}\mathrm{Pin}(2)_+.
    \]
    By the naturality of Adams maps, applied to the $\mathrm{Pin}(2)$-equivariant map 
    \[
    \Sigma^{n\R-V}j_+:\Sigma^{n\R-V}\mathrm{Pin}(2)_+\rightarrow \Sigma^{n\R-V}E\mathrm{Pin}(2)_+,
    \]
    we get the following homotopy-commutative square.
    \[
    \xymatrix{
    S^1 \ar[rr]^{e_V} \ar[d]_{\tilde\tau_V^0} && \mathrm{Th}(\xi_{\wt\R}-\xi_V+n\underline{\R}) \ar[d]^{\tilde\tau_V} \\
    \Sigma^{n\R-V}\mathrm{Pin}(2)_+ \ar[rr]^{\Sigma^{n\R-V} j_+} && \Sigma^{n\R-V}E\mathrm{Pin}(2)_+
    }
    \]
    In particular, we have
    \[
    \tilde\tau_V\circ e_V\circ \Sigma\mathrm{Res}^{\mathrm{Pin}(2)}_{\{1\}}(f) \sim \Sigma^{n\R-V} j_+ \circ \tilde\tau^0_V \circ \Sigma \mathrm{Res}^{\mathrm{Pin}(2)}_{\{1\}} (f).
    \]
    Now, by the definition of Adams map, we see that $\tilde\tau^0_V$ is given, up to suspension, by the composition
    \[
    \mathrm{Res}^{\mathrm{Pin}(2)}_{\{1\}}(n\R-V)^+ \wedge S^1 \xrightarrow{\mathrm{id}\wedge p_\mu^!} (n\R-V)^+ \wedge \mathrm{Pin}(2)_+ \xrightarrow{u_{n\R-V} \wedge \mathrm{id}} (n\R-V)^+ \wedge \mathrm{Pin}(2)_+,
    \]
    where $u_W:\mathrm{Res}^{\mathrm{Pin}(2)}_{\{1\}}W^+ \wedge \mathrm{Pin}(2)_+\rightarrow W^+\wedge \mathrm{Pin}(2)_+$ is given by $u(x,g)=(gx,g)$ for any virtual $G$-representation $W$. Hence we see that $\tilde\tau^0_V \circ \Sigma \mathrm{Res}^{\mathrm{Pin}(2)}_{\{1\}} (f)$ is homotopic to the following composition:
    \[
    (n\R-V)^+ \wedge V^+ \wedge S^1 \xrightarrow{\mathrm{id}\wedge \mathrm{Res}^{\mathrm{Pin}(2)}_{\{1\}}(f)\wedge \mathrm{id}} (n\R-V)^+ \wedge S^1 \xrightarrow{\mathrm{id}\wedge p^!_\mu} (n\R-V)^+\wedge \mathrm{Pin}(2)_+ \xrightarrow{u_{n\R-V}} (n\R-V)^+\wedge \mathrm{Pin}(2)_+ .
    \]
    Clearly, $u$ homotopy-commutes with any stable $\mathrm{Pin}(2)$-equivariant maps between virtual representation spheres in the following sense:
    \[
    u_{V'}\circ (\mathrm{Res}^{\mathrm{Pin}(2)}_{\{1\}} F \wedge \mathrm{id}_{\mathrm{Pin}(2)_+}) \sim (F \wedge \mathrm{id}_{\mathrm{Pin}(2)_+})\circ u_V \quad \text{for}\quad F:V^+ \rightarrow (V')^+.
    \]
    So the above map can be re-written as
    \[
    (n\R)^+ \wedge S^1 \xrightarrow{\mathrm{id}\wedge p^!_\mu} (n\R)^+ \wedge \mathrm{Pin}(2)_+ \xrightarrow{u_{n\R}} (n\R)^+ \wedge \mathrm{Pin}(2)_+ \xrightarrow {\Sigma^{n\R-V}f \wedge \mathrm{id}}(n\R-V)^+ \wedge \mathrm{Pin}(2)_+.
    \]
    Since $n\R$ is a trivial $\mathrm{Pin}(2)$-representation and thus $u_{n\R}\sim \mathrm{id}$, we deduce that postcomposing this by $\Sigma^{n\R-V} j_+$ gives $\widehat{f\alpha}$. Therefore we have
    \[
    \tilde\tau_V\circ e_V\circ \Sigma\mathrm{Res}^{\mathrm{Pin}(2)}_{\{1\}}(f) \sim \widehat{f\alpha},
    \]
    as desired.
\end{proof}

\begin{lem} \label{lem: nonvanishing K3 class}
    Let $\eta\in \pi^{\mathrm{Pin}(2)}_{\H-3\wt\R} (S^0)$ be any element that is nonequivariantly stably homotopic to the unique generator (i.e. Hopf element) $\eta$ of $\pi_1^s S^0 \simeq \Z_2$. Then $\eta\circ\alpha \in \pi_{\R+\H-3\wt\R}^{\mathrm{Pin}(2)}(S^0)$ is not (stably) $\mathrm{Pin}(2)$-equivariantly nullhomotopic.
\end{lem}
\begin{proof}
    By \Cref{cor: our case injective,lem: Adams isomorphism and alpha}, it suffices to show that the (nonequivariant) stable map
    \[
    e_V\circ \Sigma \eta :S^2 \rightarrow \mathrm{Th}(\xi_{4\wt\R-\H}+\underline{\R})
    \]
    is not stably nullhomotopic. Observe that $B\mathrm{Pin}(2)$ has a CW-complex structure that has cells in every dimension except $4k+3$ for all $k\ge 0$. In particular, its 3-skeleton is the same as its 2-skeleton, which is easily seen to be $\R P^2$, as also observed in the proof of \cite[Proposition 4.10]{lin2025family}. Since both $\xi_{4\wt\R}$ and $\xi_\H$ admit quaternionic line bundle structures and every quaternionic line bundle over $\R P^2$ is trivial, we see that
    \[
    \mathrm{Th}((\xi_{4\wt\R-\H}+\underline{\R})\vert_{\R P^2}) \simeq \mathrm{Th}(\underline{\R}_{\R {P^2}}) \simeq \Sigma \R P^2_+.
    \]
    Hence we see that the inclusion $\Sigma \R P^2_+ \simeq \mathrm{Th}((\xi_{4\wt\R-\H}+\underline{\R})\vert_{\R P^2})\hookrightarrow \mathrm{Th}(\xi_{4\wt\R-\H}+\underline{\R})$ induces an isomorphism in $\pi_2^s$.

    Next, observe that the inclusion $e_V:S^1 \rightarrow \mathrm{Th}(\xi_{4\wt\R-\H}+\underline{\R})$ factors as follows:
    \[
    S^1 \xrightarrow{e_V^0} \Sigma \R P^2_+ \simeq  \mathrm{Th}((\xi_{4\wt\R-\H}+\underline{\R})\vert_{\R P^2}) \hookrightarrow \mathrm{Th}(\xi_{4\wt\R-\H}+\underline{\R}).
    \]
    The map $e_V^0:S^1 \rightarrow \Sigma \R P^2_+$ clearly has a stable homotopy left inverse, given by the suspension of the collapsing map $\R P^2_+ \rightarrow S^0$. Hence we see that it induces injective maps on all stable homotopy groups. Thus it follows that the map
    \[
    e_V\circ -:\pi_2 ^s (S^1)\rightarrow \pi_2 ^s \mathrm{Th}(\xi_{4\wt\R-\H}+\underline{\R})
    \]
    is also injective. Since $\Sigma \eta \in \pi_2 ^s (S^1)$ is not stably nullhomotopic, we deduce that $e_V \circ \Sigma \eta$ is also not stably nullhomotopic. Therefore $\eta \circ \alpha$ is not stably $\mathrm{Pin}(2)$-equivariantly nullhomotopic.
\end{proof}

We are now ready to prove \Cref{thm: K3}, whose statement we recall.

\MainKThreeTheorem*

\begin{proof}
Let $Z$ be any smooth homotopy K3 surface.
By the proof of \Cref{thm: main}, there exists
$f\in \mathrm{Diff}^+(S^4)$ supported in the interior
of a standard $4$-ball such that
$\mathfrak{LM}_{S^4}(f)=1$.
In particular, $f$ fixes an open ball
$U\subset S^4$ pointwise.

Thus, by \Cref{lem: connected sum FBF}, we have the connected sum diffeomorphism
\[
f\# \mathrm{id} \in \mathrm{Diff}^+(S^4 \# Z) = \mathrm{Diff}^+(Z),
\]
which admits a spin lift $F$ such that the reduced component $FBF(F)$
of the family Bauer--Furuta invariant $BF_{\mathrm{Torus}(F)}$ of the
spin mapping torus $\mathrm{Torus}(F)\rightarrow S^1$ is given by
\[
FBF(F) = BF_{Z,\s_Z}\circ \alpha,
\]
where $\s_Z$ denotes the unique spin structure on $Z$.

Recall that the identity map on $Z$ has two spin lifts, namely
$\mathrm{id}_{Z}$ and $-\mathrm{id}_{Z}$.
We trivialize $H^+$ as in \cite[Section~3]{Kronheimer-Mrowka:2020-1},
using a fixed oriented maximal positive subspace $P\subset H^2(Z;\R)$
for the intersection form. The mapping tori of $F$, $\mathrm{id}_Z$,
and $-\mathrm{id}_Z$ have trivial monodromy on $H^2(Z;\R)$, so the
projection $P\oplus P^\perp\to P$, $p+q\mapsto p$, identifies each
$H^+$-bundle with $S^1\times P$.
We know from \cite[Proposition 4.1]{Kronheimer-Mrowka:2020-1} that
the reduced parts $FBF(\mathrm{id}_{Z})$ and $FBF(-\mathrm{id}_{Z})$
of the family Bauer--Furuta invariants of spin mapping tori
$\mathrm{Torus}(\mathrm{id}_{Z})$ and $\mathrm{Torus}(-\mathrm{id}_{Z})$
satisfies
\[
FBF(\mathrm{id}_{Z}) \sim_{\mathrm{Pin}(2)} 0,\quad
FBF(-\mathrm{id}_{Z}) \sim \eta_0^2,
\]
where $\eta_0 \in \pi_1^s(S^0)$ is the nonequivariant stable Hopf element.
Furthermore, it is shown in \cite[Section 6]{Kronheimer-Mrowka:2020-1}
that (nonequivariantly) $BF_{Z,\s_Z} \sim \eta$; see
\Cref{lem: nonvanishing K3 class} for the definition of $\eta$.

Suppose that $f\# \mathrm{id}$ is smoothly isotopic to identity.
Then its spin lift $F$ is smoothly isotopic to either $\mathrm{id}_{Z}$
or $-\mathrm{id}_{Z}$. Taking the reduced part of the family
Bauer--Furuta invariants of their spin mapping tori then gives
\[
BF_{Z,\s_Z}\circ\alpha\sim_{\mathrm{Pin}(2)} 0
\quad \text{or}\quad
BF_{Z,\s_Z}\circ\alpha \sim \eta_0^2
\]
Since $\eta_0^2$ generates $\pi_2^s(S^0)\cong \Z_2$ and $\alpha$ is
nonequivariantly nullhomotopic, we see that $BF_{Z,\s_Z}\circ\alpha$
cannot be nonequivariantly homotopic to $\eta_0^2$, and thus it should
be $\mathrm{Pin}(2)$-equivariantly nullhomotopic, contradicting
\Cref{lem: nonvanishing K3 class}. The theorem follows.
\end{proof}

\subsection{Lin--Mukherjee invariant of half-unknotted Montesinos twists}

In order to compute the Lin--Mukherjee invariant of the exotic diffeomorphism (which we denote by $f_{GH}$) of $S^4$ claimed in \cite{gabai2025pseudo}, we start by recalling the definition of Montesinos twist diffeomorphisms in \cite{gay2025relations}. A \emph{Montesinos twin} in $S^4$ is an ordered pair $W=(R,S)$
of oriented smoothly embedded $2$-spheres that meet transversely
in two points. Since $H_2(S^4;\Z)=0$, both normal bundles are
trivial and the two intersection points have opposite signs.
We call $W$ \emph{half-unknotted} if either $R$ or $S$ bounds a
smoothly embedded $3$-ball in $S^4$.

Let $\nu(W)$ be a closed regular neighborhood of $R\cup S$.
Following \cite[Definitions~1--3]{gay2025relations}, write
\[
\partial\nu(W)\cong S^1_\ell\times S^1_R\times S^1_S,
\]
where $S^1_R$ and $S^1_S$ are the positive meridians of $R$ and $S$,
respectively, and the longitude $S^1_\ell$ represents a generator of
\[
\ker\!\left(
H_1(\partial\nu(W);\Z)\longrightarrow
H_1(S^4\setminus\operatorname{int}\nu(W);\Z)
\right).
\]
Orient $S^1_\ell$ so that the product orientation agrees with the
boundary orientation of $\nu(W)$. Take a collar of $\partial\nu(W)$
with coordinate $t\in[-1,1]$ increasing outwards from $\nu(W)$.
On this collar, the \emph{twist along $W$}, denoted by $\tau_W$,
is the product of a right-handed Dehn twist on $[-1,1]\times S^1_\ell$
with the identity on $S^1_R\times S^1_S$. Extending by the identity
over the rest of $S^4$ defines $\tau_W$ up to smooth isotopy.

Consider the subgroup $\mathcal{M}_0 \subset \pi_0 \mathrm{Diff}^+(S^4)$ generated by twists along half-unknotted Montesinos twins in $S^4$. It is proven in \cite[Theorem 4]{gay2025relations} that $\mathcal{M}_0$ is either trivial or isomorphic to $\Z_2$, and generated by the twist $\tau_{W(1)}$; see \cite[Figure 1]{gay2025relations} for the definition of $W(i)$ for each integer $i$. In fact, $\tau_{W(1)}$ is exactly the diffeomorphism $f_{GH}$ that is claimed to be not smoothly isotopic to identity in \cite[Theorem 0.11]{gabai2025pseudo}.

By \cite[Lemma 6]{gay2025relations}, we know that for any half-unknotted Montesinos twin $W$ in $S^4$, there exists an associated loop $\alpha_W = \{\alpha_t\}_{t\in [0,1]}$ of smooth embeddings $\alpha_t:S^1 \hookrightarrow S^1 \times S^3$, based at $S^1 \times \{p\}$ for a chosen basepoint $p\in S^3$, such that $\tau_W$ can be recovered from $\alpha_W$ in the following way.
\begin{itemize}
    \item Extend $\{\alpha_t\}_{t\in [0,1]}$ to a loop of framed embeddings;
    \item Choose an ambient isotopy $\{\psi_t\}_{t\in [0,1]}$ realizing this loop of framed embeddings;
    \item Performing surgeries along each framed loop $\alpha_t$ turns each $\psi_t$ into a diffeomorphism $\psi'_t$ from the surgered fiber at $t=0$, i.e. $S^4$, to the surgered fiber at each $t$;
    \item Then the orientation-preserving self-diffeomorphism $\mathcal{H}(\alpha_W):=\psi_1'$ of $S^4$ is well-defined up to smooth isotopy and depends only on the based homotopy class of $\alpha_W$ in $\mathrm{Emb}(S^1,S^1\times S^3)$;
    \item The Montesinos twist $\tau_W$ is smoothly isotopic to $\mathcal{H}(\alpha_W)$.
\end{itemize}

\begin{prop} \label{prop: Montesinos twist vanishing}
    Let $\alpha = \{\alpha_t\}$ be a loop of smooth embeddings $\alpha_t:S^1 \hookrightarrow S^1 \times S^3$ based at $S^1 \times \{p\}$. Then $\mathfrak{LM}_{S^4}(\mathcal{H}(\alpha))=0$. Consequently, for any half-unknotted Montesinos twin $W$ in $S^4$, we have $\mathfrak{LM}_{S^4}(\tau_W)=0$.
\end{prop}
\begin{proof}
    Consider the trivial smooth bundle $E :=S^1 \times (S^1 \times S^3) \rightarrow S^1$. Then the loop $\alpha$ defines a family $C_\alpha\subset E$ of framed smoothly embedded loops, chosen as in the construction of $\mathcal{H}(\alpha)$. The oriented smooth $S^4$-bundle over $S^1$ obtained by performing the fiberwise surgery of $E$ along $C_\alpha$ is orientation-preserving smoothly isomorphic to the mapping torus family $\mathrm{Torus}(\mathcal{H}(\alpha))\rightarrow S^1$ of the orientation-preserving self-diffeomorphism $\mathcal{H}(\alpha)$.
    
    Clearly, $S^1 \times S^3$ admits a psc metric, which implies that $E$ has a fiberwise psc metric. Since $C_\alpha$ has codimension 3 in $E$, it follows from \cite[Theorem 3.10]{Walsh2011} and \cite[Theorem 1.1]{Chernysh2004} that the fiberwise surgery of $E$ along $C_\alpha$, and thus the mapping torus family $\mathrm{Torus}(\mathcal{H}(\alpha))\rightarrow S^1$, also admits a fiberwise psc metric. Therefore we deduce that $\mathfrak{LM}_{S^4}(\mathcal{H}(\alpha))=0$ by \Cref{cor: fiberwise psc metric implies vanishing of LM}.
\end{proof}

\begin{cor} \label{cor: GGH vanishing}
    We have $\mathfrak{LM}_{S^4}(f_{GH})=0$.
\end{cor}
\begin{proof}
    As discussed above, the diffeomorphism $f_{GH}$ is given by the twist $\tau_{W(1)}$ along the Montesinos twin $W(1)$ in $S^4$. Since $W(1)$ is half-unknotted\footnote{In fact, for any integer $i$, both spheres in the Montesinos twin $W(i)$ are unknotted.}, we get $\mathfrak{LM}_{S^4}(\tau_{W(1)})=0$, as desired.
\end{proof}

\appendix

\section{An explicit proof of $T(e)=\alpha$} \label{sec: appendix A}

In this section, we give an alternative proof of \Cref{lem: computation of T}, using explicitly defined homotopies. This gives a more geometric, rather than functorial, understanding of the equality $T(e)=\alpha$.

\begin{proof}[Proof of $T(e)=\alpha$ in \Cref{lem: computation of T}]
In the coordinates of \Cref{table: universe summands},
$\tilde\tau(x)=x(-j)$. Under the $G_\mathrm{odd}$-action,
the point $1\in\H$ has orbit $\mathrm{Pin}(2)$ and stabilizer
$\langle I\rangle$.
Fix $0<\delta<1/4$ and write $x=g(1+r)$ in the tubular
neighborhood of $\mathrm{Pin}(2)$, where $g\in\mathrm{Pin}(2)$,
$r\perp i$, and $|r|<\delta$.
Let $c\colon\H\to\H^+$ be the radial collapse
\[
 c(v)=
 \begin{cases}
 v/(\delta^2-|v|^2), & |v|<\delta,\\
 \infty, & |v|\ge\delta.
 \end{cases}
\]

Write $S^1=[0,1]/(0\sim1)$, based at $[0]$. For $s,t\in[0,1]$,
define $H_{s,t}\colon\H^+\to\H^+$ by
\[
 H_{s,t}\bigl(g(1+r)\bigr)
 =g\,c\bigl(r+s(t-\tfrac12)i\bigr), \quad g \in \mathrm{Pin}(2), \ r \in \R\langle 1, j,k\rangle 
\]
on the tube and by $\infty$ elsewhere. Since $r\perp i$, the
formula extends continuously across the boundary of the tube.
Conjugation by $-j$ preserves $i^\perp$ and sends $i$ to $-i$.
Thus radiality of $c$ gives
\[
 H_{s,1}(x)(-j)=H_{s,0}\bigl(x(-j)\bigr),
\]
so $H$ defines a $\mathrm{Pin}(2)$-equivariant fiberwise homotopy
on the mapping torus of $\tilde\tau$. At $s=0$ the fiber map
$g(1+r)\mapsto g\,c(r)$ is the Pontryagin--Thom representative
of $e=\chi_{G_\mathrm{odd}}(G_\mathrm{odd}/\langle I\rangle)$
obtained by inserting zero in the tangent direction. Its
mapping-torus family therefore has reduced stable component $T(e)$.

At $s=1$ the map takes finite values precisely where
\[
 |r|^2+(t-\tfrac12)^2<\delta^2,
\]
as shown on the slice $r=u\cdot1$ in \Cref{fig:te-alpha-homotopy}.
Choose a path $b_t$ in $\mathrm{Sp}(1)$ from $1$ to $-j$ with
$b_t=1$ for $t\le 3/4$. The trivialization
$[t,x]\mapsto([t],xb_t)$ is then the identity on this region,
and the map is constant at $\infty$ near $[0]$.
With $a=t-\tfrac12$, the induced map $S^1\wedge\H^+\to\H^+$
is the Pontryagin--Thom map of
$\{[\tfrac12]\}\times\mathrm{Pin}(2)$ with normal framing
\[
 (a,gr)\longmapsto g(r+ai).
\]
This framing uses the left-invariant tangent vector $gi$ at $g$,
so the map represents $\alpha$ by \cite[Equation~(32)]{lin2025family}.
Taking reduced stable components gives $T(e)=\alpha$.
\end{proof}
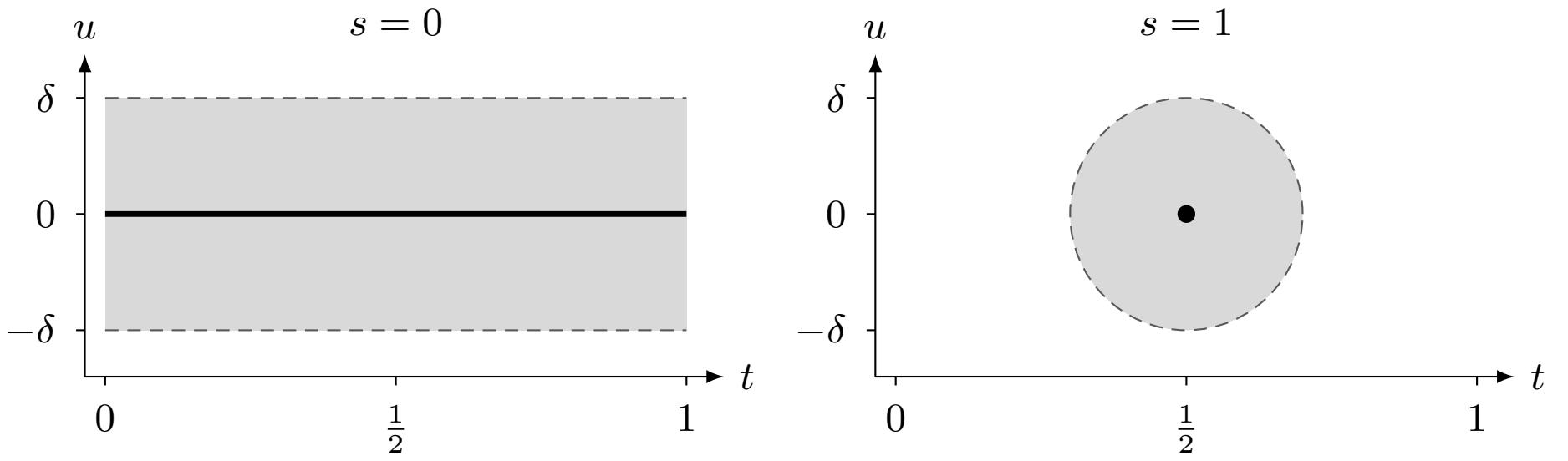
\begin{figure}[htbp]
\centering
\begin{tikzpicture}[
  x=4.6cm,y=4.6cm,font=\small,>=latex,
  edge of tube/.style={draw=black!65,thin,densely dashed}
]
  \fill[black!15] (0,-.2) rectangle (1,.2);
  \draw[edge of tube] (0,-.2)--(1,-.2) (0,.2)--(1,.2);
  \draw[line width=1.3pt] (0,0)--(1,0);

  \begin{scope}[shift={(1.36,0)}]
    \filldraw[fill=black!15,edge of tube]
      (.5,0) circle[radius=.2];
    \fill (.5,0) circle[radius=2pt];
  \end{scope}

  \foreach \xoffset/\stageindex in {0/0,1.36/1} {
    \begin{scope}[shift={(\xoffset,0)}]
      \node at (.5,.33) {$s=\stageindex$};
      \draw[->] (-.035,-.28)--(1.065,-.28) node[right] {$t$};
      \draw[->] (-.035,-.28)--(-.035,.275) node[above] {$u$};
      \foreach \tickpos/\ticktext in {0/0,.5/{\tfrac12},1/1}
        \draw (\tickpos,-.28)--(\tickpos,-.295)
          node[below=1pt] {$\ticktext$};
      \foreach \tickpos/\ticktext in {-.2/{-\delta},0/0,.2/\delta}
        \draw (-.035,\tickpos)--(-.05,\tickpos)
          node[left=1pt] {$\ticktext$};
    \end{scope}
  }
\end{tikzpicture}
\caption{On the slice $r=u\cdot1$, each point represents a
$\mathrm{Pin}(2)$-orbit. Shading indicates finite values.
The black segment represents one orbit of zeros in each fiber,
with the ends identified by $g\mapsto g(-j)$.
The dot represents the single orbit of zeros at $t=\tfrac12$.}
\label{fig:te-alpha-homotopy}
\end{figure}

\bibliographystyle{alpha}
\bibliography{tex}

@article {adams1962vector,
    AUTHOR = {Adams, John F.},
     TITLE = {Vector fields on spheres},
   JOURNAL = {Ann. of Math. (2)},
  FJOURNAL = {Annals of Mathematics. Second Series},
    VOLUME = {75},
      YEAR = {1962},
     PAGES = {603--632},
      ISSN = {0003-486X},
   MRCLASS = {57.30},
  MRNUMBER = {139178},
MRREVIEWER = {M.\ F.\ Atiyah},
       DOI = {10.2307/1970213},
       URL = {https://doi.org/10.2307/1970213},
}

@incollection {Aitchison-Rubinstein:1984,
    AUTHOR = {Aitchison, Iain R. and Rubinstein, Joachim H.},
     TITLE = {Fibered knots and involutions on homotopy spheres},
 BOOKTITLE = {Four-manifold theory ({D}urham, {N}.{H}., 1982)},
    SERIES = {Contemp. Math.},
    VOLUME = {35},
     PAGES = {1--74},
 PUBLISHER = {Amer. Math. Soc., Providence, RI},
      YEAR = {1984},
      ISBN = {0-8218-5033-4},
   MRCLASS = {57N13 (57Q45)},
  MRNUMBER = {780575},
MRREVIEWER = {Selman\ Akbulut},
       DOI = {10.1090/conm/035/780575},
       URL = {https://doi.org/10.1090/conm/035/780575},
}

@article {Akbulut-Kirby:1979,
    AUTHOR = {Akbulut, Selman and Kirby, Robion C.},
     TITLE = {An exotic involution of {$S\sp{4}$}},
   JOURNAL = {Topology},
  FJOURNAL = {Topology. An International Journal of Mathematics},
    VOLUME = {18},
      YEAR = {1979},
    NUMBER = {1},
     PAGES = {75--81},
      ISSN = {0040-9383},
   MRCLASS = {57N15},
  MRNUMBER = {528237},
MRREVIEWER = {A.\ Chernavski\u i},
       DOI = {10.1016/0040-9383(79)90015-6},
       URL = {https://doi.org/10.1016/0040-9383(79)90015-6},
}

@article {Akbulut-Kirby:1985,
    AUTHOR = {Akbulut, Selman and Kirby, Robion C.},
     TITLE = {A potential smooth counterexample in dimension {$4$} to the
              {P}oincar\'e{} conjecture, the {S}choenflies conjecture, and
              the {A}ndrews-{C}urtis conjecture},
   JOURNAL = {Topology},
  FJOURNAL = {Topology. An International Journal of Mathematics},
    VOLUME = {24},
      YEAR = {1985},
    NUMBER = {4},
     PAGES = {375--390},
      ISSN = {0040-9383},
   MRCLASS = {57R60 (57R10)},
  MRNUMBER = {816520},
MRREVIEWER = {Martin\ Scharlemann},
       DOI = {10.1016/0040-9383(85)90010-2},
       URL = {https://doi.org/10.1016/0040-9383(85)90010-2},
}

@misc{Lin-Wu-Xie-Zhang:2025,
  author        = {Lin, Jianfeng and Wu, Weiwei and Xie, Yi
                   and Zhang, Boyu},
  title         = {{Dax} invariants, light bulbs, and isotopies
                   of symplectic structures},
  year          = {2025},
  note          = {Preprint, arXiv:2501.16083},
  eprint        = {2501.16083},
  archivePrefix = {arXiv},
  primaryClass  = {math.GT},
  url           = {https://arxiv.org/abs/2501.16083}
}

@unpublished{Kronheimer:families,
  author        = {Kronheimer, P. B.},
  title         = {Some non-trivial families of symplectic structures},
  note          = {Preprint},
  url           = {https://people.math.harvard.edu/~kronheim/diffsymp.pdf}
}

@misc{Watanabe:2018,
  author        = {Watanabe, Tadayuki},
  title         = {Some exotic nontrivial elements of the rational
                   homotopy groups of {$\mathrm{Diff}(S^4)$}},
  year          = {2018},
  eprint        = {1812.02448},
  archivePrefix = {arXiv},
  primaryClass  = {math.GT},
  note          = {Version 3, 2019}
}

@article{Wang:2026,
  author  = {Wang, Luya},
  title   = {A note on a nontrivial connected component of the
             space of symplectic structures},
  journal = {Algebr. Geom. Topol.},
  volume  = {26},
  number  = {7},
  year    = {2026},
  pages   = {2705--2708}
}

@article{Gromov:1985,
  author  = {Gromov, M.},
  title   = {Pseudo holomorphic curves in symplectic manifolds},
  journal = {Invent. Math.},
  volume  = {82},
  number  = {2},
  year    = {1985},
  pages   = {307--347},
  doi     = {10.1007/BF01388806}
}

@article {akbulut2010cappell,
    AUTHOR = {Akbulut, Selman},
     TITLE = {Cappell-{S}haneson homotopy spheres are standard},
   JOURNAL = {Ann. of Math. (2)},
  FJOURNAL = {Annals of Mathematics. Second Series},
    VOLUME = {171},
      YEAR = {2010},
    NUMBER = {3},
     PAGES = {2171--2175},
      ISSN = {0003-486X,1939-8980},
   MRCLASS = {57N13},
  MRNUMBER = {2680408},
MRREVIEWER = {Martin\ Scharlemann},
       DOI = {10.4007/annals.2010.171.2171},
       URL = {https://doi.org/10.4007/annals.2010.171.2171},
}

@article {atiyah1968lefschetz,
    AUTHOR = {Atiyah, Michael F. and Bott, Raoul},
     TITLE = {A {L}efschetz fixed point formula for elliptic complexes.
              {II}. {A}pplications},
   JOURNAL = {Ann. of Math. (2)},
  FJOURNAL = {Annals of Mathematics. Second Series},
    VOLUME = {88},
      YEAR = {1968},
     PAGES = {451--491},
      ISSN = {0003-486X},
   MRCLASS = {57.50},
  MRNUMBER = {232406},
MRREVIEWER = {R.\ S.\ Palais},
       DOI = {10.2307/1970721},
       URL = {https://doi.org/10.2307/1970721},
}

@article {Auckly-Kim-Melvin-Ruberman:2015,
    AUTHOR = {Auckly, David R. and Kim, Hee Jung and Melvin, Paul M. and Ruberman, Daniel},
     TITLE = {Stable isotopy in four dimensions},
   JOURNAL = {J. Lond. Math. Soc. (2)},
  FJOURNAL = {Journal of the London Mathematical Society. Second Series},
    VOLUME = {91},
      YEAR = {2015},
    NUMBER = {2},
     PAGES = {439--463},
      ISSN = {0024-6107,1469-7750},
   MRCLASS = {57R52 (57Q60 57R40 57R65)},
  MRNUMBER = {3355110},
MRREVIEWER = {Minoru\ Yamamoto},
       DOI = {10.1112/jlms/jdu075},
       URL = {https://doi.org/10.1112/jlms/jdu075},
}

@misc{Auckly-Ruberman:2025,
    AUTHOR = {Auckly, David R. and Ruberman, Daniel},
     TITLE = {Families of diffeomorphisms, embeddings, and positive scalar curvature metrics via {Seiberg--Witten} theory},
      YEAR = {2025},
    EPRINT = {2501.11892},
ARCHIVEPREFIX = {arXiv},
      NOTE = {arXiv:2501.11892},
}

@article {Ba21,
    AUTHOR = {Baraglia, David},
     TITLE = {Constraints on families of smooth 4-manifolds from
              {B}auer-{F}uruta invariants},
   JOURNAL = {Algebr. Geom. Topol.},
  FJOURNAL = {Algebraic \& Geometric Topology},
    VOLUME = {21},
      YEAR = {2021},
    NUMBER = {1},
     PAGES = {317--349},
      ISSN = {1472-2747,1472-2739},
   MRCLASS = {57R57 (57R22 57R50)},
  MRNUMBER = {4224743},
MRREVIEWER = {Dave\ Auckly},
       DOI = {10.2140/agt.2021.21.317},
       URL = {https://doi.org/10.2140/agt.2021.21.317},
}

@misc{bamler2019ricci,
    AUTHOR = {Bamler, Richard H. and Kleiner, Bruce},
     TITLE = {{Ricci} flow and contractibility of spaces of metrics},
      YEAR = {2019},
    EPRINT = {1909.08710},
ARCHIVEPREFIX = {arXiv},
      NOTE = {arXiv:1909.08710},
}

@article {Baraglia-Konno:2020,
    AUTHOR = {Baraglia, David and Konno, Hokuto},
     TITLE = {A gluing formula for families {S}eiberg-{W}itten invariants},
   JOURNAL = {Geom. Topol.},
  FJOURNAL = {Geometry \& Topology},
    VOLUME = {24},
      YEAR = {2020},
    NUMBER = {3},
     PAGES = {1381--1456},
      ISSN = {1465-3060,1364-0380},
   MRCLASS = {57R57 (57M60 57R22)},
  MRNUMBER = {4157556},
MRREVIEWER = {Dave\ Auckly},
       DOI = {10.2140/gt.2020.24.1381},
       URL = {https://doi.org/10.2140/gt.2020.24.1381},
}

@article {Baraglia-Konno:2026,
    AUTHOR = {Baraglia, David and Konno, Hokuto},
     TITLE = {Irreducible 4-manifolds can admit exotic diffeomorphisms},
   JOURNAL = {Duke Math. J.},
  FJOURNAL = {Duke Mathematical Journal},
    VOLUME = {175},
      YEAR = {2026},
    NUMBER = {4},
     PAGES = {717--733},
      ISSN = {0012-7094,1547-7398},
   MRCLASS = {57K41 (57R50 57S05)},
  MRNUMBER = {5053586},
MRREVIEWER = {Marja\ K.\ Kankaanrinta},
       DOI = {10.1215/00127094-2025-0040},
       URL = {https://doi.org/10.1215/00127094-2025-0040},
}

@article{Baraglia-Tomlin:2026,
    AUTHOR = {Baraglia, David and Tomlin, Joshua},
     TITLE = {Exotic diffeomorphisms of 4-manifolds with {$b_+=2$}},
   JOURNAL = {Q. J. Math.},
      YEAR = {2026},
      NOTE = {Published online, article haag031},
       DOI = {10.1093/qmath/haag031},
       URL = {https://doi.org/10.1093/qmath/haag031},
}

@article {baraglia2022bauer,
    AUTHOR = {Baraglia, David and Konno, Hokuto},
     TITLE = {On the {B}auer-{F}uruta and {S}eiberg-{W}itten invariants of
              families of 4-manifolds},
   JOURNAL = {J. Topol.},
  FJOURNAL = {Journal of Topology},
    VOLUME = {15},
      YEAR = {2022},
    NUMBER = {2},
     PAGES = {505--586},
      ISSN = {1753-8416,1753-8424},
   MRCLASS = {57K41 (57R22 57R57)},
  MRNUMBER = {4441598},
       DOI = {10.1112/topo.12229},
       URL = {https://doi.org/10.1112/topo.12229},
}

@misc{baraglia2026exotic,
    AUTHOR = {Baraglia, David},
     TITLE = {Exotic diffeomorphisms of reducible {$4$}-manifolds with odd {$b_+$}},
      YEAR = {2026},
    EPRINT = {2606.31295},
ARCHIVEPREFIX = {arXiv},
      NOTE = {arXiv:2606.31295},
}

@article {Baraglia:2026-mapping,
    AUTHOR = {Baraglia, David},
     TITLE = {On the mapping class groups of simply connected smooth
              4-manifolds},
   JOURNAL = {Algebr. Geom. Topol.},
  FJOURNAL = {Algebraic \& Geometric Topology},
    VOLUME = {26},
      YEAR = {2026},
    NUMBER = {5},
     PAGES = {1635--1653},
      ISSN = {1472-2747,1472-2739},
   MRCLASS = {57K41 (57R50)},
  MRNUMBER = {5079109},
MRREVIEWER = {Allan\ Edmonds},
       DOI = {10.2140/agt.2026.26.1635},
       URL = {https://doi.org/10.2140/agt.2026.26.1635},
}

@article {BECKER19751,
    AUTHOR = {Becker, James C. and Gottlieb, Daniel H.},
     TITLE = {The transfer map and fiber bundles},
   JOURNAL = {Topology},
  FJOURNAL = {Topology. An International Journal of Mathematics},
    VOLUME = {14},
      YEAR = {1975},
     PAGES = {1--12},
      ISSN = {0040-9383},
   MRCLASS = {55E50 (55F10)},
  MRNUMBER = {377873},
MRREVIEWER = {J.\ P.\ May},
       DOI = {10.1016/0040-9383(75)90029-4},
       URL = {https://doi.org/10.1016/0040-9383(75)90029-4},
}

@article {BERW17,
    AUTHOR = {Botvinnik, Boris and Ebert, Johannes and Randal-Williams,
              Oscar},
     TITLE = {Infinite loop spaces and positive scalar curvature},
   JOURNAL = {Invent. Math.},
  FJOURNAL = {Inventiones Mathematicae},
    VOLUME = {209},
      YEAR = {2017},
    NUMBER = {3},
     PAGES = {749--835},
      ISSN = {0020-9910,1432-1297},
   MRCLASS = {53C27 (19D06 19K56 55P47 55R35 55S35 57R22 58J20)},
  MRNUMBER = {3681394},
MRREVIEWER = {Peter\ Hochs},
       DOI = {10.1007/s00222-017-0719-3},
       URL = {https://doi.org/10.1007/s00222-017-0719-3},
}

@article {BF04,
    AUTHOR = {Bauer, Stefan and Furuta, Mikio},
     TITLE = {A stable cohomotopy refinement of {S}eiberg-{W}itten
              invariants. {I}},
   JOURNAL = {Invent. Math.},
  FJOURNAL = {Inventiones Mathematicae},
    VOLUME = {155},
      YEAR = {2004},
    NUMBER = {1},
     PAGES = {1--19},
      ISSN = {0020-9910,1432-1297},
   MRCLASS = {57R57 (55Q55 57R55)},
  MRNUMBER = {2025298},
MRREVIEWER = {Vicente\ Mu\~noz},
       DOI = {10.1007/s00222-003-0288-5},
       URL = {https://doi.org/10.1007/s00222-003-0288-5},
}

@book {BKR2026,
     TITLE = {K3---a new problem list in low-dimensional topology},
    SERIES = {Mathematical Surveys and Monographs},
    VOLUME = {295},
    EDITOR = {Baykur, Refik {\.I}. and Kirby, Robion C. and Ruberman, Daniel},
 PUBLISHER = {American Mathematical Society, Providence, RI},
      YEAR = {2026},
     PAGES = {xi+430},
      ISBN = {978-1-4704-8433-0},
   MRCLASS = {57-02 (57Kxx 57Mxx)},
  MRNUMBER = {5077826},
       DOI = {10.1090/surv/295},
       URL = {https://doi.org/10.1090/surv/295},
}

@article {cappellshaneson76,
    AUTHOR = {Cappell, Sylvain E. and Shaneson, Julius L.},
     TITLE = {Some new four-manifolds},
   JOURNAL = {Ann. of Math. (2)},
  FJOURNAL = {Annals of Mathematics. Second Series},
    VOLUME = {104},
      YEAR = {1976},
    NUMBER = {1},
     PAGES = {61--72},
      ISSN = {0003-486X},
   MRCLASS = {57D65 (57E25)},
  MRNUMBER = {418125},
MRREVIEWER = {Allan\ Edmonds},
       DOI = {10.2307/1971056},
       URL = {https://doi.org/10.2307/1971056},
}

@article {carmeli2025characters,
    AUTHOR = {Carmeli, Shachar and Cnossen, Bastiaan and Ramzi, Maxime and
              Yanovski, Lior},
     TITLE = {Characters and transfer maps via categorified traces},
   JOURNAL = {Forum Math. Sigma},
  FJOURNAL = {Forum of Mathematics. Sigma},
    VOLUME = {13},
      YEAR = {2025},
     PAGES = {Paper No. e93, 84},
      ISSN = {2050-5094},
   MRCLASS = {55P42 (18N65)},
  MRNUMBER = {4914895},
MRREVIEWER = {Geoffrey\ M. L. Powell},
       DOI = {10.1017/fms.2025.23},
       URL = {https://doi.org/10.1017/fms.2025.23},
}

@article {carr1988construction,
    AUTHOR = {Carr, Rodney},
     TITLE = {Construction of manifolds of positive scalar curvature},
   JOURNAL = {Trans. Amer. Math. Soc.},
  FJOURNAL = {Transactions of the American Mathematical Society},
    VOLUME = {307},
      YEAR = {1988},
    NUMBER = {1},
     PAGES = {63--74},
      ISSN = {0002-9947,1088-6850},
   MRCLASS = {53C20},
  MRNUMBER = {936805},
MRREVIEWER = {Viktor\ Schroeder},
       DOI = {10.2307/2000751},
       URL = {https://doi.org/10.2307/2000751},
}

@misc{Chernysh2004,
    AUTHOR = {Chernysh, Vladislav},
     TITLE = {On the homotopy type of the space {$\mathcal{R}^{+}(M)$}},
      YEAR = {2004},
    EPRINT = {math/0405235},
ARCHIVEPREFIX = {arXiv},
      NOTE = {arXiv:math/0405235},
}

@book {CrabbJames1998,
    AUTHOR = {Crabb, Michael C. and James, Ioan M.},
     TITLE = {Fibrewise homotopy theory},
    SERIES = {Springer Monographs in Mathematics},
 PUBLISHER = {Springer-Verlag London, Ltd., London},
      YEAR = {1998},
     PAGES = {viii+341},
      ISBN = {1-85233-014-7},
   MRCLASS = {55-02 (55M20 55P42 55R99 55T10 55T20)},
  MRNUMBER = {1646248},
MRREVIEWER = {J.\ P. C. Greenlees},
       DOI = {10.1007/978-1-4471-1265-5},
       URL = {https://doi.org/10.1007/978-1-4471-1265-5},
}

@article {cs1976inequivalent,
    AUTHOR = {Cappell, Sylvain E. and Shaneson, Julius L.},
     TITLE = {There exist inequivalent knots with the same complement},
   JOURNAL = {Ann. of Math. (2)},
  FJOURNAL = {Annals of Mathematics. Second Series},
    VOLUME = {103},
      YEAR = {1976},
    NUMBER = {3},
     PAGES = {349--353},
      ISSN = {0003-486X},
   MRCLASS = {57C45},
  MRNUMBER = {413117},
MRREVIEWER = {J.\ P.\ Levine},
       DOI = {10.2307/1970942},
       URL = {https://doi.org/10.2307/1970942},
}

@article {ebert2022positive,
    AUTHOR = {Ebert, Johannes and Randal-Williams, Oscar},
     TITLE = {The positive scalar curvature cobordism category},
   JOURNAL = {Duke Math. J.},
  FJOURNAL = {Duke Mathematical Journal},
    VOLUME = {171},
      YEAR = {2022},
    NUMBER = {11},
     PAGES = {2275--2406},
      ISSN = {0012-7094,1547-7398},
   MRCLASS = {53C29 (19K56 55P47 57R90 58D17)},
  MRNUMBER = {4466049},
MRREVIEWER = {Lorenz\ J.\ Schwachh\"ofer},
       DOI = {10.1215/00127094-2022-0023},
       URL = {https://doi.org/10.1215/00127094-2022-0023},
}

@article {fausk2008survey,
    AUTHOR = {Fausk, Halvard},
     TITLE = {Survey on the {B}urnside ring of compact {L}ie groups},
   JOURNAL = {J. Lie Theory},
  FJOURNAL = {Journal of Lie Theory},
    VOLUME = {18},
      YEAR = {2008},
    NUMBER = {2},
     PAGES = {351--368},
      ISSN = {0949-5932},
   MRCLASS = {19A22},
  MRNUMBER = {2431120},
}

@article {Fu01,
    AUTHOR = {Furuta, Mikio},
     TITLE = {Monopole equation and the {$\frac{11}8$}-conjecture},
   JOURNAL = {Math. Res. Lett.},
  FJOURNAL = {Mathematical Research Letters},
    VOLUME = {8},
      YEAR = {2001},
    NUMBER = {3},
     PAGES = {279--291},
      ISSN = {1073-2780},
   MRCLASS = {57R15 (53C27 57R57)},
  MRNUMBER = {1839478},
       DOI = {10.4310/MRL.2001.v8.n3.a5},
       URL = {https://doi.org/10.4310/MRL.2001.v8.n3.a5},
}

@misc{gabai2025pseudo,
    AUTHOR = {Gabai, David and Gay, David T. and Hartman, Daniel},
     TITLE = {Pseudo-isotopy and diffeomorphisms of the 4-sphere {I}: loops of spheres},
      YEAR = {2025},
    EPRINT = {2505.12088},
ARCHIVEPREFIX = {arXiv},
      NOTE = {arXiv:2505.12088},
}

@article {gay2025diffeomorphisms,
    AUTHOR = {Gay, David T.},
     TITLE = {Diffeomorphisms of the 4-sphere, {C}erf theory and
              {M}ontesinos twins},
   JOURNAL = {Algebr. Geom. Topol.},
  FJOURNAL = {Algebraic \& Geometric Topology},
    VOLUME = {25},
      YEAR = {2025},
    NUMBER = {5},
     PAGES = {2817--2849},
      ISSN = {1472-2747,1472-2739},
   MRCLASS = {57K40 (57R40 57R65)},
  MRNUMBER = {4996240},
MRREVIEWER = {Inasa\ Nakamura},
       DOI = {10.2140/agt.2025.25.2817},
       URL = {https://doi.org/10.2140/agt.2025.25.2817},
}

@article {gay2025relations,
    AUTHOR = {Gay, David T. and Hartman, Daniel},
     TITLE = {Relations amongst twists along {M}ontesinos twins in the
              4-sphere},
   JOURNAL = {Algebr. Geom. Topol.},
  FJOURNAL = {Algebraic \& Geometric Topology},
    VOLUME = {25},
      YEAR = {2025},
    NUMBER = {1},
     PAGES = {287--299},
      ISSN = {1472-2747,1472-2739},
   MRCLASS = {57K40 (57K45 57N37)},
  MRNUMBER = {4877258},
MRREVIEWER = {William\ Rushworth},
       DOI = {10.2140/agt.2025.25.287},
       URL = {https://doi.org/10.2140/agt.2025.25.287},
}

@misc{Gay2026PCMI,
    Author = {Gay, David T.},
    Title = {Surfaces in and diffeomorphisms of 4-manifolds},
    Year = {2026},
    Note = {{PCMI} minicourse lecture notes, Lecture 2},
    URL = {https://www.ias.edu/sites/default/files/PCMIMiniCourseLecture2.pdf},
}

@article {GGHKP2026,
    AUTHOR = {Gabai, David and Gay, David T. and Hartman, Daniel and Krushkal, Vyacheslav and Powell, Mark},
     TITLE = {Pseudo-isotopies of simply connected 4-manifolds},
   JOURNAL = {Forum Math. Pi},
  FJOURNAL = {Forum of Mathematics. Pi},
    VOLUME = {14},
      YEAR = {2026},
     PAGES = {Paper No. e9, 25},
      ISSN = {2050-5086},
   MRCLASS = {57K40 (57N37)},
  MRNUMBER = {5038934},
MRREVIEWER = {Marja\ K.\ Kankaanrinta},
       DOI = {10.1017/fmp.2026.10026},
       URL = {https://doi.org/10.1017/fmp.2026.10026},
}

@article {gompf2010more,
    AUTHOR = {Gompf, Robert E.},
     TITLE = {More {C}appell-{S}haneson spheres are standard},
   JOURNAL = {Algebr. Geom. Topol.},
  FJOURNAL = {Algebraic \& Geometric Topology},
    VOLUME = {10},
      YEAR = {2010},
    NUMBER = {3},
     PAGES = {1665--1681},
      ISSN = {1472-2747,1472-2739},
   MRCLASS = {57R60},
  MRNUMBER = {2683748},
MRREVIEWER = {Karl\ Heinz\ Dovermann},
       DOI = {10.2140/agt.2010.10.1665},
       URL = {https://doi.org/10.2140/agt.2010.10.1665},
}

@article {Gompf:1991,
    AUTHOR = {Gompf, Robert E.},
     TITLE = {Killing the {A}kbulut-{K}irby {$4$}-sphere, with relevance to
              the {A}ndrews-{C}urtis and {S}choenflies problems},
   JOURNAL = {Topology},
  FJOURNAL = {Topology. An International Journal of Mathematics},
    VOLUME = {30},
      YEAR = {1991},
    NUMBER = {1},
     PAGES = {97--115},
      ISSN = {0040-9383},
   MRCLASS = {57Q25 (57N13)},
  MRNUMBER = {1081936},
MRREVIEWER = {Martin\ Scharlemann},
       DOI = {10.1016/0040-9383(91)90036-4},
       URL = {https://doi.org/10.1016/0040-9383(91)90036-4},
}

@article {gromov1983positive,
    AUTHOR = {Gromov, Mikhael and Lawson, Jr., Herbert B.},
     TITLE = {Positive scalar curvature and the {D}irac operator on complete
              {R}iemannian manifolds},
   JOURNAL = {Inst. Hautes \'Etudes Sci. Publ. Math.},
  FJOURNAL = {Institut des Hautes \'Etudes Scientifiques. Publications
              Math\'ematiques},
    NUMBER = {58},
      YEAR = {1983},
     PAGES = {83--196},
      ISSN = {0073-8301,1618-1913},
   MRCLASS = {58G10 (53C20 57N10 57R99)},
  MRNUMBER = {720933},
MRREVIEWER = {Howard\ D.\ Rees},
       URL = {http://www.numdam.org/item?id=PMIHES_1983__58__83_0},
}

@article {Hillman-Wilson:1989,
    AUTHOR = {Hillman, Jonathan A. and Wilson, Stephen M. J.},
     TITLE = {On the reflexivity of {C}appell-{S}haneson {$2$}-knots},
   JOURNAL = {Bull. London Math. Soc.},
  FJOURNAL = {The Bulletin of the London Mathematical Society},
    VOLUME = {21},
      YEAR = {1989},
    NUMBER = {6},
     PAGES = {591--593},
      ISSN = {0024-6093,1469-2120},
   MRCLASS = {57Q45},
  MRNUMBER = {1018209},
MRREVIEWER = {Mitsuyoshi\ Kato},
       DOI = {10.1112/blms/21.6.591},
       URL = {https://doi.org/10.1112/blms/21.6.591},
}

@book {hirsch1976differential,
    AUTHOR = {Hirsch, Morris W.},
     TITLE = {Differential topology},
    SERIES = {Graduate Texts in Mathematics},
    VOLUME = {No. 33},
 PUBLISHER = {Springer-Verlag, New York-Heidelberg},
      YEAR = {1976},
     PAGES = {x+221},
   MRCLASS = {57DXX (58-01)},
  MRNUMBER = {448362},
MRREVIEWER = {Ulrich\ Koschorke},
}

@article {hitchin1974harmonic,
    AUTHOR = {Hitchin, Nigel J.},
     TITLE = {Harmonic spinors},
   JOURNAL = {Advances in Math.},
  FJOURNAL = {Advances in Mathematics},
    VOLUME = {14},
      YEAR = {1974},
     PAGES = {1--55},
      ISSN = {0001-8708},
   MRCLASS = {58G10 (14C30 32J25)},
  MRNUMBER = {358873},
MRREVIEWER = {F.\ Hirzebruch},
       DOI = {10.1016/0001-8708(74)90021-8},
       URL = {https://doi.org/10.1016/0001-8708(74)90021-8},
}

@article {hopkins2022intersection,
    AUTHOR = {Hopkins, Michael J. and Lin, Jianfeng and Shi, XiaoLin D. and Xu, Zhouli},
     TITLE = {Intersection forms of spin 4-manifolds and the {${\rm
              Pin}(2)$}-equivariant {M}ahowald invariant},
   JOURNAL = {Commun. Am. Math. Soc.},
  FJOURNAL = {Communications of the American Mathematical Society},
    VOLUME = {2},
      YEAR = {2022},
     PAGES = {22--132},
      ISSN = {2692-3688},
   MRCLASS = {55P91 (57K40 57K41)},
  MRNUMBER = {4385297},
MRREVIEWER = {J.\ P. C. Greenlees},
       DOI = {10.1090/cams/4},
       URL = {https://doi.org/10.1090/cams/4},
}

@article {Iida-Konno-Mukherjee-Taniguchi:2025,
    AUTHOR = {Iida, Nobuo and Konno, Hokuto and Mukherjee, Anubhav and
              Taniguchi, Masaki},
     TITLE = {Diffeomorphisms of 4-manifolds with boundary and exotic
              embeddings},
   JOURNAL = {Math. Ann.},
  FJOURNAL = {Mathematische Annalen},
    VOLUME = {391},
      YEAR = {2025},
    NUMBER = {2},
     PAGES = {1845--1897},
      ISSN = {0025-5831,1432-1807},
   MRCLASS = {57R50 (57K41)},
  MRNUMBER = {4853007},
MRREVIEWER = {Sergey\ M.\ Finashin},
       DOI = {10.1007/s00208-024-02974-x},
       URL = {https://doi.org/10.1007/s00208-024-02974-x},
}

@article {Iwaki:2025,
    AUTHOR = {Iwaki, Kazunori},
     TITLE = {Infinite families of standard {C}appell-{S}haneson homotopy
              4-spheres},
   JOURNAL = {Topology Appl.},
  FJOURNAL = {Topology and its Applications},
    VOLUME = {366},
      YEAR = {2025},
     PAGES = {Paper No. 109293, 17},
      ISSN = {0166-8641,1879-3207},
   MRCLASS = {57K40 (55P10)},
  MRNUMBER = {4873792},
       DOI = {10.1016/j.topol.2025.109293},
       URL = {https://doi.org/10.1016/j.topol.2025.109293},
}

@misc{kang2025exotic,
    AUTHOR = {Kang, Sungkyung and Park, JungHwan and Taniguchi, Masaki},
     TITLE = {Equivariant lattice homology for {Seifert} homology 3-spheres},
      YEAR = {2025},
    EPRINT = {2510.12394},
ARCHIVEPREFIX = {arXiv},
      NOTE = {arXiv:2510.12394},
}

@article {Kato:2022,
    AUTHOR = {Kato, Yuya},
     TITLE = {Nonsmoothable actions of {$\Bbb Z_2 \times \Bbb Z_2$} on spin
              four-manifolds},
   JOURNAL = {Topology Appl.},
  FJOURNAL = {Topology and its Applications},
    VOLUME = {307},
      YEAR = {2022},
     PAGES = {Paper No. 107868, 13},
      ISSN = {0166-8641,1879-3207},
   MRCLASS = {57M60},
  MRNUMBER = {4365044},
       DOI = {10.1016/j.topol.2021.107868},
       URL = {https://doi.org/10.1016/j.topol.2021.107868},
}

@article {Kim-Yamada:2023,
    AUTHOR = {Kim, Min Hoon and Yamada, Shohei},
     TITLE = {Ideal classes and {C}appell-{S}haneson homotopy 4-spheres},
   JOURNAL = {Kyungpook Math. J.},
  FJOURNAL = {Kyungpook Mathematical Journal},
    VOLUME = {63},
      YEAR = {2023},
    NUMBER = {3},
     PAGES = {373--411},
      ISSN = {1225-6951,0454-8124},
   MRCLASS = {58D27 (58A05)},
  MRNUMBER = {4655393},
MRREVIEWER = {Krzysztof\ M.\ Pawa\l owski},
       DOI = {10.5666/KMJ.2023.63.3.373},
       URL = {https://doi.org/10.5666/KMJ.2023.63.3.373},
}

@misc{KLMM:2024-Milnor,
    AUTHOR = {Konno, Hokuto and Lin, Jianfeng and Mukherjee, Anubhav and Mu{\~n}oz-Ech{\'a}niz, Juan},
     TITLE = {On four-dimensional {Dehn} twists and {Milnor} fibrations},
      YEAR = {2024},
    EPRINT = {2409.11961},
ARCHIVEPREFIX = {arXiv},
      NOTE = {arXiv:2409.11961, to appear in Duke Math. J.},
}

@misc{KLMM:2024-monodromy,
    AUTHOR = {Konno, Hokuto and Lin, Jianfeng and Mukherjee, Anubhav and Mu{\~n}oz-Ech{\'a}niz, Juan},
     TITLE = {The monodromy diffeomorphism of weighted singularities and {Seiberg--Witten} theory},
      YEAR = {2024},
    EPRINT = {2411.12202},
ARCHIVEPREFIX = {arXiv},
      NOTE = {arXiv:2411.12202},
}

@misc{KLMM:2025-Lefschetz,
    AUTHOR = {Konno, Hokuto and Lin, Jianfeng and Mukherjee, Anubhav and Mu{\~n}oz-Ech{\'a}niz, Juan},
     TITLE = {Constraints on {Lefschetz} fibrations with four-dimensional fibers from {Seiberg--Witten} theory},
      YEAR = {2025},
    EPRINT = {2511.00601},
ARCHIVEPREFIX = {arXiv},
      NOTE = {arXiv:2511.00601},
}

@misc{KLMM:2026-aspherical,
    AUTHOR = {Konno, Hokuto and Lin, Jianfeng and Mukherjee, Anubhav and Mu{\~n}oz-Ech{\'a}niz, Juan},
     TITLE = {Exotic diffeomorphisms on closed aspherical 4-manifolds},
      YEAR = {2026},
    EPRINT = {2610.21083},
ARCHIVEPREFIX = {arXiv},
      NOTE = {arXiv:2610.21083},
}

@article {Konno-Lin:2026,
    AUTHOR = {Konno, Hokuto and Lin, Jianfeng},
     TITLE = {Homological instability for moduli spaces of smooth
              4-manifolds},
   JOURNAL = {Geom. Topol.},
  FJOURNAL = {Geometry \& Topology},
    VOLUME = {30},
      YEAR = {2026},
    NUMBER = {8},
     PAGES = {2993--3068},
      ISSN = {1465-3060,1364-0380},
   MRCLASS = {57K41 (57S05)},
  MRNUMBER = {5124776},
       DOI = {10.2140/gt.2026.30.2993},
       URL = {https://doi.org/10.2140/gt.2026.30.2993},
}

@misc{Konno-Mallick-Taniguchi:2025-absolute,
    AUTHOR = {Konno, Hokuto and Mallick, Abhishek and Taniguchi, Masaki},
     TITLE = {On absolutely exotic diffeomorphisms of 4-manifolds},
      YEAR = {2025},
    EPRINT = {2510.05505},
ARCHIVEPREFIX = {arXiv},
      NOTE = {arXiv:2510.05505},
}

@article {Konno-Mallick-Taniguchi:2026,
    AUTHOR = {Konno, Hokuto and Mallick, Abhishek and Taniguchi, Masaki},
     TITLE = {Exotic {D}ehn twists on 4-manifolds},
   JOURNAL = {Geom. Topol.},
  FJOURNAL = {Geometry \& Topology},
    VOLUME = {30},
      YEAR = {2026},
    NUMBER = {6},
     PAGES = {2395--2429},
      ISSN = {1465-3060,1364-0380},
   MRCLASS = {57K41 (57R50 57R58 58D05)},
  MRNUMBER = {5111601},
       DOI = {10.2140/gt.2026.30.2395},
       URL = {https://doi.org/10.2140/gt.2026.30.2395},
}

@article {Konno-Miyazawa-Taniguchi:2023,
    AUTHOR = {Konno, Hokuto and Miyazawa, Jin and Taniguchi, Masaki},
     TITLE = {Involutions, links, and {F}loer cohomologies},
   JOURNAL = {J. Topol.},
  FJOURNAL = {Journal of Topology},
    VOLUME = {17},
      YEAR = {2024},
    NUMBER = {2},
     PAGES = {Paper No. e12340, 47},
      ISSN = {1753-8416,1753-8424},
   MRCLASS = {57K10 (57K41 57M60 57R58)},
  MRNUMBER = {4821360},
MRREVIEWER = {Nikolai\ N.\ Saveliev},
       DOI = {10.1112/topo.12340},
       URL = {https://doi.org/10.1112/topo.12340},
}

@article {Konno-Miyazawa-Taniguchi:2025,
    AUTHOR = {Konno, Hokuto and Miyazawa, Jin and Taniguchi, Masaki},
     TITLE = {Involutions, knots, and {F}loer {$K$}-theory},
   JOURNAL = {Compos. Math.},
  FJOURNAL = {Compositio Mathematica},
    VOLUME = {161},
      YEAR = {2025},
    NUMBER = {11},
     PAGES = {2852--2910},
      ISSN = {0010-437X,1570-5846},
   MRCLASS = {57K10 (57R58)},
  MRNUMBER = {5002153},
MRREVIEWER = {Mattia\ Mecchia},
       DOI = {10.1112/S0010437X25102820},
       URL = {https://doi.org/10.1112/S0010437X25102820},
}

@article {Konno-Taniguchi:2020-PSC,
    AUTHOR = {Konno, Hokuto and Taniguchi, Masaki},
     TITLE = {Positive scalar curvature and 10/8-type inequalities on
              4-manifolds with periodic ends},
   JOURNAL = {Invent. Math.},
  FJOURNAL = {Inventiones Mathematicae},
    VOLUME = {222},
      YEAR = {2020},
    NUMBER = {3},
     PAGES = {833--880},
      ISSN = {0020-9910,1432-1297},
   MRCLASS = {53C21 (53C20 57K41)},
  MRNUMBER = {4169052},
MRREVIEWER = {Thomas\ Schick},
       DOI = {10.1007/s00222-020-00979-2},
       URL = {https://doi.org/10.1007/s00222-020-00979-2},
}

@article {Konno:2024-moduli,
    AUTHOR = {Konno, Hokuto},
     TITLE = {The homology of moduli spaces of 4-manifolds may be infinitely
              generated},
   JOURNAL = {Forum Math. Pi},
  FJOURNAL = {Forum of Mathematics. Pi},
    VOLUME = {12},
      YEAR = {2024},
     PAGES = {Paper No. e25, 18},
      ISSN = {2050-5086},
   MRCLASS = {57S05 (57K41 57R50 57R52)},
  MRNUMBER = {4831140},
MRREVIEWER = {Rainer\ L\"owen},
       DOI = {10.1017/fmp.2024.26},
       URL = {https://doi.org/10.1017/fmp.2024.26},
}

@article {KPT-cables,
    AUTHOR = {Kang, Sungkyung and Park, JungHwan and Taniguchi, Masaki},
     TITLE = {Cables of the figure-eight knot via real {F}r\o yshov
              invariants},
   JOURNAL = {Geom. Topol.},
  FJOURNAL = {Geometry \& Topology},
    VOLUME = {30},
      YEAR = {2026},
    NUMBER = {2},
     PAGES = {781--833},
      ISSN = {1465-3060,1364-0380},
   MRCLASS = {57K10 (57K41)},
  MRNUMBER = {5047786},
MRREVIEWER = {Allan\ Edmonds},
       DOI = {10.2140/gt.2026.30.781},
       URL = {https://doi.org/10.2140/gt.2026.30.781},
}

@article {KPT26,
    AUTHOR = {Kang, Sungkyung and Park, JungHwan and Taniguchi, Masaki},
     TITLE = {Exotic {D}ehn twists and homotopy coherent group actions},
   JOURNAL = {Invent. Math.},
  FJOURNAL = {Inventiones Mathematicae},
    VOLUME = {243},
      YEAR = {2026},
    NUMBER = {1},
     PAGES = {209--241},
      ISSN = {0020-9910,1432-1297},
   MRCLASS = {57K41 (57R50 57R85)},
  MRNUMBER = {5008150},
MRREVIEWER = {Marja\ K.\ Kankaanrinta},
       DOI = {10.1007/s00222-025-01378-1},
       URL = {https://doi.org/10.1007/s00222-025-01378-1},
}

@incollection {Kreck1979,
    AUTHOR = {Kreck, Matthias},
     TITLE = {Isotopy classes of diffeomorphisms of {$(k-1)$}-connected
              almost-parallelizable {$2k$}-manifolds},
 BOOKTITLE = {Algebraic topology, {A}arhus 1978 ({P}roc. {S}ympos., {U}niv.
              {A}arhus, {A}arhus, 1978)},
    SERIES = {Lecture Notes in Math.},
    VOLUME = {763},
     PAGES = {643--663},
 PUBLISHER = {Springer, Berlin},
      YEAR = {1979},
      ISBN = {3-540-09721-X},
   MRCLASS = {57R50},
  MRNUMBER = {561244},
MRREVIEWER = {E.\ C.\ Turner},
}

@article {Kronheimer-Mrowka:2020-1,
    AUTHOR = {Kronheimer, Peter B. and Mrowka, Tomasz S.},
     TITLE = {The {D}ehn twist on a sum of two {$K3$} surfaces},
   JOURNAL = {Math. Res. Lett.},
  FJOURNAL = {Mathematical Research Letters},
    VOLUME = {27},
      YEAR = {2020},
    NUMBER = {6},
     PAGES = {1767--1783},
      ISSN = {1073-2780,1945-001X},
   MRCLASS = {57R50 (57K40)},
  MRNUMBER = {4216604},
       DOI = {10.4310/MRL.2020.v27.n6.a8},
       URL = {https://doi.org/10.4310/MRL.2020.v27.n6.a8},
}

@article {Krushkal-Mukherjee-Powell-Warren:2026,
    AUTHOR = {Krushkal, Vyacheslav and Mukherjee, Anubhav and Powell, Mark
              and Warren, Terrin},
     TITLE = {Corks for exotic diffeomorphisms},
   JOURNAL = {Geom. Topol.},
  FJOURNAL = {Geometry \& Topology},
    VOLUME = {30},
      YEAR = {2026},
    NUMBER = {6},
     PAGES = {2297--2335},
      ISSN = {1465-3060,1364-0380},
   MRCLASS = {57K40 (57N37 57R58)},
  MRNUMBER = {5111598},
       DOI = {10.2140/gt.2026.30.2297},
       URL = {https://doi.org/10.2140/gt.2026.30.2297},
}

@article {lee1995homology,
    AUTHOR = {Lee, Ronnie and Weintraub, Steven H.},
     TITLE = {On the homology of double branched covers},
   JOURNAL = {Proc. Amer. Math. Soc.},
  FJOURNAL = {Proceedings of the American Mathematical Society},
    VOLUME = {123},
      YEAR = {1995},
    NUMBER = {4},
     PAGES = {1263--1266},
      ISSN = {0002-9939,1088-6826},
   MRCLASS = {57M12},
  MRNUMBER = {1224618},
MRREVIEWER = {Neal\ Brand},
       DOI = {10.2307/2160729},
       URL = {https://doi.org/10.2307/2160729},
}

@article {lewis1995change,
    AUTHOR = {Lewis, Jr., LeMoine G.},
     TITLE = {Change of universe functors in equivariant stable homotopy
              theory},
   JOURNAL = {Fund. Math.},
  FJOURNAL = {Fundamenta Mathematicae},
    VOLUME = {148},
      YEAR = {1995},
    NUMBER = {2},
     PAGES = {117--158},
      ISSN = {0016-2736,1730-6329},
   MRCLASS = {55P91 (55P20 55Q91)},
  MRNUMBER = {1360142},
MRREVIEWER = {J.\ M.\ Boardman},
       DOI = {10.4064/fm-148-2-117-158},
       URL = {https://doi.org/10.4064/fm-148-2-117-158},
}

@article {lewis2000splitting,
    AUTHOR = {Lewis, Jr., LeMoine G.},
     TITLE = {Splitting theorems for certain equivariant spectra},
   JOURNAL = {Mem. Amer. Math. Soc.},
  FJOURNAL = {Memoirs of the American Mathematical Society},
    VOLUME = {144},
      YEAR = {2000},
    NUMBER = {686},
     PAGES = {x+89},
      ISSN = {0065-9266,1947-6221},
   MRCLASS = {55P42 (55M35 55N20 55P91 55R12)},
  MRNUMBER = {1679450},
MRREVIEWER = {Steven\ R.\ Costenoble},
       DOI = {10.1090/memo/0686},
       URL = {https://doi.org/10.1090/memo/0686},
}

@book {LewisMaySteinberger,
    AUTHOR = {Lewis, Jr., LeMoine G. and May, Jon P. and Steinberger, Mark and McClure, James E.},
     TITLE = {Equivariant stable homotopy theory},
    SERIES = {Lecture Notes in Mathematics},
    VOLUME = {1213},
      NOTE = {With contributions by J. E. McClure},
 PUBLISHER = {Springer-Verlag, Berlin},
      YEAR = {1986},
     PAGES = {x+538},
      ISBN = {3-540-16820-6},
   MRCLASS = {55-02 (55Nxx 55Pxx 57S99)},
  MRNUMBER = {866482},
MRREVIEWER = {T.\ tom Dieck},
       DOI = {10.1007/BFb0075778},
       URL = {https://doi.org/10.1007/BFb0075778},
}

@article {lin2025family,
    AUTHOR = {Lin, Jianfeng and Mukherjee, Anubhav},
     TITLE = {Family {B}auer-{F}uruta invariant, exotic surfaces and {S}male
              conjecture},
   JOURNAL = {J. Assoc. Math. Res.},
  FJOURNAL = {Journal of the Association for Mathematical Research},
    VOLUME = {3},
      YEAR = {2025},
    NUMBER = {2},
     PAGES = {237--275},
      ISSN = {2998-4114},
   MRCLASS = {57R57 (55Q91 57K41 57R50)},
  MRNUMBER = {4958155},
       DOI = {10.56994/jamr.003.002.003},
       URL = {https://doi.org/10.56994/jamr.003.002.003},
}

@article {Lin:2023-1,
    AUTHOR = {Lin, Jianfeng},
     TITLE = {Isotopy of the {D}ehn twist on {$K3\, \#\, K3$} after a single
              stabilization},
   JOURNAL = {Geom. Topol.},
  FJOURNAL = {Geometry \& Topology},
    VOLUME = {27},
      YEAR = {2023},
    NUMBER = {5},
     PAGES = {1987--2012},
      ISSN = {1465-3060,1364-0380},
   MRCLASS = {57R52 (55P91 57R50 57R57)},
  MRNUMBER = {4621924},
MRREVIEWER = {Andrey\ Gogolev},
       DOI = {10.2140/gt.2023.27.1987},
       URL = {https://doi.org/10.2140/gt.2023.27.1987},
}

@article {Ma16,
    AUTHOR = {Manolescu, Ciprian},
     TITLE = {Pin(2)-equivariant {S}eiberg-{W}itten {F}loer homology and the
              triangulation conjecture},
   JOURNAL = {J. Amer. Math. Soc.},
  FJOURNAL = {Journal of the American Mathematical Society},
    VOLUME = {29},
      YEAR = {2016},
    NUMBER = {1},
     PAGES = {147--176},
      ISSN = {0894-0347,1088-6834},
   MRCLASS = {57R58 (57M27 57Q15)},
  MRNUMBER = {3402697},
MRREVIEWER = {Tirasan\ Khandhawit},
       DOI = {10.1090/jams829},
       URL = {https://doi.org/10.1090/jams829},
}

@article {manolescu2007gluing,
    AUTHOR = {Manolescu, Ciprian},
     TITLE = {A gluing theorem for the relative {B}auer-{F}uruta invariants},
   JOURNAL = {J. Differential Geom.},
  FJOURNAL = {Journal of Differential Geometry},
    VOLUME = {76},
      YEAR = {2007},
    NUMBER = {1},
     PAGES = {117--153},
      ISSN = {0022-040X,1945-743X},
   MRCLASS = {57R57 (57N13 57R56 57R58)},
  MRNUMBER = {2312050},
MRREVIEWER = {Rostislav\ Matveyev},
       URL = {http://projecteuclid.org/euclid.jdg/1180135667},
}

@article {marques2012deforming,
    AUTHOR = {Cod{\'a} Marques, Fernando},
     TITLE = {Deforming three-manifolds with positive scalar curvature},
   JOURNAL = {Ann. of Math. (2)},
  FJOURNAL = {Annals of Mathematics. Second Series},
    VOLUME = {176},
      YEAR = {2012},
    NUMBER = {2},
     PAGES = {815--863},
      ISSN = {0003-486X,1939-8980},
   MRCLASS = {53C20 (53C44 57M50)},
  MRNUMBER = {2950765},
MRREVIEWER = {David\ J.\ Wraith},
       DOI = {10.4007/annals.2012.176.2.3},
       URL = {https://doi.org/10.4007/annals.2012.176.2.3},
}

@book {may1996equivariant,
    AUTHOR = {May, Jon P.},
     TITLE = {Equivariant homotopy and cohomology theory},
    SERIES = {CBMS Regional Conference Series in Mathematics},
    VOLUME = {91},
      NOTE = {With contributions by M. Cole, G. Comeza\~na, S. Costenoble,
              A. D. Elmendorf, J. P. C. Greenlees, L. G. Lewis, Jr., R. J.
              Piacenza, G. Triantafillou, and S. Waner},
 PUBLISHER = {Conference Board of the Mathematical Sciences, Washington, DC;
              by the American Mathematical Society, Providence, RI},
      YEAR = {1996},
     PAGES = {xiv+366},
      ISBN = {0-8218-0319-0},
   MRCLASS = {55P91 (18G99 55N91 55U35)},
  MRNUMBER = {1413302},
MRREVIEWER = {Ian\ Hambleton},
       DOI = {10.1090/cbms/091},
       URL = {https://doi.org/10.1090/cbms/091},
}

@article {may2003wirthmuller,
    AUTHOR = {May, Jon P.},
     TITLE = {The {W}irthm\"uller isomorphism revisited},
   JOURNAL = {Theory Appl. Categ.},
  FJOURNAL = {Theory and Applications of Categories},
    VOLUME = {11},
      YEAR = {2003},
     PAGES = {No. 5, 132--142},
      ISSN = {1201-561X},
   MRCLASS = {55P42 (55P91)},
  MRNUMBER = {1988073},
       DOI = {10.70930/tac/o56xigni},
       URL = {https://doi.org/10.70930/tac/o56xigni},
}

@book {may2006parametrized,
    AUTHOR = {May, Jon P. and Sigurdsson, Johann},
     TITLE = {Parametrized homotopy theory},
    SERIES = {Mathematical Surveys and Monographs},
    VOLUME = {132},
 PUBLISHER = {American Mathematical Society, Providence, RI},
      YEAR = {2006},
     PAGES = {x+441},
      ISBN = {978-0-8218-3922-5; 0-8218-3922-5},
   MRCLASS = {55P42 (19L99 55N20 55N22 55P91)},
  MRNUMBER = {2271789},
MRREVIEWER = {A.\ A.\ Ranicki},
       DOI = {10.1090/surv/132},
       URL = {https://doi.org/10.1090/surv/132},
}

@article {Meier-Zupan:2022,
    AUTHOR = {Meier, Jeffrey and Zupan, Alexander},
     TITLE = {Generalized square knots and homotopy 4-spheres},
   JOURNAL = {J. Differential Geom.},
  FJOURNAL = {Journal of Differential Geometry},
    VOLUME = {122},
      YEAR = {2022},
    NUMBER = {1},
     PAGES = {69--129},
      ISSN = {0022-040X,1945-743X},
   MRCLASS = {57K10},
  MRNUMBER = {4507471},
MRREVIEWER = {Masakazu\ Teragaito},
       DOI = {10.4310/jdg/1668186788},
       URL = {https://doi.org/10.4310/jdg/1668186788},
}

@misc{miyazawa2023gauge,
    AUTHOR = {Miyazawa, Jin},
     TITLE = {A gauge theoretic invariant of embedded surfaces in 4-manifolds and exotic {$P^2$}-knots},
      YEAR = {2023},
    EPRINT = {2312.02041},
ARCHIVEPREFIX = {arXiv},
      NOTE = {arXiv:2312.02041, to appear in Ann. of Math. (2)},
}

@misc{miyazawa2025satellite,
    AUTHOR = {Miyazawa, Jin and Park, JungHwan and Taniguchi, Masaki},
     TITLE = {A satellite formula for real {Seiberg--Witten Floer} homotopy types},
      YEAR = {2025},
    EPRINT = {2504.03270},
ARCHIVEPREFIX = {arXiv},
      NOTE = {arXiv:2504.03270},
}

@misc{Miyazawa:2024-boundary,
    AUTHOR = {Miyazawa, Jin},
     TITLE = {Boundary {Dehn} twists on {Milnor} fibers and {Family} {Bauer--Furuta} invariants},
      YEAR = {2024},
    EPRINT = {2410.21742},
ARCHIVEPREFIX = {arXiv},
      NOTE = {arXiv:2410.21742},
}

@article {MorganMrowkaSzabo1997,
    AUTHOR = {Morgan, John W. and Mrowka, Tomasz S. and Szab{\'o}, Zolt{\'a}n},
     TITLE = {Product formulas along {$T^3$} for {S}eiberg-{W}itten
              invariants},
   JOURNAL = {Math. Res. Lett.},
  FJOURNAL = {Mathematical Research Letters},
    VOLUME = {4},
      YEAR = {1997},
    NUMBER = {6},
     PAGES = {915--929},
      ISSN = {1073-2780},
   MRCLASS = {57R57 (57N13)},
  MRNUMBER = {1492130},
MRREVIEWER = {Dave\ Auckly},
       DOI = {10.4310/MRL.1997.v4.n6.a11},
       URL = {https://doi.org/10.4310/MRL.1997.v4.n6.a11},
}

@article {nakamura2015mathrm,
    AUTHOR = {Nakamura, Nobuhiro},
     TITLE = {{$\mathrm{Pin}^-(2)$}-monopole invariants},
   JOURNAL = {J. Differential Geom.},
  FJOURNAL = {Journal of Differential Geometry},
    VOLUME = {101},
      YEAR = {2015},
    NUMBER = {3},
     PAGES = {507--549},
      ISSN = {0022-040X,1945-743X},
   MRCLASS = {57N13 (57R57)},
  MRNUMBER = {3415770},
MRREVIEWER = {Raphael\ Zentner},
       URL = {http://projecteuclid.org/euclid.jdg/1445518922},
}

@article {Nakamura:2013,
    AUTHOR = {Nakamura, Nobuhiro},
     TITLE = {{$\rm{Pin}^-(2)$}-monopole equations and intersection forms
              with local coefficients of four-manifolds},
   JOURNAL = {Math. Ann.},
  FJOURNAL = {Mathematische Annalen},
    VOLUME = {357},
      YEAR = {2013},
    NUMBER = {3},
     PAGES = {915--939},
      ISSN = {0025-5831,1432-1807},
   MRCLASS = {57R57},
  MRNUMBER = {3118618},
MRREVIEWER = {Thomas\ Gibbs\ Leness},
       DOI = {10.1007/s00208-013-0924-3},
       URL = {https://doi.org/10.1007/s00208-013-0924-3},
}

@misc{Per17,
    AUTHOR = {Perlmutter, Nathan},
     TITLE = {Parametrized {Morse} theory and positive scalar curvature},
      YEAR = {2017},
    EPRINT = {1705.02754},
ARCHIVEPREFIX = {arXiv},
      NOTE = {arXiv:1705.02754},
}

@misc{Qiu:2024-exotic,
    AUTHOR = {Qiu, Haochen},
     TITLE = {Exotic diffeomorphisms on 4-manifolds with {$b_2^+=2$}},
      YEAR = {2024},
    EPRINT = {2409.07009},
ARCHIVEPREFIX = {arXiv},
      NOTE = {arXiv:2409.07009},
}

@article {Quinn1986,
    AUTHOR = {Quinn, Frank},
     TITLE = {Isotopy of {$4$}-manifolds},
   JOURNAL = {J. Differential Geom.},
  FJOURNAL = {Journal of Differential Geometry},
    VOLUME = {24},
      YEAR = {1986},
    NUMBER = {3},
     PAGES = {343--372},
      ISSN = {0022-040X,1945-743X},
   MRCLASS = {57N13 (57N37 57S05)},
  MRNUMBER = {868975},
MRREVIEWER = {Laurence\ R.\ Taylor},
       URL = {http://projecteuclid.org/euclid.jdg/1214440552},
}

@article {ruberman2002positive,
    AUTHOR = {Ruberman, Daniel},
     TITLE = {Positive scalar curvature, diffeomorphisms and the
              {S}eiberg-{W}itten invariants},
   JOURNAL = {Geom. Topol.},
  FJOURNAL = {Geometry and Topology},
    VOLUME = {5},
      YEAR = {2001},
     PAGES = {895--924},
      ISSN = {1465-3060,1364-0380},
   MRCLASS = {57R57 (53C21)},
  MRNUMBER = {1874146},
MRREVIEWER = {Stefano\ Vidussi},
       DOI = {10.2140/gt.2001.5.895},
       URL = {https://doi.org/10.2140/gt.2001.5.895},
}

@article {Ruberman:1998,
    AUTHOR = {Ruberman, Daniel},
     TITLE = {An obstruction to smooth isotopy in dimension {$4$}},
   JOURNAL = {Math. Res. Lett.},
  FJOURNAL = {Mathematical Research Letters},
    VOLUME = {5},
      YEAR = {1998},
    NUMBER = {6},
     PAGES = {743--758},
      ISSN = {1073-2780},
   MRCLASS = {57R52 (57N37 57R57)},
  MRNUMBER = {1671187},
MRREVIEWER = {Ignasi\ Mundet-Riera},
       DOI = {10.4310/MRL.1998.v5.n6.a5},
       URL = {https://doi.org/10.4310/MRL.1998.v5.n6.a5},
}

@inproceedings {Ruberman:1999,
    AUTHOR = {Ruberman, Daniel},
     TITLE = {A polynomial invariant of diffeomorphisms of 4-manifolds},
 BOOKTITLE = {Proceedings of the {K}irbyfest ({B}erkeley, {CA}, 1998)},
    SERIES = {Geom. Topol. Monogr.},
    VOLUME = {2},
     PAGES = {473--488},
 PUBLISHER = {Geom. Topol. Publ., Coventry},
      YEAR = {1999},
   MRCLASS = {57R57 (57N13 57R50)},
  MRNUMBER = {1734421},
MRREVIEWER = {Ignasi\ Mundet-Riera},
       DOI = {10.2140/gtm.1999.2.473},
       URL = {https://doi.org/10.2140/gtm.1999.2.473},
}

@incollection {stolz2014metrics,
    AUTHOR = {Rosenberg, Jonathan M. and Stolz, Stephan},
     TITLE = {Metrics of positive scalar curvature and connections with
              surgery},
 BOOKTITLE = {Surveys on surgery theory, {V}ol. 2},
    SERIES = {Ann. of Math. Stud.},
    VOLUME = {149},
     PAGES = {353--386},
 PUBLISHER = {Princeton Univ. Press, Princeton, NJ},
      YEAR = {2001},
      ISBN = {0-691-08814-4; 0-691-08815-2},
   MRCLASS = {53C21 (53C27 57R65)},
  MRNUMBER = {1818778},
MRREVIEWER = {Thomas\ Schick},
}

@article {szymik2010characteristic,
    AUTHOR = {Szymik, Markus},
     TITLE = {Characteristic cohomotopy classes for families of 4-manifolds},
   JOURNAL = {Forum Math.},
  FJOURNAL = {Forum Mathematicum},
    VOLUME = {22},
      YEAR = {2010},
    NUMBER = {3},
     PAGES = {509--523},
      ISSN = {0933-7741,1435-5337},
   MRCLASS = {57R57 (55P42 55Q55 55R70 57R20)},
  MRNUMBER = {2652709},
MRREVIEWER = {Vicente\ Mu\~noz},
       DOI = {10.1515/FORUM.2010.027},
       URL = {https://doi.org/10.1515/FORUM.2010.027},
}

@article {Tian-Wang:2009,
    AUTHOR = {Tian, Gang and Wang, Shuguang},
     TITLE = {Orientability and real {S}eiberg-{W}itten invariants},
   JOURNAL = {Internat. J. Math.},
  FJOURNAL = {International Journal of Mathematics},
    VOLUME = {20},
      YEAR = {2009},
    NUMBER = {5},
     PAGES = {573--604},
      ISSN = {0129-167X,1793-6519},
   MRCLASS = {57R57 (14P25 53C27)},
  MRNUMBER = {2526308},
MRREVIEWER = {Hans\ U.\ Boden},
       DOI = {10.1142/S0129167X09005455},
       URL = {https://doi.org/10.1142/S0129167X09005455},
}

@article {tomDieck1975Orbittypen,
    AUTHOR = {tom Dieck, Tammo},
     TITLE = {Orbittypen und \"aquivariante {H}omologie. {II}},
   JOURNAL = {Arch. Math. (Basel)},
  FJOURNAL = {Archiv der Mathematik},
    VOLUME = {26},
      YEAR = {1975},
    NUMBER = {6},
     PAGES = {650--662},
      ISSN = {0003-889X,1420-8938},
   MRCLASS = {57F15},
  MRNUMBER = {436177},
MRREVIEWER = {H.\ Scheerer},
       DOI = {10.1007/BF01229795},
       URL = {https://doi.org/10.1007/BF01229795},
}

@article {Walsh2011,
    AUTHOR = {Walsh, Mark G.},
     TITLE = {Metrics of positive scalar curvature and generalised {M}orse
              functions, {P}art {I}},
   JOURNAL = {Mem. Amer. Math. Soc.},
  FJOURNAL = {Memoirs of the American Mathematical Society},
    VOLUME = {209},
      YEAR = {2011},
    NUMBER = {983},
     PAGES = {xviii+80},
      ISSN = {0065-9266,1947-6221},
      ISBN = {978-0-8218-5304-7},
   MRCLASS = {53C20 (53C21 57R65 58E05)},
  MRNUMBER = {2789750},
MRREVIEWER = {David\ J.\ Wraith},
       DOI = {10.1090/S0065-9266-10-00622-8},
       URL = {https://doi.org/10.1090/S0065-9266-10-00622-8},
}

@article {walsh2018aspects,
    AUTHOR = {Walsh, Mark G.},
     TITLE = {Aspects of positive scalar curvature and topology {II}},
   JOURNAL = {Irish Math. Soc. Bull.},
  FJOURNAL = {Irish Mathematical Society Bulletin},
    NUMBER = {81},
      YEAR = {2018},
     PAGES = {57--95},
      ISSN = {0791-5578},
   MRCLASS = {58D17 (53C20)},
  MRNUMBER = {3839561},
MRREVIEWER = {Gang\ Li},
}

\end{document}